# On the classification of unitary representations of reductive Lie groups

By Susana A. Salamanca-Riba and David A. Vogan, Jr.*

### Table of Contents

### Introduction

Suppose $G$ is a real reductive Lie group in Harish-Chandra's class (see [6, §3]). We propose here a structure for the set $\Pi_u(G)$ of equivalence classes of irreducible unitary representations of $G$. (The subscript $u$ will be used throughout to indicate structures related to unitary representations.) We decompose $\Pi_u(G)$ into disjoint subsets with a (very explicit) discrete parameter set $\Lambda_u$:

$$(0.1) \qquad \Pi_u(G) = \bigcup_{\lambda_u \in \Lambda_u} \Pi_u^{\lambda_u}(G).$$

Each subset is identified conjecturally (Conjecture 0.6) with a collection of unitary representations of a certain subgroup $G(\lambda_u)$ of $G$. (We will give strong evidence and partial results for this conjecture.) In this way the problem of classifying $\Pi_u(G)$ would be reduced (by induction on the dimension of $G$) to the case $G(\lambda_u) = G$.

Before considering the general program in more detail, we describe it in the familiar case $G = \mathrm{SL}(2, \mathbb{R})$. (This example will be treated more com-

*Research by Susana Salamanca-Riba is partially supported by NSF Grant #DMS-9510608. Research by David Vogan is partially supported by NSF Grant #DMS-9402994.

pletely in Example 6.2 below.) The parameter set $\Lambda_u$ is $\mathbb{Z}$. If $n$ is a nonzero integer, then $G(n) = \mathrm{SO}(2)$, and $\Pi_u^n(G)$ consists of a single discrete series representation (with Harish-Chandra parameter $n + \operatorname{sgn} n$ and lowest $\mathrm{SO}(2)$-type $n + 2\operatorname{sgn} n$). The corresponding representation of $\mathrm{SO}(2)$ is the character usually parametrized by $n$. If $n = 0$, then $G(0) = \mathrm{SL}(2, \mathbb{R})$, and $\Pi_u^0(G)$ consists of all the remaining unitary representations: those containing an $\mathrm{SO}(2)$ character parametrized by $0$, $\pm 1$, or $\pm 2$. (These are the principal series, the complementary series, the trivial representation, and the first two discrete series.)

At first glance it might appear natural to separate also the first two discrete series as isolated points in the parametrization. But recall that these representations are not isolated in the unitary dual of $\mathrm{SL}(2, \mathbb{R})$; rather they are limits of complementary series. The decomposition we have given is consistent with the topology of $\Pi_u(G)$, and separating the first two discrete series would not be.

As the example of $\mathrm{SL}(2, \mathbb{R})$ indicates, our "reduction" to the case $G(\lambda_u)$ $= G$ is in some ways not very powerful; almost all of the interesting phenomena are swept into that case. Nevertheless there is some content in these ideas. In a separate paper, the first author will show how to use them to deduce an old conjecture classifying unitary representations with sufficiently regular infinitesimal character.

We turn now to a description of the decomposition (??) in general. Choose a Cartan involution $\theta$ of $G$, with fixed points $K$ a maximal compact subgroup of $G$. The real Lie algebra $\mathfrak{g}_0$ then acquires a Cartan decomposition

$$(0.2a) \qquad \mathfrak{g}_0 = \mathfrak{k}_0 + \mathfrak{p}_0,$$

with $\mathfrak{p}_0$ the $-1$ eigenspace of (the differential of) $\theta$. The complexification of $\mathfrak{g}_0$ is written

$$(0.2b) \qquad \mathfrak{g} = \mathfrak{g}_0 \otimes_{\mathbb{R}} \mathbb{C};$$

similar notation is used for other groups. Fix a nondegenerate $G$-invariant, $\theta$-invariant symmetric bilinear form $\langle \ , \ \rangle$ on $\mathfrak{g}_0$, negative definite on $\mathfrak{k}_0$ and positive definite on $\mathfrak{p}_0$. We will use the same notation for the complexification of $\langle \ , \ \rangle$ (a complex-valued form on $\mathfrak{g}$), as well as its various restrictions and dualizations.

Our partition of representations of $G$ is in terms of their restrictions to $K$; so we begin by recalling the highest weight theory for representations of $K$.

Fix a maximal torus $T^c \subset K$. The abelian group $\widehat{T^c}$ of characters of $T^c$ will be identified (by passing to differentials) with a lattice of imaginary-valued linear functionals on the real Lie algebra:

$$(0.3a) \qquad \widehat{T^c} \subset i(\mathfrak{t}_0^c)^*.$$

Because the form $\langle \ , \ \rangle$ was chosen to be negative definite on $\mathfrak{k}_0$, it is positive definite on $i(\mathfrak{t}_0^c)^*$. The roots of $T^c$ in $\mathfrak{k}$ and the weights of $T^c$ in $\mathfrak{g}$ are characters of $T^c$, and so will be regarded as elements of $i(\mathfrak{t}_0^c)^*$. Choose a system of positive roots

$$(0.3b) \qquad \Delta^+(\mathfrak{k}, \mathfrak{t}^c) \subset \Delta(\mathfrak{k}, \mathfrak{t}^c).$$

The Weyl group of $T^c$ in $K$ is by definition the normalizer of $T^c$ in $K$ divided by the centralizer:

$$(0.3c) \qquad W(K, T^c) = N_K(T^c)/Z_K(T^c).$$

Since $K$ need not be connected, this group may be larger than the Weyl group $W(\mathfrak{k}, \mathfrak{t}^c)$ of the root system. Obviously $W(K, T^c)$ acts on the root system $\Delta(\mathfrak{k}, \mathfrak{t}^c)$, so we can define

$$(0.3d) \qquad R(G) = \{w \in W(K, T^c) \mid w\Delta^+(\mathfrak{k}, \mathfrak{t}^c) = \Delta^+(\mathfrak{k}, \mathfrak{t}^c)\}$$

(see [16, Def. 5.1.2]). The notation is chosen to reflect the fact that $R(G)$ is closely related to the $R$-group of Knapp and Stein. Notice that if $G$ is connected, then $R(G)$ is trivial. There is a semidirect product decomposition

$$(0.3e) \qquad W(K, T^c) = R(G)W(\mathfrak{k}, \mathfrak{t}^c)$$

with the second factor normal. Each root $\alpha$ in $\Delta(\mathfrak{k}, \mathfrak{t}^c)$ corresponds to an eigenspace $\mathfrak{k}_\alpha$ of $\operatorname{Ad}(T^c)$ on the complexified Lie algebra. We define

$$(0.3f) \qquad \mathfrak{n}_\mathfrak{k} = \sum_{\alpha \in \Delta^+(\mathfrak{k}, \mathfrak{t}^c)} \mathfrak{k}_\alpha, \qquad \mathfrak{b}_\mathfrak{k} = \mathfrak{n}_\mathfrak{k} + \mathfrak{t}^c;$$

this is a Borel subalgebra of $\mathfrak{k}$.

LEMMA 0.1 (See [16], Prop. 5.1.9). *Using the notation of* (0.3), *suppose $V$ is a finite-dimensional complex representation of the Lie algebra $\mathfrak{k}$. The set of highest weights of $V$ is by definition the set of weights of $\mathfrak{t}^c$ acting on $V^{\mathfrak{n}_\mathfrak{k}}$, a subset of $(\mathfrak{t}^c)^*$.*

a) *If $V$ is an irreducible representation of $\mathfrak{k}$, then $V^{\mathfrak{n}_\mathfrak{k}}$ has dimension one, so $V$ has exactly one highest weight.*

b) *The element $\mu \in (\mathfrak{t}^c)^*$ is the highest weight of an irreducible representation of $\mathfrak{k}$ if and only if $\mu$ is dominant integral for $\mathfrak{k}$: that is,*

$$2\langle\mu, \alpha\rangle/\langle\alpha, \alpha\rangle \in \mathbb{N}$$

*for all $\alpha \in \Delta^+(\mathfrak{k}, \mathfrak{t}^c)$.*

c) *Suppose $V$ is a finite-dimensional irreducible representation of $K$. Then the collection of highest weights of $V$ is a single orbit of $R(G)$ on $\widehat{T}^c \subset i(\mathfrak{t}_0^c)^*$ (cf.* (0.3d)).

Write

$$(0.4a) \qquad \mathfrak{h}_0^c = \mathfrak{g}_0^{T^c} = \mathfrak{k}_0^{T^c} + \mathfrak{p}_0^{T^c} = \mathfrak{t}_0^c + \mathfrak{a}_0^c,$$

a fundamental Cartan subalgebra of $\mathfrak{g}_0$. We will be interested not so much in the roots of $\mathfrak{h}^c$ in $\mathfrak{g}$ as in their restrictions to $\mathfrak{t}^c$: the nonzero weights of $T^c$ on $\mathfrak{g}$. Define

$$(0.4b) \qquad \Delta(\mathfrak{p}, \mathfrak{t}^c) = \text{ set of nonzero weights of } T^c \text{ on } \mathfrak{p}.$$

These weights all occur with multiplicity one. Define

$$(0.4c) \qquad \Delta(\mathfrak{g}, \mathfrak{t}^c) = \Delta(\mathfrak{p}, \mathfrak{t}^c) \cup \Delta(\mathfrak{k}, \mathfrak{t}^c),$$

a subset with multiplicities in $i(\mathfrak{t}_0^c)^*$. One can find fairly detailed structural information about this set for example in [15, §2]. We mention here only a few highlights. First, $\Delta(\mathfrak{g}, \mathfrak{t}^c)$ is a (possibly nonreduced) root system. The elements having multiplicity one are the restrictions of imaginary roots of $\mathfrak{h}^c$; each is compact or noncompact. The roots having multiplicity two are the restrictions of complex roots; each corresponds to one root in $\mathfrak{k}$ and one in $\mathfrak{p}$. The Weyl group $W(\mathfrak{g}, \mathfrak{t}^c)$ of $\Delta(\mathfrak{g}, \mathfrak{t}^c)$ may be identified with the stabilizer of $\mathfrak{t}^c$ in $W(\mathfrak{g}, \mathfrak{h}^c)$; it contains $W(K, T^c)$.

Whenever $\mathfrak{u} \subset \mathfrak{g}$ is a $T^c$-invariant subspace, we write $\Delta(\mathfrak{u}, \mathfrak{t}^c)$ for the set of weights of $T^c$ in $\mathfrak{u}$, counted with multiplicities. Define

$$(0.4d) \qquad 2\rho(\mathfrak{u}) = \sum_{\gamma \in \Delta(\mathfrak{u}, \mathfrak{t}^c)} \gamma \in i(\mathfrak{t}_0^c)^*.$$

As a character of $T^c$, this is the determinant of the adjoint action on $\mathfrak{u}$. As a weight for the Lie algebra $\mathfrak{t}^c$, it is the trace of the adjoint action on $\mathfrak{u}$. We will make particular use of

$$(0.4e) \qquad 2\rho(\mathfrak{n}_{\mathfrak{k}}) = 2\rho_c,$$

the sum of the positive roots of $T^c$ in $\mathfrak{k}$ (see (0.3)).

We turn now to the definition of the partition 0.1 of the unitary dual. We will say more about the motivations behind the construction in Section 2, but here we mention just three primary sources: the notion of lowest $K$-type in [15]; the geometric reformulation of that notion by Carmona in [3]; and a characterization used by the first author of the lowest $K$-types of the representations $A_{\mathfrak{q}}(\lambda)$ studied in [18].

Fix a highest weight $\mu \in i(\mathfrak{t}_0^c)^*$. Choose a system of positive roots $\Delta^+(\mathfrak{g}, \mathfrak{t}^c)$ for $\Delta(\mathfrak{g}, \mathfrak{t}^c)$ making $\mu + 2\rho_c$ dominant:

$$(0.5a) \qquad \langle \mu + 2\rho_c, \gamma \rangle \geq 0, \qquad \gamma \in \Delta^+(\mathfrak{g}, \mathfrak{t}^c).$$

Because $\mu$ is dominant for $K$ and $2\rho_c$ is dominant and regular, $\Delta^+(\mathfrak{g}, \mathfrak{t}^c)$ necessarily contains $\Delta^+(\mathfrak{k}, \mathfrak{t}^c)$. Define

$$(0.5\text{b}) \qquad\qquad 2\rho = \sum_{\gamma \in \Delta^+(\mathfrak{g}, \mathfrak{t}^c)} \gamma,$$

the sum of the positive roots. This positive system defines a positive Weyl chamber

$$(0.5\text{c}) \qquad C = \{\phi \in i(\mathfrak{t}_0^c)^* \mid \langle \phi, \gamma \rangle \geq 0 \quad (\gamma \in \Delta^+(\mathfrak{g}, \mathfrak{t}^c))\}.$$

This is a closed convex cone in the Euclidean space $i(\mathfrak{t}_0^c)^*$. We may therefore define a projection $P$ from $i(\mathfrak{t}_0^c)^*$ onto $C$: $P\phi$ is the unique element of $C$ closest to $\phi$. (This classical construction will be considered carefully in Section 1.) Using $P$, we can now define

$$(0.5\text{d}) \qquad\qquad \lambda_u(\mu) = P(\mu + 2\rho_c - 2\rho).$$

This is an element of the positive Weyl chamber $C \subset i\mathfrak{t}_0^*$; we will see in Proposition 1.4 that it is independent of the choice of positive root system subject to (0.5a). By (0.5a), $\mu + 2\rho_c$ lies in $C$; so by definition of $P$,

$$|\lambda_u(\mu) - (\mu + 2\rho_c)| \leq |2\rho|.$$

Consequently the distance from $\mu$ to $\lambda_u(\mu)$ may be bounded independently of $\mu$ (say by $|2\rho| + |2\rho_c|$). Because the weights $\mu$ lie in a (discrete) lattice in $i(\mathfrak{t}_0^c)^*$, it follows that the set

$$(0.5\text{e}) \qquad\qquad \Lambda_u = \{\lambda_u(\mu), \ \mu \in \widehat{T}^c \text{ dominant}\}$$

is discrete in $i(\mathfrak{t}_0^c)^*$. It is very easy to check that the group $R(G)$ acts on $\Lambda_u$.

We can now define the decomposition (0.1).

*Definition* 0.2. Fix $\lambda_u \in \Lambda_u$. The set of *K-types attached to* $\lambda_u$ is

$$B_u^{\lambda_u}(G) = \{\delta \in \widehat{K} \mid \delta \text{ has a highest weight } \mu \in \widehat{T}^c \text{ such that } \lambda_u(\mu) = \lambda_u\}.$$

The set of *unitary representations of $G$ attached to* $\lambda_u$ is

$$\Pi_u^{\lambda_u}(G) = \{(\pi, \mathcal{H}_\pi) \in \Pi_u(G) \mid \pi \text{ has a lowest } K\text{-type in } B_u^{\lambda_u}(G)\}.$$

It will sometimes be convenient (especially in our consideration of the Fell topology in Section 8) to abuse notation by writing

$$\lambda_u(\pi) = \lambda_u$$

when $\pi \in \Pi_u^{\lambda_u}(G)$. The difficulty is that $\pi$ may belong to several different $\Pi_u^{\lambda_u}(G)$. Corollary 0.4 below says that $\lambda_u(\pi)$ is well-defined as an orbit of $R(G)$, but not (unless $G$ is connected) as an element of $i(\mathfrak{t}_0^c)^*$.

The definition of lowest $K$-type is given in [15, Def. 5.1], or [16, Def. 5.4.18]; we will recall it in Section 2.

PROPOSITION 0.3.    *Suppose $\mu$ and $\mu'$ are highest weights of lowest $K$-types of a single irreducible representation of $G$. Then the weights $\lambda_u(\mu)$ and $\lambda_u(\mu')$ lie in the same orbit of $R(G)$ (cf. (0.3)).*

We will prove this at the end of Section 3.

COROLLARY 0.4.    *Suppose $\lambda_u$ and $\lambda'_u$ belong to $\Lambda_u$ (cf. (0.5)). If $\lambda_u$ and $\lambda'_u$ belong to the same orbit of $R(G)$, then $\Pi_u^{\lambda_u}(G) = \Pi_u^{\lambda'_u}(G)$. If not, then $\Pi_u^{\lambda_u}(G)$ and $\Pi_u^{\lambda'_u}(G)$ are disjoint.*

To get a partition as described in (0.1), we could, of course, replace $\Lambda_u$ by a fundamental domain for the action of $R(G)$; but since there is no natural choice for such a fundamental domain, it is more convenient for us to leave matters in this form.

Next, we want to describe the group $G(\lambda_u)$ discussed after (0.1). The root space decomposition of $\mathfrak{k}$ implies that

$$\mathfrak{k}_0 = \mathfrak{t}_0^c \oplus [\mathfrak{t}_0^c, \mathfrak{k}_0],$$

and therefore that

(0.6a)                $$\mathfrak{g}_0 = \mathfrak{t}_0^c \oplus [\mathfrak{t}_0^c, \mathfrak{k}_0] \oplus \mathfrak{p}_0.$$

By means of this decomposition any linear functional on $\mathfrak{t}_0^c$ extends canonically to $\mathfrak{g}_0$, by making it zero on the other summands. In this way we can identify $i(\mathfrak{t}_0^c)^*$ naturally as a subspace of $i\mathfrak{g}_0^*$:

(0.6b)                $$i(\mathfrak{t}_0^c)^* \hookrightarrow i\mathfrak{g}_0^*.$$

The group $G$ acts on $i\mathfrak{g}_0^*$ by the coadjoint action. In the same way $T^c$ (or even its normalizer in $G$) acts on $i\mathfrak{t}_0^*$, and the inclusion 0.6(b) is $T^c$-equivariant. Since $T^c$ is abelian, it follows that $T^c$ acts trivially on the image. For any $\lambda \in i(\mathfrak{t}_0^c)^* \subset i\mathfrak{g}_0^*$, define

(0.6c)          $$G(\lambda) = \text{isotropy group at } \lambda \text{ for the } G \text{ action};$$

this is a subgroup of $G$ containing $T^c$. It is easy to check that the Cartan involution $\theta$ preserves $G(\lambda)$; and that $G(\lambda)$ is a reductive group in Harish-Chandra's class, with Cartan involution $\theta|_{G(\lambda)}$.

The set of roots of $T^c$ in $G(\lambda)$ is

(0.6d)          $$\Delta(\mathfrak{g}(\lambda), \mathfrak{t}^c) = \{\gamma \in \Delta(\mathfrak{g}, \mathfrak{t}^c) \mid \langle \gamma, \lambda \rangle = 0\}.$$

We can therefore construct a $\theta$-stable parabolic subalgebra $\mathfrak{q}(\lambda) = \mathfrak{g}(\lambda) + \mathfrak{u}(\lambda)$ by the following requirement (compare [16, Def. 5.2.1], and [9, Prop. 4.76]):

(0.6e)          $$\Delta(\mathfrak{u}(\lambda), \mathfrak{t}^c) = \{\gamma \in \Delta(\mathfrak{g}, \mathfrak{t}^c) \mid \langle \gamma, \lambda \rangle > 0\}.$$

Using the parabolic subalgebra $\mathfrak{q}(\lambda)$, we define functors of *cohomological parabolic induction* $\mathcal{L}_j(\lambda)$ carrying $(\mathfrak{g}(\lambda), G(\lambda) \cap K)$-modules to $(\mathfrak{g}, K)$-modules (see [9, §5.1], or [16, Def. 6.3.1]). These are most interesting in the degree

$$(0.6f) \qquad\qquad S = \dim(\mathfrak{u}(\lambda) \cap \mathfrak{k}).$$

At the same time, we have $\mathcal{L}_j^K(\lambda)$ carrying $(G(\lambda) \cap K)$-modules to $K$-modules ([9, §5.6]).

PROPOSITION 0.5. *With notation as in Definition* 0.2 *and* (0.6), *the functor* $\mathcal{L}_S^K(\lambda_u)$ *implements a bijection from* $B_u^{\lambda_u}(G(\lambda_u))$ *onto* $B_u^{\lambda_u}(G)$.

We will prove this in Section 3 (Proposition 3.1). The aim of this paper is a corresponding result for unitary representations:

*Conjecture* 0.6. With notation as in Definition 0.2 and (0.6), the functor $\mathcal{L}_S(\lambda_u)$ implements a bijection from $\Pi_u^{\lambda_u}(G(\lambda_u))$ onto $\Pi_u^{\lambda_u}(G)$.

In Section 5 we will explain a program for proving this conjecture, reducing it to the more elementary Conjecture 5.7′. We present strong evidence for the validity of this latter conjecture, including Proposition 7.12 (which is a partial result towards it) and Theorem 5.9 (which provides moral support for the formulation of the conjecture). The program leads to an alternative description of the subsets $\Pi_u^{\lambda_u}(G)$ (Theorem 5.8(c)). This description will imply that each is an open and closed subset of $\Pi_u(G)$, and that the bijection in Conjecture 0.6 is a homeomorphism in the Fell topology. We will examine these issues in Section 8.

## 1. Projections on convex cones

In this section we recall some elementary general results about projections on convex cones. They will be applied to positive Weyl chambers in $i(\mathfrak{t}_0^{\mathfrak{c}})^*$, with notation as described in the introduction. The ideas and the exposition are taken from [3]. Throughout this section, $V$ will be a finite-dimensional real vector space, endowed with a positive definite inner product $\langle \, , \, \rangle$. Write $|v| = \langle v, v \rangle^{1/2}$ for the corresponding Hilbert space norm on $V$. Let $C$ be a closed convex cone in $V$. Write

$$(1.1) \qquad\qquad C^o = \{v \in V \mid \langle v, c \rangle \geq 0, \quad \text{for all } c \in C\}$$

for the dual cone; it is again a closed convex cone in $V$.

PROPOSITION 1.1 (See [3, Prop. 1.2]). *Suppose* $v \in V$. *Then there is a unique element* $c_0$ *of* $C$ *closest to* $v$. *It may be characterized by any of the following equivalent conditions.*

a) *For any $c \in C$, $|v - c_0| \leq |v - c|$.*
b) *For any $c \in C$, $\langle v - c_0, c - c_0 \rangle \leq 0$.*
c) *The vector $c_0 - v$ belongs to the dual cone $C^o$, and $\langle c_0 - v, c_0 \rangle = 0$.*

We omit the elementary proof.

*Definition* 1.2. In the setting of Proposition 1.1, the element $c_0$ of $C$ satisfying the equivalent conditions (a)–(c) is written $Pv$, the *projection of $v$ on $C$*.

We will need some additional properties of the projection $P$ that depend on special properties of Weyl chambers. We therefore assume for the rest of this section that we are given a (possibly nonreduced) root system

$$(1.2a) \qquad \Delta \subset V,$$

and a set of positive roots

$$(1.2b) \qquad \Delta^+ \subset \Delta.$$

We do not assume that $\Delta$ spans $V$. Write $V_s$ for the linear span of $\Delta$ in $V$, and $V_z$ for its orthogonal complement. (The subscripts $s$ and $z$ stand for "semisimple" and "central.") Write

$$(1.2c) \qquad \Pi = \{\alpha_1, \ldots, \alpha_l\} \subset \Delta^+$$

for the simple roots; these form a basis of $V_s$. Write $\xi_i$ for the dual basis of fundamental weights:

$$(1.2d) \qquad \langle \xi_i, \alpha_j \rangle = \delta_{ij}, \qquad \langle \xi_i, v \rangle = 0 \quad (v \in V_z).$$

The *closed positive Weyl chamber* is by definition

$$(1.2e) \qquad C = \{v \in V \mid \langle v, \alpha \rangle \geq 0 \quad (\alpha \in \Delta^+)\}$$

$$= V_z + \sum_{i=1}^{l} \mathbb{R}^{\geq 0} \xi_i.$$

Occasionally it will be useful to emphasize the dependence on the choice of positive root system; in that case we will write $C(\Delta^+)$. The dual cone (cf. (1.1)) is

$$(1.2f) \qquad C^o = \sum_{i=1}^{l} \mathbb{R}^{\geq 0} \alpha_i,$$

the *cone of positive roots*. We will refer to this collection of notation and assumptions as "the setting 1.2."

The main construction we use in this paper is (0.5d): to begin with a dominant weight, subtract a dominant weight, then project onto the positive

Weyl chamber. The following lemma will lead to several important properties of this construction.

LEMMA 1.3. *In the setting 1.2, suppose $\gamma \in C$ is a dominant weight, $\alpha \in \Pi$ is a simple root, and $v \in V$. If $\langle v, \alpha \rangle \leq 0$, then*

$$\langle P(v - \gamma), \alpha \rangle = 0.$$

*Proof.* By Proposition 1.1, we may write

$$(1.3a) \qquad v - \gamma = P(v - \gamma) - e(v, \gamma),$$

with $e(v, \gamma) \in C^o$ orthogonal to $P(v - \gamma)$. By 1.2(f), this means

$$(1.3b) \qquad e(v, \gamma) = \sum_i c_i \alpha_i \qquad (c_i \geq 0),$$

and $c_i > 0$ implies $\langle \alpha_i, P(v - \gamma) \rangle = 0$. Now form the inner product of (1.3a) with $\alpha$. The hypothesis on $v$ and the dominance of $\gamma$ make the left side nonpositive. The first term on the right is nonnegative. If it is zero we are done; so we may assume it is positive. Then

$$(1.3c) \qquad \langle \alpha, e(v, \gamma) \rangle > 0.$$

Distinct simple roots have nonpositive inner product; so (1.3b) and (1.3c) imply that the coefficient $c_i$ of $\alpha$ in (1.3b) must be strictly positive. We have already seen that this implies $\alpha$ is orthogonal to $P(v - \gamma)$, as we wished to show. $\square$

PROPOSITION 1.4. *In the setting 1.2, fix a dominant weight $\gamma \in C$. Define a map $T_\gamma$ from $V$ to $V$ as follows. Given $v \in V$, choose a positive root system $\Delta^+(v)$ making $v$ dominant; write $\Delta^+(v) = \sigma \Delta^+$, with $\sigma \in W(\Delta)$ (the Weyl group). Define*

$$T_\gamma(v) = P(\Delta^+(v))(v - \sigma\gamma).$$

*Then $T_\gamma(v)$ is well-defined (independent of the choice of positive root system making $v$ dominant).*

Recall that this result was used in (0.5d).

*Proof.* Suppose $w\Delta^+(v)$ is another positive system making $v$ dominant. Then $v$ and $w^{-1}v$ are both dominant for $\Delta^+(v)$, and are conjugate by the Weyl group; so they coincide. Consequently $w$ fixes $v$. Using $w\Delta^+(v)$ to define $T_\gamma(v)$ leads to

$$P(w\Delta^+(v))(v - w\sigma\gamma) = P(w\Delta^+(v))(wv - w\sigma\gamma) = w(P(\Delta^+(v))(v - \sigma\gamma)).$$

By Chevalley's theorem (see for example [9, Prop. 4.146]) $w$ is a product of reflections in simple roots orthogonal to $v$. By Lemma 1.3, each such simple root is orthogonal to $P(\Delta^+(v))(v - \sigma\gamma)$. So $w$ fixes $P(\Delta^+(v))(v - \sigma\gamma)$, and we get finally

$$P(w\Delta^+(v))(v - w\sigma\gamma) = P(\Delta^+(v))(v - \sigma\gamma),$$

as we wished to show.     □

The next result will be used in Section 2 to relate the parameter $\lambda_u(\mu)$ defined in (0.5) to the parameter $\lambda(\mu)$ constructed in [15, Prop. 4.1] or [16, Prop. 5.3.3]. This will lead directly to Proposition 0.3, and to a variety of stronger properties of the sets $\Pi_u^{\lambda_u}(G)$.

PROPOSITION 1.5.     *In the setting* 1.2, *suppose* $v$, $\gamma$, *and* $\delta$ *are all dominant weights in* $C$. *Then*

$$P(v - \gamma - \delta) = P(P(v - \gamma) - \delta).$$

*Proof.* Write

(1.4a)                    $P(v - \gamma) = (v - \gamma) + e(v, \gamma),$

with $e(v, \gamma) \in C^o$ orthogonal to $P(v - \gamma)$ (Proposition 1.1(c)). Next, write

(1.4b)          $P(P(v - \gamma) - \delta) = (P(v - \gamma) - \delta) + e(P(v - \gamma), \delta),$

with $e(P(v - \gamma), \delta) \in C^o$ orthogonal to $P(P(v - \gamma) - \delta)$. Combining these two equations gives

(1.4c)          $P(P(v - \gamma) - \delta) = (v - \gamma - \delta) + e(v, \gamma) + e(P(v - \gamma), \delta).$

The last term on the right is orthogonal to $P(P(v - \gamma) - \delta)$ by 1.4(b). The second is orthogonal to $P(v - \gamma)$ by 1.4(a), and therefore also to $P(P(v - \gamma) - \delta)$ by Lemma 1.3. We have therefore shown that $e(v, \gamma) + e(P(v - \gamma), \delta)$ is an element of $C^o$ orthogonal to $P(P(v - \gamma) - \delta)$. By Proposition 1.1(c), the formula in 1.4(c) now shows that

$$P(v - \gamma - \delta) = P(P(v - \gamma) - \delta),$$

as desired.     □

COROLLARY 1.6.     *In the setting* 1.2, *suppose* $\gamma$ *and* $\delta$ *are dominant weights in* $C$. *Define maps* $T_\gamma$, *etc., as in Proposition* 1.4. *Then*

$$T_{\gamma+\delta} = T_\gamma \circ T_\delta.$$

The next result will be used in the characterization of the $K$-types attached to the parameter 0 (Definition 0.2), for which Conjecture 0.6 gives no information. It is exactly analogous to [3, Cor. 2.13].

PROPOSITION 1.7.       *In the setting 1.2, suppose $v$ and $\gamma$ are dominant weights in $C$. The following conditions are equivalent.*

a) $P(v - \gamma) = 0$.
b) *The weight $\gamma - v$ belongs to $C^o$; that is, it is a sum of positive roots with nonnegative coefficients.*
c) *The weight $v$ belongs to the convex hull of the Weyl group orbit of $\gamma$.*

*Proof.* The equivalence of (a) and (b) follows from Proposition 1.1(c). That (b) and (c) are equivalent is a well-known aspect of the Cartan-Weyl theory of finite-dimensional representations; but since it is not easy to find a convenient reference, we give a proof. The main step is

LEMMA 1.8 ([16, Lemma 6.3.28]).       *Suppose $\gamma \in C$ is dominant and $\sigma \in W(\Delta)$. Then there is an element $e(\sigma) \in C^o$ so that*

$$\sigma \cdot \gamma = \gamma - e(\sigma).$$

Using this lemma, we proceed with the proof of the proposition. Suppose first that (c) holds; say $v = \sum_{\sigma \in W} c_\sigma \sigma \cdot \gamma$, with $c_\sigma \geq 0$ and $\sum c_\sigma = 1$. Adding the equations in Lemma 1.8 gives

$$v = \gamma - \sum_\sigma c_\sigma e(\sigma).$$

Since $C^o$ is a convex cone, the last sum belongs to $C^o$, as we wished to show.

Conversely, suppose $v$ is *not* in the convex hull of the Weyl group orbit of $\gamma$. Then it is separated from the convex hull by a hyperplane. That is, there is an element $\lambda' \in V$ so that

$$\langle v, \lambda' \rangle > \langle \sigma \cdot \gamma, \lambda' \rangle \qquad (\sigma \in W).$$

This inequality may be written as

$$\langle v, \lambda' \rangle > \langle \gamma, \sigma \cdot \lambda' \rangle \qquad (\sigma \in W).$$

(We have replaced the variable $\sigma$ by $\sigma^{-1}$; since the inequality holds for all $\sigma \in W$, this is permissible.) Now write $\lambda$ for the dominant element of the Weyl group orbit of $\lambda'$, and $\lambda' = \sigma_0 \cdot \lambda$. Then we have

$$\langle v, \sigma_0 \cdot \lambda \rangle > \langle \gamma, \sigma \cdot \lambda \rangle \qquad (\sigma \in W).$$

(This time we have replaced $\sigma$ by $\sigma \sigma_0^{-1}$.) Applying Lemma 1.8 to the dominant weight $\lambda$, we find $\sigma_0 \cdot \lambda = \lambda - e(\sigma_0)$, with $e(\sigma_0) \in C^o$. Since $v$ is dominant, it has a nonnegative inner product with $e(\sigma_0)$. So we get

$$\langle v, \lambda \rangle > \langle \gamma, \sigma \cdot \lambda \rangle \qquad (\sigma \in W).$$

Taking $\sigma = 1$, we find in particular

$$\langle v - \gamma, \lambda \rangle > 0.$$

Since $\lambda$ is dominant, this inequality implies that $\gamma - v \notin C^o$, as we wished to show.                                                                                      □

We can now give a convenient characterization of the map $v \mapsto P(v - \gamma)$.

COROLLARY 1.9.    *In the setting 1.2, suppose $v$ and $\gamma$ are dominant weights in $C$. Let $W_0$ be the Weyl group of roots orthogonal to $P(v - \gamma)$ (the stabilizer in $W$ of $P(v - \gamma)$). Then $v - P(v - \gamma)$ is dominant, and belongs to the convex hull of the $W_0$ orbit of $\gamma$.*

*Conversely, suppose that $v_0$ and $\gamma$ are dominant weights in $C$. Let $W_0$ be the Weyl group of the roots orthogonal to $v_0$. Suppose $\gamma_0$ is a dominant weight in the convex hull of the $W_0$ orbit of $\gamma$. Then $v = v_0 + \gamma_0$ is dominant, and $P(v - \gamma) = v_0$.*

*Proof.* Write $\Delta_0$ for the system of roots orthogonal to $P(v - \gamma)$. By Proposition 1.1(c), $v - \gamma - P(v - \gamma)$ is a sum of negative roots orthogonal to $P(v - \gamma)$; that is, a sum of roots in $\Delta_0^-$. Now apply Proposition 1.7 to the root system $\Delta_0$ and the $\Delta_0^+$-dominant weight $v - P(v - \gamma)$. The conclusion is that $v - P(v - \gamma)$ belongs to the convex hull of the $W_0$ orbit of $\gamma$. Clearly $v - P(v - \gamma)$ is dominant for $\Delta_0^+$. The remaining set $\Delta^+ - \Delta_0^+$ is permuted by $W_0$. Since $\gamma$ is dominant for these roots, the convex hull of $W_0 \cdot \gamma$ must be as well.

The converse is similar but easier, and we leave it to the reader.          □

Because we will most often use these results in the case $\gamma = \rho$ (half the sum of the positive roots), we will need a precise description of the convex hull in that case.

PROPOSITION 1.10 ([10, Lemma 5.9]).    *In the setting of 1.2, write $\rho$ for half the sum of the positive roots. Then the convex hull of the Weyl group orbit of $\rho$ coincides with the set $R$ of all weights of the form*

$$r = \sum_{\alpha \in \Delta^+} c_\alpha \alpha, \qquad (-1/2 \le c_\alpha \le 1/2).$$

*Equivalently, these are the weights*

$$r = \rho - \sum_{\alpha \in \Delta^+} b_\alpha \alpha, \qquad (0 \le b_\alpha \le 1).$$

## 2. Classification by lowest $K$-types

In this section we recall part of the Langlands classification of admissible irreducible representations of $G$, as reformulated in [15], [16], and [3].

*Definition* 2.1.   Suppose $G$ is a reductive group in Harish-Chandra's class. The *admissible dual* of $G$ is the set $\Pi_a(G)$ of infinitesimal equivalence classes of irreducible admissible representations of $G$. Equivalently, $\Pi_a(G)$ may be identified with the set of equivalence classes of irreducible $(\mathfrak{g}, K)$-modules ([16, Def. 0.3.8]; here $K$ is a maximal compact subgroup of $G$ as in 0.2.).

Most of the equivalence of these two definitions of $\Pi_a(G)$ is proved in [5]. What is missing there is the proof that every irreducible $(\mathfrak{g}, K)$-module appears as the space of $K$-finite vectors in an admissible representation. This was proved for linear groups by Harish-Chandra, and in general by Lepowsky and Rader independently (see for example [19, Th. 3.5.6]).

*Definition* 2.2 ([15, Def. 5.1] or [16, Def. 5.4.18]).   Suppose $X \in \Pi_a(G)$ is an irreducible $(\mathfrak{g}, K)$-module. A *lowest $K$-type of $X$* is an irreducible representation $\delta \in \widehat{K}$ such that

a) $\delta$ appears in the restriction of $X$ to $K$, and
b) if $\mu \in i(\mathfrak{t}_0^c)^*$ is a highest weight of $\delta$ (Lemma 0.4), then $\langle \mu + 2\rho_c, \mu + 2\rho_c \rangle$ (notation 0.4(e)) is minimal subject to (a).

Obviously the norm in Definition 2.2(b) is independent of the choice of highest weight of $\delta$.

The classification of admissible representations in [15] is based on the following fact.

PROPOSITION 2.3 ([15, Prop. 4.1] or [16, Prop. 5.3.3]).   *Suppose $\mu \in i(\mathfrak{t}_0^c)^*$ is dominant integral for $\Delta^+(\mathfrak{k}, \mathfrak{t}^c)$. Choose a positive root system $\Delta^+ = \Delta^+(\mathfrak{g}, \mathfrak{t}^c)$ making $\mu + 2\rho_c$ dominant (cf. (0.5a)). Write $\rho$ for the corresponding half sum of positive roots. Then there is a weight $\lambda_a(\mu) \in i\mathfrak{t}_0^*$ with the following properties.*

a) *The weight $\lambda_a(\mu)$ is dominant for $\Delta^+$.*
b) *There is an orthogonal collection $\beta_i$ of positive imaginary roots such that*
   i) $\lambda_a(\mu) = \mu + 2\rho_c - \rho + \sum_i c_i\beta_i$     $(0 \le c_i \le 1/2)$, *and*
   ii) $\langle \lambda_a(\mu), \beta_i \rangle = 0$.

COROLLARY 2.4 ([3, Th. 2.10]).   *In the notation of Proposition 2.3, the weight $\lambda_a(\mu)$ is the projection of $\mu + 2\rho_c - \rho$ on the cone $C$ of dominant weights. In the notation of Proposition 1.4, $\lambda_a(\mu) = T_\rho(\mu + 2\rho_c)$. Consequently $\lambda_a(\mu)$ is independent of the choice of $\Delta^+$ making $\mu + 2\rho_c$ dominant.*

*Proof.*  Proposition 2.3 shows that $\lambda_a(\mu)$ satisfies the conditions in Proposition 1.1(c) characterizing $P(\mu + 2\rho_c - \rho)$. The remaining assertions now follow from Proposition 1.4. $\qquad\square$

In analogy with (0.5e), we may now define

$$(2.1) \qquad\qquad \Lambda_a = \{\lambda_a(\mu), \ \mu \in \widehat{T}^c \text{ dominant}\}.$$

The analogue of Definition 0.2 for admissible representations is:

*Definition 2.5.* Fix $\lambda_a \in \Lambda_a$. The set of $K$-types attached to $\lambda_a$ is

$$B_a^{\lambda_a}(G) = \{\delta \in \widehat{K} \mid \delta \text{ has a highest weight } \mu \in \widehat{T}^c \text{ such that } \lambda_a(\mu) = \lambda_a\}.$$

The set of *admissible representations of $G$ attached to $\lambda_a$* is

$$\Pi_a^{\lambda_a}(G) = \{\pi \in \Pi_a(G) \mid \pi \text{ has a lowest } K\text{-type in } B_a^{\lambda_a}(G)\}.$$

Just as in the case of Definition 0.2, we will abuse notation in Section 8 by writing

$$\lambda_a(\pi) = \lambda_a$$

when $\pi \in \Pi_a^{\lambda_a}(G)$. Again the difficulty is that $\pi$ may belong to several different $\Pi_a^{\lambda_a}(G)$. It follows from Theorem 2.9 below that $\lambda_a(\pi)$ is well-defined as an orbit of $R(G)$.

PROPOSITION 2.6. *With notation as in Definition 2.5, (0.6), and (0.4d), write $S_a = \dim \mathfrak{u}(\lambda_a) \cap \mathfrak{k}$. Then the functor $\mathcal{L}_{S_a}^K(\lambda_a)$ implements a bijection from $B_a^{\lambda_a - \rho(\mathfrak{u}(\lambda_a))}(G(\lambda_a))$ onto $B_a^{\lambda_a}(G)$.*

This is essentially proved in [16], although some careful inspection of definitions is needed to verify that. Typical of the complications is the fact that the functors $\mathcal{L}_S^K(\lambda)$ are not defined in [16]. Fortunately they are easy to describe in terms of highest weights: roughly speaking, they add $2\rho(\mathfrak{u}(\lambda) \cap \mathfrak{p})$ to the highest weight. Here is a more precise statement.

LEMMA 2.7 ([9, Cor. 5.72]). *Suppose $\lambda \in i(\mathfrak{t}_0^c)^*$ is dominant for $\Delta^+(\mathfrak{k}, \mathfrak{t}^c)$, so that $\mathfrak{u}(\lambda) \cap \mathfrak{k} \subset \mathfrak{n}_\mathfrak{k}$ (notation (0.6e) and (0.3f)). Write $R = \dim \mathfrak{u}(\lambda) \cap \mathfrak{p}$, so that $\wedge^R(\mathfrak{u}(\lambda) \cap \mathfrak{p})$ is a one-dimensional representation of $G(\lambda) \cap K$. Suppose $Z$ is an irreducible representation of $G(\lambda) \cap K$, of highest weight denoted $\mu_{G(\lambda)} \in \widehat{T}^c$, and $V$ is an irreducible representation of $K$. Then*

$$\text{Hom}_K(\mathcal{L}_S^K(\lambda)(Z), V) \simeq \text{Hom}_{G(\lambda) \cap K}(Z \otimes \wedge^R(\mathfrak{u}(\lambda) \cap \mathfrak{p}), V^{\mathfrak{u}(\lambda) \cap \mathfrak{k}}).$$

*Define $\mu_G = \mu_{G(\lambda)} + 2\rho(\mathfrak{u}(\lambda) \cap \mathfrak{p})$. If $\mu_G$ is dominant for $K$, then every irreducible constituent of $\mathcal{L}_S^K(\lambda)(Z)$ has highest weight $\mu_G$. If $\mu_G$ is not dominant for $K$, then $\mathcal{L}_S^K(\lambda)(Z) = 0$.*

Given this description of the functor $\mathcal{L}_S^K(\lambda)$, the translation from Lemma 6.5.4 in [16] to Proposition 2.6 above is more or less straightforward; we omit the details.

In order to formulate the classification of admissible representations, we need one more ingredient.

PROPOSITION 2.8 ([9, Th. 5.80]).     *In the setting* (0.6), *suppose $Z$ is a* $(\mathfrak{g}(\lambda), G(\lambda) \cap K)$-*module. Then there is a natural injective map of $K$-modules*

$$\mathcal{B}_Z \colon \mathcal{L}_S^K(\lambda)(Z) \to \mathcal{L}_S(\lambda)(Z).$$

This map is called the *bottom layer map* for $Z$, for reasons explained in Section 5.6 of [9].

THEOREM 2.9 ([16, Th. 6.5.12]).     *Suppose $\lambda_a \in \Lambda_a$ (cf.* (2.1)). *Then there is a natural bijection*

$$(2.2) \qquad \Pi_a^{\lambda_a - \rho(\mathfrak{u}(\lambda_a))}(G(\lambda_a)) \to \Pi_a^{\lambda_a}(G)$$

(Def. 2.5) *defined as follows. Suppose $Z_a \in \Pi_a^{\lambda_a - \rho(\mathfrak{u}(\lambda_a))}(G(\lambda_a))$ is an irreducible $(\mathfrak{g}(\lambda_a), G(\lambda_a) \cap K)$-module.*

  a)  *The set of lowest $G(\lambda_a) \cap K$-types of $Z_a$ consists precisely of those in $B_a^{\lambda_a - \rho(\mathfrak{u}(\lambda_a))}(G(\lambda_a))$.*

  b)  *The set of lowest $K$-types of $\mathcal{L}_{S_a}(\lambda_a)(Z_a)$ consists precisely of the image of the bottom layer map composed with $\mathcal{L}_{S_a}^K(\lambda_a)$ on the lowest $G(\lambda_a) \cap K$-types of $Z_a$.*

Let $\delta_{G(\lambda_a)} \in B_a^{\lambda_a - \rho(\mathfrak{u}(\lambda_a))}(G(\lambda_a))$ be an irreducible representation of the group $G(\lambda_a) \cap K$ appearing in $Z_a$, and let $\delta = \mathcal{L}_{S_a}^K(\lambda_a)(\delta_{G(\lambda_a)})$ be the corresponding representation of $K$ in $B_a^{\lambda_a}(G)$ (Prop. 2.6). By Proposition 2.8, $\delta$ must appear in $\mathcal{L}_{S_a}(\lambda_a)(Z_a)$.

  c)  *There is a unique irreducible subquotient $X$ of $\mathcal{L}_{S_a}(\lambda_a)(Z_a)$ containing the $K$-type $\delta$.*

  d)  *The multiplicity of $\delta$ in $X$ is equal to the multiplicity of $\delta$ in $\mathcal{L}_{S_a}(\lambda_a)(Z_a)$, and to the multiplicity of $\delta_{G(\lambda_a)}$ in $Z_a$.*

  e)  *The subquotient $X$ is independent of the choice of the representation $\delta_{G(\lambda_a)}$ in $B_a^{\lambda_a - \rho(\mathfrak{u}(\lambda_a))}(G(\lambda_a))$, subject to the requirement that $\delta_{G(\lambda_a)}$ occur in $Z_a$. The bijection* (2.2) *sends $Z_a$ to $X$.*

To construct the inverse map, suppose $X$ is an irreducible $(\mathfrak{g}, K)$-module in $\Pi_a^{\lambda_a}(G)$, having a lowest $K$-type $\delta \in B_a^{\lambda_a}(G)$. Let $\delta_{G(\lambda_a)} \in B_a^{\lambda_a - \rho(\mathfrak{u}(\lambda_a))}(G(\lambda_a))$ be the corresponding irreducible representation of $G(\lambda_a) \cap K$ (Prop. 2.6). Write $R_a = \dim(\mathfrak{u}(\lambda_a) \cap \mathfrak{p})$ as in Lemma 2.7.

  f)  *The representation $\delta_{G(\lambda_a)}$ appears in the cohomology space $H^*(\mathfrak{u}(\lambda_a), X)$ (a $(\mathfrak{g}(\lambda_a), G(\lambda_a) \cap K)$)-module) exactly as often as $\delta$ appears in $X$; all occurrences are in degree $R_a$.*

  g)  *There is a unique irreducible subquotient $Z_a$ of $H^{R_a}(\mathfrak{u}(\lambda_a), X)$ containing the $G(\lambda_a) \cap K$-type $\delta_{G(\lambda_a)}$.*

h) *The subquotient $Z_a$ is independent of the choice of $\delta \in B_a^{\lambda_a}(G)$, subject to the requirement that $\delta$ occur in $X$.*

*The inverse of (2.2) sends $X$ to $Z_a$.*

Just as for Proposition 2.6, some translation is required to pass from the results in [16] to this one; but we omit the details.

Because we are interested in unitary representations, we need to know how the bijection of Theorem 2.9 affects Hermitian forms. Here is the underlying machinery.

PROPOSITION 2.10 ([9, Th. 6.34]). *In the setting of (0.6), suppose $Z$ is a representation of $G(\lambda) \cap K$. Then any invariant Hermitian form on $Z$ induces in a natural way a $K$-invariant Hermitian form on $\mathcal{L}_S^K(\lambda)(Z)$. This correspondence carries positive definite forms to positive definite forms, and negative definite forms to negative definite forms.*

*Suppose $Z$ is actually a $(\mathfrak{g}(\lambda), G(\lambda) \cap K)$-module. Then any invariant Hermitian form on $Z$ induces in a natural way a $(\mathfrak{g}, K)$-invariant Hermitian form on $\mathcal{L}_S(\lambda)(Z)$. The bottom layer map $\mathcal{B}_Z$ of Proposition 2.8 is unitary; that is, it respects these forms.*

THEOREM 2.11. *The bijection of Theorem 2.9 preserves Hermitian representations. Explicitly, suppose $Z_a$ is an irreducible Hermitian $(\mathfrak{g}(\lambda_a), G(\lambda_a) \cap K)$-module. Then the induced Hermitian form on $\mathcal{L}_{S_a}(\lambda_a)(Z_a)$ descends to a nondegenerate form on the subquotient $X$ containing the lowest $K$-types.*

Unfortunately, the analogue of Theorem 2.11 with "Hermitian" replaced by "unitary" is false. We know no examples in which $Z_a$ is unitary and $X$ is not; but the other direction fails spectacularly and often. One of the simplest examples has $G = \mathrm{SO}_e(4,1)$ and $\mu = 0$ the highest weight of the trivial representation of $K$. In this case $G(\lambda_a) = \mathrm{SO}(2) \times \mathrm{SO}_e(2,1)$, and Theorem 2.9 relates spherical representations of $G$ to representations of $G(\lambda_a)$ containing the $\mathrm{SO}(2) \times \mathrm{SO}(2)$ representation $\mathbb{C}_{-1} \otimes \mathbb{C}_0$. (The subscripts refer to the representations of $\mathrm{SO}(2)$, which we index by $\mathbb{Z}$ as usual.) Then the representation $Z_a = \mathbb{C}_{-1} \otimes \mathbb{C}^3$ of $G(\lambda_a)$ carries a Hermitian form of signature $(1,2)$; but the corresponding irreducible representation $X$ of $G$ is the trivial representation.

Here is how we propose to circumvent this problem. Suppose $\mathfrak{q} = \mathfrak{l} + \mathfrak{u}$ is a $\theta$-stable parabolic subalgebra containing $\mathfrak{q}(\lambda_a)$. Using "induction by stages" theorems, we can get analogues of Theorems 2.9 and 2.11 relating certain representations of $G$ to certain representations of $L$. These results are in some sense weaker than those above, because they reduce problems about $G$ not to $G(\lambda_a)$, but only to the larger group $L$. On the other hand, and for the same reason, a representation of $L$ "remembers" more about the corresponding representation of $G$ than the representation of $G(\lambda_a)$ does. For that reason, a unitary version of Theorem 2.11 has a better chance to be true. An extreme possibility is to take $\mathfrak{q} = \mathfrak{g}$; then $L = G$, so the "reduction" from $G$ to $L$ is perfectly behaved but worthless. Our goal is to choose $\mathfrak{q}$ large enough that the reduction is well-behaved, but small enough that it is interesting.

We begin with formalities. Fix $\lambda_a \in \Lambda_a$. A weight $\lambda \in i(\mathfrak{t}_0^c)^*$ is called a *singularization of* $\lambda_a$ if for every root $\alpha \in \Delta(\mathfrak{g}, \mathfrak{t}^c)$,

$$\text{(2.3a)} \qquad\qquad \langle \lambda, \alpha \rangle > 0 \Rightarrow \langle \lambda_a, \alpha \rangle > 0.$$

That is, we require that $\lambda$ belong to the same closed positive Weyl chamber as $\lambda_a$, and that it lie on at least as many walls. Another way to say the same thing is: $\lambda$ should belong to every closed positive Weyl chamber containing $\lambda_a$. A consequence of (2.3a) is

$$\text{(2.3b)} \qquad\qquad \langle \lambda_a, \alpha \rangle = 0 \Rightarrow \langle \lambda, \alpha \rangle = 0.$$

By (0.6),

$$\text{(2.3c)} \qquad \mathfrak{u}(\lambda) \subset \mathfrak{u}(\lambda_a), \qquad \mathfrak{g}(\lambda_a) \subset \mathfrak{g}(\lambda), \qquad \mathfrak{q}(\lambda_a) \subset \mathfrak{q}(\lambda).$$

Because $G$ is in Harish-Chandra's class, it follows also that

$$\text{(2.3d)} \qquad\qquad G(\lambda_a) \subset G(\lambda).$$

(To work with groups not in Harish-Chandra's class, this condition should be included as an additional assumption on $\lambda$. We will eventually apply these formalities with $\lambda = \lambda_u$, the weight defined by (0.5) for some $\mu$ giving rise to $\lambda_a$. In that case (2.3d) is automatically satisfied.) We have direct sum decompositions like

$$\text{(2.3e)} \qquad\qquad \mathfrak{u}(\lambda_a) = \mathfrak{u}(\lambda) + (\mathfrak{u}(\lambda_a) \cap \mathfrak{g}(\lambda)),$$

and similarly after intersection with $\mathfrak{k}$. We need two consequences of this decomposition. First,

$$\text{(2.3f)} \qquad\qquad \rho(\mathfrak{u}(\lambda_a)) = \rho(\mathfrak{u}(\lambda)) + \rho(\mathfrak{u}(\lambda_a) \cap \mathfrak{g}(\lambda))$$

as weights in $i(\mathfrak{t}_0^c)^*$ (cf. (0.4d)). Next, if we define

$$\text{(2.3g)} \qquad S_a = \dim \mathfrak{u}(\lambda_a) \cap \mathfrak{k}, \qquad S = \dim \mathfrak{u}(\lambda) \cap \mathfrak{k},$$

then

$$\text{(2.3h)} \qquad\qquad \dim(\mathfrak{u}(\lambda_a) \cap \mathfrak{g}(\lambda) \cap \mathfrak{k}) = S_a - S.$$

We are going to be interested in cohomological induction functors related to three parabolic subalgebras:

$$\text{(2.4a)} \qquad \mathcal{L}_{S_a}(\lambda_a) \colon (\mathfrak{g}(\lambda_a), G(\lambda_a) \cap K)\text{-modules} \rightarrow (\mathfrak{g}, K)\text{-modules}$$

(defined in terms of the parabolic subalgebra $\mathfrak{q}(\lambda_a)$);

(2.4b)    $\mathcal{L}_S(\lambda)\colon (\mathfrak{g}(\lambda), G(\lambda) \cap K)\text{-modules} \to (\mathfrak{g}, K)\text{-modules}$

(defined in terms of the parabolic subalgebra $\mathfrak{q}(\lambda)$); and

(2.4c)  $\mathcal{L}_{S_a - S}(\lambda_a, \lambda)\colon (\mathfrak{g}(\lambda_a), G(\lambda_a) \cap K)\text{-modules} \to (\mathfrak{g}(\lambda), G(\lambda) \cap K)\text{-modules}$

(defined in terms of the parabolic subalgebra $\mathfrak{q}(\lambda_a) \cap \mathfrak{g}(\lambda)$). Similarly, we have the corresponding functors for compact group representations.

PROPOSITION 2.12 ([9, Th. 11.77].    *In the setting of* (2.3) *and* (2.4), *there are natural equivalences*

$$\mathcal{L}_{S_a}(\lambda_a) \simeq \mathcal{L}_S(\lambda) \circ \mathcal{L}_{S_a - S}(\lambda_a, \lambda), \qquad \mathcal{L}_{S_a}^K(\lambda_a) \simeq \mathcal{L}_S^K(\lambda) \circ \mathcal{L}_{S_a - S}^{G(\lambda) \cap K}(\lambda_a, \lambda).$$

*These equivalences respect the bottom layer maps of Proposition* 2.8 *and the Hermitian forms of Proposition* 2.10.

What is established in [9] is a composition-of-functors spectral sequence. The equivalence here appears as a corner of that spectral sequence, because the various cohomological induction functors all vanish above the degrees in the proposition. (A careful reader may notice that the assertions about bottom layer maps and Hermitian forms are not actually stated or proved in [9]. In any case they are stated here only for moral support; we will make no use of them.)

Using Proposition 2.12, it is a simple matter to extend Proposition 2.6 and Theorems 2.9 and 2.11 to apply to $\mathfrak{q}(\lambda)$ instead of $\mathfrak{q}(\lambda_a)$.

THEOREM 2.13.    *Suppose* $\lambda_a \in \Lambda_a$, *and there is fixed a singularization* $\lambda$ *of* $\lambda_a$ (*cf.* (2.3) *and* (2.4)).

a) *The functor* $\mathcal{L}_S^K(\lambda)$ *implements a bijection from* $B_a^{\lambda_a - \rho(\mathfrak{u}(\lambda))}(G(\lambda))$ *onto* $B_a^{\lambda_a}(G)$.

b) *There is a natural bijection*

$$\Pi_a^{\lambda_a - \rho(\mathfrak{u}(\lambda))}(G(\lambda)) \to \Pi_a^{\lambda_a}(G)$$

*defined as follows. Suppose that* $Z \in \Pi_a^{\lambda_a - \rho(\mathfrak{u}(\lambda))}(G(\lambda))$ *is an irreducible* $(\mathfrak{g}(\lambda), G(\lambda) \cap K)$-*module. Then the corresponding* $(\mathfrak{g}, K)$-*module* $X$ *is the unique irreducible subquotient of* $\mathcal{L}_S(\lambda)(Z)$ *containing* $K$-*types in* $B_a^{\lambda_a}(G)$. *Conversely, given* $X \in \Pi_a^{\lambda_a}(G)$, *put* $R = \dim(\mathfrak{u}(\lambda) \cap \mathfrak{p})$. *Then we assign to* $X$ *the unique irreducible subquotient* $Z$ *of* $H^*(\mathfrak{u}(\lambda), X)$ *containing* $G(\lambda) \cap K$-*types in* $B_a^{\lambda_a - \rho(\mathfrak{u}(\lambda))}(G(\lambda))$; *it appears in degree* $R$.

c) *The bijection of* (b) *preserves Hermitian representations.*

d) *Suppose* $Z \in \Pi_a^{\lambda - \rho(\mathfrak{u}(\lambda))}(G(\lambda))$ *is an irreducible Hermitian* $(\mathfrak{g}(\lambda), G(\lambda) \cap K)$-*module. Then the image of the bottom layer map* $\mathcal{B}_Z$ (*Prop.* 2.8) *maps* *injectively to the irreducible subquotient X described in* (b). *On these K-types, the signature of the Hermitian form on X therefore corresponds precisely to the signature on Z.*

e) *Suppose* $Z \in \Pi_a^{\lambda - \rho(\mathfrak{u}(\lambda))}(G(\lambda))$ *is an irreducible* $(\mathfrak{g}(\lambda), G(\lambda) \cap K)$-*module, and the infinitesimal character of Z corresponds by the Harish-Chandra isomorphism to a weight* $\phi \in \mathfrak{h}^*$. *Assume that*

$$\operatorname{Re} \langle \phi + \rho(\mathfrak{u}(\lambda)), \alpha \rangle \geq 0 \qquad (\alpha \in \Delta(\mathfrak{u}(\lambda), \mathfrak{h})).$$

*Then* $\mathcal{L}_S(\lambda)(Z)$ *is irreducible.*

f) *Under the assumptions of* (e), *suppose also that Z is unitary. Then* $\mathcal{L}_S(\lambda)(Z)$ *is unitary.*

g) *Under the assumptions of* (e), *suppose also that*

$$\langle \phi + \rho(\mathfrak{u}(\lambda)), \alpha \rangle \neq 0 \qquad (\alpha \in \Delta(\mathfrak{u}(\lambda), \mathfrak{h}))$$

*and that* $\mathcal{L}_S(\lambda)(Z)$ *is unitary. Then Z is unitary.*

*Sketch of proof.* For (a), Proposition 2.6 (applied to the groups $G$ and $G(\lambda)$) provides bijections

$$\mathcal{L}_{S_a}^K : B_a^{\lambda_a - \rho(\mathfrak{u}(\lambda_a))}(G(\lambda_a)) \to B_a^{\lambda_a}(G),$$
$$\mathcal{L}_{S_a - S}^{G(\lambda) \cap K} : B_a^{\lambda_a - \rho(\mathfrak{u}(\lambda_a))}(G(\lambda_a)) \to B_a^{\lambda_a - \rho(\mathfrak{u}(\lambda))}(G(\lambda)).$$

Now Proposition 2.12 implies that $\mathcal{L}_S^K(\lambda)$ provides a bijection between the ranges of these two maps. The proofs of (b) and (c) are similar; we omit the details. Part (d) is the Signature Theorem of [9, Cor. 11.228]; (e) is the Irreducibility Theorem of [9, Th. 8.7]; and (f) is the Unitarizability Theorem of [9, Th. 9.1] or [17, Th. 1.3(a)]. Part (g) is [17, Th. 1.3(b)].

### 3. Relation to the unitary classification

In this section we will relate the various sets defined for unitary representations in Definition 0.2 with the classification by lowest $K$-types described in Section 2. We begin with a $K$-dominant weight $\mu \in \widehat{T^c}$ as in (0.5), and a positive root system $\Delta^+(\mathfrak{g}, \mathfrak{t}^c)$ making $\mu + 2\rho_c$ dominant. Define $\rho$ as in (0.5b), half the sum of the positive roots. Recall from Corollary 2.4 that in the notation of Proposition 1.4,

$$(3.1a) \qquad \lambda_a(\mu) = T_\rho(\mu + 2\rho_c).$$

The definition in (0.5d) may be written as

$$(3.1b) \qquad \lambda_u(\mu) = T_{2\rho}(\mu + 2\rho_c).$$

By Corollary 1.6,

$$(3.1c) \qquad \lambda_u(\mu) = T_\rho \circ T_\rho(\mu + 2\rho_c) = T_\rho(\lambda_a(\mu)).$$

Lemma 1.3 (or Proposition 1.4) now implies that

$$(3.1d) \qquad \lambda_u(\mu) \text{ is a singularization of } \lambda_a(\mu) \text{ (notation 2.3)}.$$

This information leads immediately to a relationship between the sets of $K$-types defined as in Definitions 0.2 and 2.5.

PROPOSITION 3.1.  *Use the notation of* (0.5), (2.1), *and Proposition* 1.4.

a)  $\Lambda_u = T_\rho(\Lambda_a)$.
b)  $B_u^{\lambda_u}(G) = \bigcup_{\lambda_a \in \Lambda_a, \; T_\rho(\lambda_a) = \lambda_u} B_a^{\lambda_a}(G)$.
c)  *In the notation of Theorem* 2.13, *the functor* $\mathcal{L}_{S_u}^K(\lambda_u)$ *implements a bijection from* $B_u^{\lambda_u}(G(\lambda_u))$ *onto* $B_u^{\lambda_u}(G)$.

*Proof.* Parts (a) and (b) are immediate from (3.1c) and the definitions. For (c), we combine the description of $B_u^{\lambda_u}(G)$ in (b) with Theorem 2.13(a) (which applies because of (3.1d)). The conclusion is that $\mathcal{L}_{S_u}^K(\lambda_u)$ implements a bijection from

$$\bigcup_{\substack{\lambda_a \in \Lambda_a \\ T_\rho(\lambda_a) = \lambda_u}} B_a^{\lambda_a - \rho(\mathfrak{u}(\lambda_u))}(G(\lambda_u))$$

onto $B_u^{\lambda_u}(G)$. To complete the proof of (c), we need

LEMMA 3.2.  *In the setting of Proposition* 3.1, *suppose* $\lambda_a' \in \Lambda_a(G(\lambda_u))$. *Then the following two properties are equivalent.*

a)  $T_{\rho(\mathfrak{g}(\lambda_u))}(\lambda_a') = \lambda_u$.
b)  *There is a* $\lambda_a \in \Lambda_a(G)$ *such that* $T_\rho(\lambda_a) = \lambda_u$, *and* $\lambda_a' = \lambda_a - \rho(\mathfrak{u}(\lambda_u))$.

Assuming this lemma for a moment, we see that the last union may be rewritten as

$$\bigcup_{\substack{\lambda_a' \in \Lambda_a(G(\lambda_u)) \\ T_{\rho(\mathfrak{g}(\lambda_u))}(\lambda_a') = \lambda_u}} B_a^{\lambda_a'}(G(\lambda_u)).$$

According to Proposition 3.1(b) (which we have already proved) this is precisely $B_u^{\lambda_u}(G(\lambda_u))$, as we wished to show.  $\square$

*Proof of Lemma* 3.2. Suppose first that (b) holds. Fix a positive root system $\Delta^+(\mathfrak{g}, \mathfrak{t}^c)$ making $\lambda_a$ dominant, and write $\rho$ for half the sum of the positive roots. Then

$$(3.2a) \qquad \rho = \rho(\mathfrak{u}(\lambda_u)) + \rho(\mathfrak{g}(\lambda_u)).$$

According to Proposition 1.1(c), the hypothesis that $T_\rho(\lambda_a) = \lambda_u$ means that $\lambda_u$ is dominant for $\Delta^+(\mathfrak{g}, \mathfrak{t}^c)$, and that

(3.2b) $$\lambda_u = (\lambda_a - \rho(\mathfrak{u}(\lambda_u))) - \rho(\mathfrak{g}(\lambda_u)) + w,$$

with $w$ a nonnegative combination of positive roots orthogonal to $\lambda_u$. But this formulation also shows that if we set $\lambda'_a = \lambda_a - \rho(\mathfrak{u}(\lambda_u))$, then

$$T_{\rho(\mathfrak{g}(\lambda_u))}(\lambda'_a) = \lambda_u.$$

So (a) is satisfied.

Conversely, assume (a) holds. Choose a positive root system $\Delta^+(\mathfrak{g}(\lambda_u), \mathfrak{t}^c)$ making $\lambda'_a$ dominant. Define

(3.2c) $$\Delta^+(\mathfrak{g}, \mathfrak{t}^c) = \Delta^+(\mathfrak{g}(\lambda_u), \mathfrak{t}) \cup \Delta(\mathfrak{u}(\lambda_u), \mathfrak{t}^c).$$

Then (3.2a) still holds. Define

(3.2d) $$\lambda_a = \lambda'_a + \rho(\mathfrak{u}(\lambda_u)).$$

By Proposition 1.1(c), the assumption (a) may be written as in (3.2b). This looks at first as if it immediately implies (b) of the proposition. What is missing is this. To compute $T_\rho(\lambda_a)$, we must subtract from $\lambda_a$ half the sum of the positive roots for a system making $\lambda_a$ dominant. We do not yet know that $\lambda_a$ is dominant for $\Delta^+(\mathfrak{g}, \mathfrak{t}^c)$. Here is a proof of that. Since $\lambda_u$ is central in $G(\lambda_u)$, we can apply Proposition 1.7 almost directly to the hypothesis (a). The conclusion is

$\lambda'_a - \lambda_u$ belongs to the convex hull of the Weyl group orbit of $\rho(\mathfrak{g}(\lambda_u))$.

Of course the Weyl group in question is for $\Delta(\mathfrak{g}(\lambda_u), \mathfrak{t}^c)$. This Weyl group fixes $\rho(\mathfrak{u}(\lambda_u))$. We conclude that $\lambda_a$ may be written as a convex combination

$$\lambda_a = \sum c_w(\lambda_u + w\rho),$$

with $w$ running over the Weyl group of $\Delta(\mathfrak{g}(\lambda_u), \mathfrak{t}^c)$. Each summand here has positive inner product with the roots in $\Delta(\mathfrak{u}(\lambda_u), \mathfrak{t}^c)$, so $\lambda_a$ does as well. On the other hand it is immediate from (3.2d) that $\lambda_a$ is dominant for $\Delta^+(\mathfrak{g}(\lambda_u), \mathfrak{t})$; so we conclude that $\lambda_a$ is actually dominant for $\Delta^+(\mathfrak{g}, \mathfrak{t}^c)$, as we wished to show. Now (3.2b) does imply that $T_\rho(\lambda_a) = \lambda_u$, which is (b). $\qquad\square$

*Proof of Proposition* 0.3. Write $\delta$ and $\delta'$ for lowest $K$-types of $\pi$ having highest weights $\mu$ and $\mu'$. Define $\lambda_a = \lambda_a(\mu)$, and $\lambda'_a = \lambda_a(\mu')$. Use Theorem 2.9 and the $K$-type $\delta$ to write $\pi$ as a subquotient of some $\mathcal{L}_{S_a}(\lambda_a)(Z_a)$. Then $\delta'$ must be a lowest $K$-type of this cohomologically induced representation. By Theorem 2.9(b), $\delta'$ must also belong to $B_a^{\lambda_a}(G)$. It is clear from Definition 2.5 (and Lemma 1.4) that this implies that $\lambda'_a = r \cdot \lambda_a$, for some $r \in R(G)$. Consequently (3.1c) implies that

$$\lambda_u(\mu') = T_\rho(\lambda'_a) = T_\rho(r \cdot \lambda_a) = r \cdot T_\rho(\lambda_a) = r \cdot \lambda_u(\mu),$$

as we wished to show.                                                                 □

## 4. Relation to ordinary parabolic induction

In this section we consider the relationship between our results about cohomological induction and those formulated with ordinary parabolic induction. We begin with the case to which Theorem 2.9 reduces the classification of admissible representations.

PROPOSITION 4.1.    *Suppose $\lambda_a \in \Lambda_a$ (cf. (2.1)), and that $G(\lambda_a) = G$. Then $G$ is quasisplit; that is, there is a Borel subgroup $B^q = T^q A^q N^q$ of $G$. The Levi factor $H^q = T^q A^q$ is a maximally split Cartan subgroup of $G$. After replacing $B^q$ by a conjugate, we may assume that the identity component $T_0^q$ is contained in $T^c$ (cf. (0.3)). Then the linear functional $\lambda_a$ vanishes on the orthogonal complement of $\mathfrak{t}^q \subset \mathfrak{t}^c$, so it makes sense to write $\lambda_a \in i(\mathfrak{t}_0^q)^*$. An irreducible representation $\gamma \in \widehat{T}^q$ is said to be fine of type $\lambda_a$ if its differential is a multiple of $\lambda_a$.*

 a) *Suppose $\delta \in B_a^{\lambda_a}(G)$. Then the restriction of $\delta$ to $T^q$ consists of fine representations of $T^q$ of type $\lambda_a$, belonging to a single orbit of $W(G, H^q)$.*
 b) *The set $B_a^{\lambda_a}(G)$ consists of all lowest $K$-types of the induced representations $\mathrm{Ind}_{T^q}^K(\gamma)$, as $\gamma$ runs over fine representations of $T^q$ of type $\lambda_a$.*
 c) *Suppose $\delta \in B_a^{\lambda_a}(G)$, and $X \in \Pi_a^{\lambda_a}(G)$ is an irreducible admissible representation of lowest $K$-type $\delta$. Let $\gamma$ be any representation of $T^q$ appearing in the restriction of $\delta$. Then there is a character $\nu \in \widehat{A}^q$ so that $X$ is a Langlands subquotient of the principal series representation $\mathrm{Ind}_{T^q A^q N^q}^G(\gamma \otimes \nu)$. In particular, the infinitesimal character of $X$ has Harish-Chandra parameter $(\lambda_a, \nu) \in (\mathfrak{t}^q)^* + (\mathfrak{a}^q)^*$.*

This is [16, Prop. 5.3.26, Th. 4.3.16, and Th. 4.4.8]. Inspection of the proofs shows that what is claimed here does not use the assumption of abelian Cartan subgroups. (That is needed for example to show that the restriction of $\delta$ to $T^q$ decomposes with multiplicity one.)

We turn now to general representations. So let $\lambda_a \in \Lambda_a(G)$. Proposition 4.1 guarantees that $G(\lambda_a)$ is quasisplit. A little more precisely, Proposition 2.6 guarantees that

$$(4.1a) \qquad\qquad \lambda'_a = \lambda_a - \rho(\mathfrak{u}(\lambda_a))$$

belongs to $\Lambda_a(G(\lambda_a))$; and evidently

(4.1b) $$G(\lambda_a)(\lambda'_a) = G(\lambda_a),$$

so we can apply Proposition 4.1. Pick a quasi-split Cartan subgroup

(4.1c) $$H(\lambda_a) = T(\lambda_a)A(\lambda_a) \subset G(\lambda_a)$$

as in Proposition 4.1. Just as in the proposition, we may regard $\lambda'_a$, $\lambda_a$, and $\rho(\mathfrak{u}(\lambda_a))$ as weights in $it(\lambda_a)_0^*$. Now consider the centralizer of $A(\lambda_a)$ in $G$. Because $A(\lambda_a)$ is the vector part of the Cartan subgroup $H(\lambda_a)$, this centralizer has a Langlands decomposition

(4.1d) $$G^{A(\lambda_a)} = M(\lambda_a)A(\lambda_a).$$

Here $M(\lambda_a)$ is a $\theta$-stable reductive subgroup of $G$ containing $T(\lambda_a)$ as a compact Cartan subgroup.

PROPOSITION 4.2.    *In the setting* 4.1, *there are natural bijections among the following three sets.*

  a) *Fine representations of $T(\lambda_a)$ of type $\lambda'_a$ (Proposition 4.1); that is, irreducible representations of $T(\lambda_a)$ of differential equal to a multiple of $\lambda'_a$.*
  b) *Irreducible representations of $T(\lambda_a)$ of differential equal to a multiple of $\lambda_a - \rho(\mathfrak{m}(\lambda))$. Here the positive root system is the one for $\mathfrak{t}(\lambda_a)$ in $\mathfrak{m}(\lambda_a)$ making $\lambda_a$ dominant.*
  c) *Discrete series representations of $M(\lambda_a)$ of Harish-Chandra parameter $\lambda_a$.*

*Proof.* This result is analogous to Proposition 6.6.2 in [16]. To construct the bijection between (a) and (b), consider the set of roots $\Delta(\mathfrak{u}(\lambda_a), \mathfrak{h}(\lambda_a))$. The imaginary roots among these are precisely the positive roots of $\mathfrak{t}(\lambda_a)$ in $\mathfrak{m}(\lambda_a)$. The remaining roots occur in pairs $(\alpha, \theta\alpha)$. Choose one root from each of these pairs, and list the choices as $\{\alpha_1, \ldots, \alpha_r\}$. Define

$$\tau = \sum_i \alpha_i \in T(\lambda_a)\widehat{\phantom{.}}.$$

Because $\alpha$ and $\theta\alpha$ have the same restriction to $T(\lambda_a)$, the character $\tau$ is independent of choices. Furthermore

$$\rho(\mathfrak{u}(\lambda_a)) = \rho(\mathfrak{m}(\lambda_a)) + d\tau.$$

The bijection from (a) to (b) is just tensoring with $\tau$. The bijection from (b) to (c) is Harish-Chandra's parametrization of the discrete series. In the context of Section 2, the set in (c) turns out to be $\Pi_a^{\lambda_a}(M(\lambda_a))$, and the bijection with (b) is precisely 2.2; that is, it is given by cohomological induction from $T(\lambda_a)$ to $M(\lambda_a)$. Details may be found in [9, §11.8].    □

We can now formulate the relationship between Theorem 2.9 and the parabolic induction version of the Langlands classification.

THEOREM 4.3 ([9 , Th. 11.225]).   *Suppose we are in the setting of (2.2), and that the irreducible $(\mathfrak{g}, K)$-module $X$ corresponds to $Z_a$. Use the notation of (4.1), and fix a character $\gamma^q \otimes \nu$ of $H(\lambda_a)$ as in Proposition 4.1 so that $Z_a$ is a Langlands quotient of $\mathrm{Ind}_{H(\lambda_a)N^q}^{G(\lambda_a)}(\gamma \otimes \nu)$. Let $\gamma$ be the discrete series representation of $M(\lambda_a)$ with Harish-Chandra parameter $\lambda_a$ that corresponds to $\gamma_q$ in the bijection from* (a) *to* (c) *of Proposition 4.2. Choose a real parabolic subgroup $P = M(\lambda_a)A(\lambda_a)N$ of $G$ making $\mathrm{Re}\,\nu$ dominant for the $\mathfrak{a}(\lambda_a)$-weights in $\mathfrak{n}$; and define $N^q = N \cap G(\lambda_a)$. Then*

$$\mathcal{L}_{S_a}(\lambda_a)\left(\mathrm{Ind}_{H(\lambda_a)N^q}^{G(\lambda_a)}(\gamma \otimes \nu)\right) \simeq \mathrm{Ind}_P^G(\gamma \otimes \nu).$$

*In particular, $X$ is a Langlands quotient of $\mathrm{Ind}_P^G(\gamma \otimes \nu)$.*

This is stated without complete proof in [16, Th. 6.6.15]; it was essentially first proved in the unpublished second half of [15].

Recall that the tempered representations constitute the closure (in the Fell topology) of the support of the Plancherel measure. A tempered irreducible representation is unitary, but unitary representations need not be tempered.

COROLLARY 4.4.   *The bijection of Theorem 2.9 preserves tempered representations. Explicitly, suppose $Z_a \in \Pi_a^{\lambda_a - \rho(\mathfrak{u}(\lambda_a))}(G(\lambda_a))$ is an irreducible tempered $(\mathfrak{g}(\lambda_a), G(\lambda_a) \cap K)$-module. Then $\mathcal{L}_{S_a}(\lambda_a)(Z_a)$ is an irreducible tempered $(\mathfrak{g}, K)$-module in $\Pi_a^{\lambda_a}(G)$; and every such tempered representation arises in this way.*

*Proof.* The nature of the Langlands classification implies that a Langlands quotient is tempered if and only if the continuous parameter $\nu$ is purely imaginary (see for example [8, Th. 8.53]). Theorem 4.3 says that the bijection of Theorem 2.9 preserves continuous parameters, so temperedness is also preserved. For the rest, the infinitesimal character of $Z_a$ is given by the weight $(\lambda_a - \rho(\mathfrak{u}(\lambda_a)), \nu) \in \mathfrak{h}(\lambda_a)^*$. When $\nu$ is purely imaginary, the dominance of $\lambda_a$ implies immediately that the condition in Theorem 2.13(e) is satisfied; so $\mathcal{L}_{S_a}(\lambda_a)(Z_a)$ is irreducible.   $\square$

## 5. The reduction step in the classification

Theorem 2.13 and Proposition 3.1 together allow us to approach Conjecture 0.6. To see what we have, we consider first the case of Hermitian representations.

*Definition* 5.1.   Suppose $G$ is a reductive group in Harish-Chandra's class. The *Hermitian dual* of $G$ is the set $\Pi_h(G)$ of equivalence classes of irreducible $(\mathfrak{g}, K)$ modules endowed with nondegenerate invariant Hermitian forms. Here $K$ is a maximal compact subgroup of $G$; and an equivalence is required to preserve the Hermitian form.

Suppose $\lambda_u \in \Lambda_u$ (cf. (0.5)). The set of *Hermitian representations of $G$ attached to $\lambda_u$* is

$$\Pi_h^{\lambda_u}(G) = \{\pi \in \Pi_h(G) \mid \pi \text{ has a lowest } K\text{-type in } B_u^{\lambda_u}(G)\}$$

(cf. Definition 0.2).

Proposition 0.3 immediately shows that these sets partition the Hermitian dual:

PROPOSITION 5.2. *Suppose $\lambda_u$ and $\lambda_u'$ belong to $\Lambda_u$ (cf. (0.5)). If $\lambda_u$ and $\lambda_u'$ belong to the same orbit of $R(G)$, then $\Pi_h^{\lambda_u}(G) = \Pi_h^{\lambda_u'}(G)$. If not, then $\Pi_h^{\lambda_u}(G)$ and $\Pi_h^{\lambda_u'}(G)$ are disjoint.*

At this point it is useful to be a little more precise than we have been about the notion of signature.

*Definition* 5.3. (See [17, Def. 1.4], or [9, Prop. 6.12].) Suppose $X$ is an admissible representation of $K$ (that is, locally finite with finite multiplicities), and $\langle \, , \, \rangle$ is a $K$-invariant Hermitian form on $X$. The *signature of $\langle \, , \, \rangle$* is a set of three functions $(p, q, z)$ from $\widehat{K}$ to the nonnegative integers, defined as follows. Suppose $(\delta, V_\delta)$ is an irreducible representation of $K$. Fix a positive-definite invariant Hermitian form $\langle \, , \, \rangle_\delta$ on $V_\delta$. Define

$$X^\delta = \mathrm{Hom}_K(V_\delta, X).$$

This is a finite-dimensional vector space, of dimension equal to the multiplicity $m(\delta)$ of $\delta$ in $X$. There is a unique Hermitian form $\langle \, , \, \rangle^\delta$ on $X^\delta$ characterized by the property

$$\langle Tv, Sw \rangle = \langle T, S \rangle^\delta \langle v, w \rangle_\delta \qquad (T, S \in X^\delta, \ v, w \in V_\delta).$$

The form $\langle \, , \, \rangle^\delta$ may be diagonalized; and we define $p(\delta)$ (respectively $q(\delta), z(\delta)$) to be the number of positive (respectively negative, zero) diagonal entries. These numbers are of course independent of the choice of diagonalization, and characterize the form $\langle \, , \, \rangle^\delta$ up to equivalence. Notice that

$$m(\delta) = p(\delta) + q(\delta) + z(\delta).$$

The form $\langle \, , \, \rangle$ on $X$ is nondegenerate if and only if $z = 0$, and positive definite if and only if $q = z = 0$.

THEOREM 5.4. *Suppose $\lambda_u \in \Lambda_u$ (cf. (0.5)).*

a) *The functor $\mathcal{L}_{S_u}^K(\lambda_u)$ provides a bijection from $B_u^{\lambda_u}(G(\lambda_u))$ onto $B_u^{\lambda_u}(G)$.*

b) *The set* $\Pi_h^{\lambda_u}(G(\lambda_u))$ *consists of all irreducible* $(\mathfrak{g}(\lambda_u), G(\lambda_u) \cap K)$-*modules endowed with an invariant Hermitian form, and containing at least one* $(G(\lambda_u) \cap K)$-*type in the set* $B_u^{\lambda_u}(G(\lambda_u))$.

c) *There is a natural bijection*

$$\Pi_h^{\lambda_u}(G(\lambda)) \to \Pi_h^{\lambda_u}(G)$$

*defined as follows. Suppose that* $Z \in \Pi_h^{\lambda_u}(G(\lambda_u))$ *is irreducible as a* $(\mathfrak{g}(\lambda_u), G(\lambda_u) \cap K)$-*module. Then the corresponding* $(\mathfrak{g}, K)$-*module* $X$ *is the unique irreducible subquotient of* $\mathcal{L}_{S_u}(\lambda_u)(Z)$ *containing* $K$-*types in* $B_a^{\lambda_a}(G)$. *The Hermitian form on* $X$ *is the one inherited from the one on* $\mathcal{L}_{S_u}(\lambda_u)(Z)$ *induced from* $Z$ *by Proposition 2.10.*

d) *The bijections of* (a) *and* (c) *preserve signatures of Hermitian forms; that is, the signature on a* $K$ *type in* $B_u^{\lambda_u}(G)$ *of the Hermitian form on* $X$ *is equal to the signature on the corresponding* $G(\lambda_u) \cap K$-*type of the Hermitian form on* $Z$.

e) *Suppose* $Z \in \Pi_h^{\lambda_u}(G(\lambda))$ *is an irreducible* $(\mathfrak{g}(\lambda_u), G(\lambda) \cap K)$-*module, and the infinitesimal character of* $Z$ *corresponds by the Harish-Chandra isomorphism to a weight* $\phi \in \mathfrak{h}^*$. *Assume that*

$$\operatorname{Re} \langle \phi + \rho(\mathfrak{u}(\lambda_u)), \alpha \rangle \geq 0 \qquad (\alpha \in \Delta(\mathfrak{u}(\lambda_u), \mathfrak{h})).$$

*Then* $\mathcal{L}_{S_u}(\lambda_u)(Z)$ *is irreducible.*

*Proof.* Part (a) is Proposition 3.1(c). For (b), the group $G$ plays no role, so we may as well assume $G = G(\lambda_u)$. One inclusion is obvious; so suppose $X$ is an irreducible $(\mathfrak{g}, K)$-module containing a $K$-type $\delta \in B_u^{\lambda_u}(G)$. We must show that the lowest $K$-type $\delta'$ also belongs to $B_u^{\lambda_u}(G)$. Write $Z$ for the identity component of the compact part of the center of $G$. Then $\mathfrak{t}^c$ is the orthogonal direct sum of $\mathfrak{z}$ and the span of the roots. Since $G(\lambda_u) = G$, the weight $\lambda_u$ is orthogonal to all the roots; so $\lambda_u \in i\mathfrak{z}_0^*$ is just the differential of the character by which $Z$ acts on $\delta$. Since $Z$ acts by the same character on $\delta'$, we see that $\lambda_u'|_{\mathfrak{z}} = \lambda_u$. This implies that

$$(5.1) \qquad\qquad |\lambda_u| \leq |\lambda_u'|.$$

Now we apply

LEMMA 5.5. *Suppose* $X$ *is an irreducible* $(\mathfrak{g}, K)$-*module of lowest* $K$-*type* $\delta' \in B_u^{\lambda_u'}(G)$, *and* $\delta$ *is any other* $K$-*type of* $X$. *Then* $\delta$ *has a highest weight* $\mu$ *with*

$$\lambda_u(\mu) = \lambda_u' + Q',$$

*with* $Q'$ *a sum with nonnegative coefficients of roots of* $\mathfrak{t}^c$ *in* $\mathfrak{q}(\lambda_u')$. *In particular,*

$$\langle \lambda_u(\mu), \lambda_u(\mu) \rangle \geq \langle \lambda_u', \lambda_u' \rangle,$$

*with equality only if* $\lambda_u(\mu) = \lambda_u'$.

We will give the proof in a moment. Assuming the result, we find that $|\lambda_u| \geq |\lambda'_u|$. By 5.1, equality must hold. The condition for equality in the Lemma is $\lambda'_u = \lambda_u$; so the lowest $K$-type $\delta'$ of $X$ belongs to $B^{\lambda_u}_u(G)$, as we wished to show.

Part (c) is (in light of (3.1)) just Theorem 2.13(b); part (d) is Theorem 2.13(d); and part (e) is Theorem 2.13(e). □

*Proof of Lemma* 5.5. This is a variation on Lemma 6.5.6 in [16], which we will also use here. Fix a highest weight $\mu'$ of $\delta'$ so that $\lambda_u(\mu') = \lambda'_u$, and define $\lambda'_a = \lambda_a(\mu')$. According to the proofs of Lemma 6.5.6 and Theorem 6.5.9 in [16] we can find a highest weight $\mu$ of $\delta$ so that

$$(5.2a) \qquad \lambda_a(\mu) = \lambda'_a + Q'_1,$$

with $Q'$ a sum with nonnegative coefficients of roots of $\mathfrak{t}^c$ in $\mathfrak{q}(\lambda'_u)$. We want a relationship between $\lambda_u(\mu)$ and $\lambda'_u$, so we rewrite both sides in those terms. Fix a positive root system $(\Delta^+)'$ making $\lambda'_a$ dominant, and write $\rho'$ for half the sum of its positive roots. Then 5.2(a) becomes

$$(5.2b) \qquad \lambda_u(\mu) + \rho_0 = \lambda'_u + Q'_1 + \rho' - w'.$$

Here $w'$ is a nonnegative combination of roots in $(\Delta^+)'$ orthogonal to $\lambda'_u$. On the left we have used instead Corollary 1.9: $\rho_0$ is a dominant weight in the convex hull of the $W_0$-orbit of $\rho$. By Proposition 1.10, we can write

$$(5.2c) \qquad \rho_0 = \rho' - \sum_{\alpha' \in (\Delta^+)'} b_{\alpha'} \alpha' \qquad (0 \leq b_{\alpha'} \leq 1).$$

Inserting this in (5.2b) gives

$$(5.2d) \qquad \lambda_u(\mu) = \lambda'_u + Q'_1 - w' + \sum_{\alpha' \in (\Delta^+)'} b_{\alpha'} \alpha'.$$

This has the form required in the lemma. The last inequality is an immediate consequence. □

We need one more definition.

*Definition* 5.6 (16, Def. 5.4.11]). Suppose $\mathfrak{h}_0$ is a $\theta$-stable Cartan subalgebra of $\mathfrak{g}_0$. Write $\mathfrak{h}_0 = \mathfrak{t}_0 + \mathfrak{a}_0$ for its decomposition into the $+1$ and $-1$ eigenspaces of $\theta$: $\mathfrak{t}_0 = \mathfrak{h}_0 \cap \mathfrak{k}_0$, and $\mathfrak{a}_0 = \mathfrak{h}_0 \cap \mathfrak{p}_0$. Any complex-valued linear functional $\gamma \in \mathfrak{h}^*$ may then be written as

$$\gamma = (\lambda, \nu) \in \mathfrak{t}^* + \mathfrak{a}^*.$$

Now make a further decomposition of $\gamma$ as

$$\gamma = \gamma_1 + i\gamma_2 = (\lambda_1, \nu_1) + i(\lambda_2, \nu_2),$$

in such a way that $\gamma_j$ takes purely real values on $\mathfrak{a}_0$, and purely imaginary values on $\mathfrak{t}_0$. That is,

$$\lambda_j \in i\mathfrak{t}_0^*, \qquad \nu_j \in \mathfrak{a}_0^*, \qquad \gamma_j \in i\mathfrak{t}_0^* + \mathfrak{a}_0^*.$$

Define the *canonical real and imaginary parts of $\gamma$* as

$$\mathrm{RE}\,\gamma = \gamma_1, \qquad \mathrm{IM}\,\gamma = \gamma_2.$$

If $\alpha \in \mathfrak{h}^*$ is a root, then

$$\mathrm{Re}\,\langle \gamma, \alpha \rangle = \langle \mathrm{RE}\,\gamma, \alpha \rangle,$$

and similarly for imaginary parts.

One checks easily that the decomposition into canonical real and imaginary parts commutes with the action of the Weyl group, and in fact with any automorphism that is inner for the complexified Lie algebra ([16, Lemma 5.4.12]).

Here is the main conjecture of this paper.

*Conjecture* 5.7.   Suppose $\lambda_u \in \Lambda_u$ (cf. (0.5)), and that $G(\lambda_u) = G$. Suppose $X \in \Pi_h^{\lambda_u}(G)$ is an irreducible Hermitian $(\mathfrak{g}, K)$-module (Definition 5.1). Let $\mathfrak{h}$ be a $\theta$-stable Cartan subalgebra of $\mathfrak{g}_0$, and $\phi \in \mathfrak{h}^*$ a weight parametrizing the infinitesimal character of $X$ by the Harish-Chandra isomorphism. Assume that the canonical real part $\mathrm{RE}\,\phi$ does *not* belong to $\lambda_u$ plus the convex hull of the Weyl group orbit of $\rho$. Then the signature of the Hermitian form on $X$ must be indefinite on $K$-types in $B_u^{\lambda_u}(G)$.

The statement of this conjecture is complicated somewhat by the central weight $\lambda_u$. Of course the center of $\mathfrak{g}$ is more or less irrelevant to interesting statements about unitary representations. In particular, Conjecture 5.7 is equivalent to the following statement about semisimple groups. (This kind of reduction is discussed again in 6.1 below.)

*Conjecture* 5.7′.   Suppose $G$ is a connected semisimple group in Harish-Chandra's class; and suppose $X \in \Pi_h^0(G)$ is an irreducible Hermitian $(\mathfrak{g}, K)$-module (Definition 5.1). Let $\mathfrak{h}$ be a $\theta$-stable Cartan subalgebra of $\mathfrak{g}_0$, and $\phi \in \mathfrak{h}^*$ a weight parametrizing the infinitesimal character of $X$ by the Harish-Chandra isomorphism. Assume that the canonical real part $\mathrm{RE}\,\phi$ does *not* belong to the convex hull of the Weyl group orbit of $\rho$. Then the signature of the Hermitian form on $X$ must be indefinite on $K$-types in $B_u^0(G)$.

We will give some supporting evidence for this conjecture in Proposition 7.12 below, which is based on Parthasarathy's Dirac operator inequality. We also formulate a conjectural sharpening of the Dirac operator inequality which would imply Conjecture 5.7′ (Conjecture 7.19).

THEOREM 5.8. *In the setting of Theorem 5.4, assume that Conjecture 5.7 is true for $G(\lambda_u)$. Then the bijection of (c) preserves unitarity. More precisely, suppose $\mathfrak{h}_0$ is a $\theta$-stable Cartan subalgebra of $\mathfrak{g}(\lambda_u)_0$. Write $W_0$ for the Weyl group of $\mathfrak{h}$ in $\mathfrak{g}(\lambda_u)$.*

a) *The cohomological parabolic induction functor $\mathcal{L}_{S_u}(\lambda_u)$ implements a bijection from $\Pi_u^{\lambda_u}(G(\lambda_u))$ onto $\Pi_u^{\lambda_u}(G)$.*

b) *Assume that $Z$ is an irreducible unitary $(\mathfrak{g}(\lambda_u), G(\lambda_u) \cap K)$-module containing a $(G(\lambda_u) \cap K)$-type in $B_u^{\lambda_u}(G(\lambda_u))$, and that the infinitesimal character of $Z$ corresponds by the Harish-Chandra isomorphism to a weight $\phi(Z) \in \mathfrak{h}^*$. Then $\mathrm{RE}\,\phi(Z) - \lambda_u$ (cf. Definition 5.6) belongs to the convex hull of the $W_0$ orbit of $\rho(\mathfrak{g}(\lambda_u))$.*

c) *Assume that $X \in \Pi_u^{\lambda_u}(G)$, and that the infinitesimal character of $X$ corresponds by the Harish-Chandra isomorphism to a weight $\phi(X) \in \mathfrak{h}^*$. Assume that $\mathrm{RE}\,\phi(X)$ is dominant for the roots of $\mathfrak{h}$ in $\mathfrak{u}(\lambda_u)$. Then $\mathrm{RE}\,\phi(X) - (\lambda_u + \rho(\mathfrak{u}(\lambda_u)))$ belongs to the convex hull of the $W_0$ orbit of $\rho(\mathfrak{g}(\lambda_u))$. In particular, $T_\rho(\mathrm{RE}\,\phi(X)) = \lambda_u$.*

The last statement of (c) appears to recover $\lambda_u$ from the infinitesimal character of any representation in $\Pi_u^{\lambda_u}(G)$. It does not quite do that, however, since we needed the weight $\phi(X)$ to satisfy a positivity condition determined by $\lambda_u$. At any rate we can recover the Weyl group orbit of $\lambda_u$.

Before embarking on the proof of this result, we need a complement to Conjecture 5.7.

THEOREM 5.9. *Suppose $\lambda_u \in \Lambda_u$ (cf. (0.5)), and that $G(\lambda_u) = G$. Suppose $X \in \Pi_h^{\lambda_u}(G)$ is an irreducible Hermitian $(\mathfrak{g}, K)$-module. Let $\mathfrak{h}$ be a $\theta$-stable Cartan subalgebra of $\mathfrak{g}_0$, and $\phi \in \mathfrak{h}^*$ a weight parametrizing the infinitesimal character of $X$ by the Harish-Chandra isomorphism. Assume that the canonical real part $\mathrm{RE}\,\phi$ belongs to $\lambda_u$ plus the convex hull of the Weyl group orbit of $\rho$. Then $X$ is unitary if and only if its Hermitian form is positive definite on $K$-types in $B_u^{\lambda_u}(G)$.*

*Proof.* There is no loss of generality in assuming that $\mathfrak{h}$ is the fundamental Cartan subalgebra. We use Theorem 1.5 of [17]. That provides irreducible tempered $(\mathfrak{g}, K)$-modules $Z_1, \ldots, Z_p$ (of real infinitesimal character) and integers $r_1^\pm, \ldots, r_p^\pm$, so that the signature of the Hermitian form on $X$ is

$$\left( \sum r_i^+ \Theta_K(Z_i), \sum r_i^- \Theta_K(Z_i) \right).$$

Here $\Theta_K(Z)$ denotes the formal $K$-character of $Z$. Fix a lowest $K$-type $\delta_i$ of $Z_i$, and a highest weight $\mu_i$ of $\delta_i$. Write $\lambda_i = \lambda_a(\mu_i)$.

Because $Z_i$ is tempered and has real infinitesimal character, Theorem 4.3 shows that the infinitesimal character is given by the weight $\lambda_i$. The proof of Theorem 1.5 in [17] shows first of all that each $Z_i$ has the same central character as $X$; and second, that the infinitesimal character $\lambda_i$ must belong to the convex hull of the Weyl group orbit of $\mathrm{RE}\,\phi$. From these two facts we deduce that $\lambda_i$ must belong to $\lambda_u$ plus the convex hull of the Weyl group orbit of $\rho$. Corollary 1.9 therefore implies that (in the notation of Proposition 1.4) $T_\rho(\lambda_i) = \lambda_u$. By Corollary 1.6, it follows that $\lambda_u(\mu_i) = \lambda_u$, and therefore that $\delta_i \in B_u^{\lambda_u}(G)$.

For the proof of the theorem, "only if" is clear; so suppose $X$ is not positive definite. This means that $\sum r_i^- \Theta_K(Z_i) \neq 0$. It follows that the lowest $K$-types of one of the $Z_i$ must appear with nonzero multiplicity in this sum; that is, that the form on $X$ is not positive on some $\delta_i$. (For this it is enough to take $|\lambda_i|$ minimal with $r_i^- \neq 0$. The absence of cancellation then follows from Theorem 6.5.9(b) in [16].) But we have just seen that $\delta_i$ belongs to $B_u^{\lambda_u}(G)$, as we wished to show.          $\square$

*Proof of Theorem* 5.8. Part (b) is an immediate consequence of Conjecture 5.7. For (a), suppose $X \in \Pi_u^{\lambda_u}(G)$; write $Z$ for the corresponding Hermitian representation in $\Pi_h^{\lambda_u}(G(\lambda_u))$ (Theorem 5.4(c)). By Theorem 5.4(d), the Hermitian form on $Z$ must be positive on the $(G(\lambda_u) \cap K)$-types in $B_u^{\lambda_u}(G(\lambda_u))$. By Conjecture 5.7, the infinitesimal character of $Z$ corresponds to a weight $\phi(Z) \in \mathfrak{h}^*$ with

$$\mathrm{RE}\,\phi(Z) = \lambda_u + w_0. \tag{5.3a}$$

Here $w_0$ is a weight in the convex hull of the $W_0$ orbit of $\rho(\mathfrak{g}(\lambda_u))$. By Theorem 5.9, $Z$ is unitary. Adding $\rho(\mathfrak{u}(\lambda_u))$ to (5.3a) gives

$$\mathrm{RE}\,\phi(Z) + \rho(\mathfrak{u}(\lambda_u)) = \lambda_u + \rho(\mathfrak{u}(\lambda_u)) + w_0. \tag{5.3b}$$

The right side here is strictly dominant on roots of $\mathfrak{h}$ in $\mathfrak{u}(\lambda_u)$ (compare the argument at the end of the proof of Lemma 3.2). By Theorem 5.4(e), $\mathcal{L}_{S_u}(\lambda_u)(Z)$ is irreducible. For the converse, suppose $Z \in \Pi_u^{\lambda_u}(G(\lambda_u))$. By Conjecture 5.7, the infinitesimal character of $Z$ satisfies (5.3a). We have seen that this implies the positivity condition in Theorem 2.13(e); so by Theorem 2.13(e) and (f), $\mathcal{L}_{S_u}(\lambda_u)(Z)$ is irreducible and unitary, as we wished to show.

For (c), we may assume by (a) that $X = \mathcal{L}_{S_u}(\lambda_u)(Z)$, and that the infinitesimal character of $Z$ corresponds to a weight $\phi(Z)$ satisfying (5.3a). Using Corollary 5.25 of [9] (or Prop. 6.3.11 of [16]) we find that the infinitesimal character of $X$ is represented by the weight

$$\phi(X) = \phi(Z) + \rho(\mathfrak{u}(\lambda_u)). \tag{5.3c}$$

We have seen above that $\mathrm{RE}\,\phi(X)$ is strictly dominant for the roots of $\mathfrak{h}$ in $\mathfrak{u}(\lambda_u)$; so this $\phi(X)$ must coincide (up to the action of $W_0$) with the representative chosen in (c). Now (5.3c) and (5.3b) together imply that $\phi(X)$ has the property asserted in (c). This in turn implies that $T_\rho(\operatorname{RE}\phi(X)) = \lambda_u$ by (for example) Corollary 1.9. □

COROLLARY 5.10. *Conjecture 5.7 (for the subgroup $G(\lambda_u)$ of $G$) implies Conjecture 0.6.*

## 6. The set $B_u^z(G)$ of unitarily small $K$-types

*Definition* 6.1. Suppose $G$ is a real reductive group in Harish-Chandra's class, and $K$ is a maximal compact subgroup. A parameter $\lambda_u \in \Lambda_u$ (cf. (0.5)) is called *unitarily small* if $G(\lambda_u) = G$ (cf. (0.6)). The set of unitarily small parameters is written $\Lambda_u^z(G)$; the $z$ stands for central. A highest weight $\mu \in \widehat{T}^c$ is called *unitarily small* if $\lambda_u(\mu)$ is unitarily small. A representation $\delta \in \widehat{K}$ is called *unitarily small* (with respect to $G$) if $\delta$ has a unitarily small highest weight. The set of all unitarily small $K$-types is written $B_u^z(G)$.

This is the analogue of "small" ([16, Def. 5.3.24], or [15, Def. 5.1]); the change is that the weight $\lambda_a(\mu)$ constructed in Proposition 2.3 has been replaced by $\lambda_u(\mu)$. Because of (3.1d), unitarily small is a weaker condition than small. It is fairly easy to see from the definition and Proposition 1.1 that the trivial representation of $K$ is always unitarily small. By contrast, the trivial representation is small in the sense of [16] if and only if $G$ is quasisplit.

Theorem 5.8 reduces the classification of unitary representations (at least modulo Conjecture 5.7) to the classification of sets $\Pi_u^{\lambda_u}(G)$, with $\lambda_u$ unitarily small. By Theorem 5.4(b), membership in such a set is equivalent to containing a $K$-type in $B_u^{\lambda_u}(G)$. So Conjecture 5.7 reduces the classification of unitary representations to the classification of those containing a unitarily small $K$-type. That is one reason to understand such $K$-types as completely as possible. A second reason is that Conjecture 5.7 is again a statement about representations containing unitarily small $K$-types.

For these reasons, we devote the next two sections to a series of characterizations of unitarily small $K$-types (Theorem 6.7, Corollary 6.12, Proposition 7.1, and Proposition 7.17). We begin with some examples of the calculation in (0.5).

*Example* 6.2. Suppose $G = \operatorname{SL}(2,\mathbb{R})$, and $K = \operatorname{SO}(2) = T^c$, so that

(a) $$\widehat{K} = \widehat{T}^c \simeq \mathbb{Z} \subset \mathbb{R} \simeq i(\mathfrak{t}_0^c)^*.$$

The set of roots of $T^c$ in $\mathfrak{g}$ is

(b) $$\Delta(\mathfrak{g}, \mathfrak{t}^c) = \{\pm 2\};$$

both of these are noncompact imaginary. Write $\mu_n$ for the irreducible representation of $K$ corresponding to the integer $n$. Because there are no compact roots, $2\rho_c = 0$. Therefore $\mu_n + 2\rho_c = \mu_n$ is dominant for the positive root system

(c)                               $\Delta^+(\mathfrak{g}, \mathfrak{t}^c) = \{(\operatorname{sgn} n) \cdot 2\};$

here we define $\operatorname{sgn} 0 = +1$. The corresponding positive Weyl chamber $C$ consists of the half line containing $\operatorname{sgn} n$. The projection $P$ on $C$ is the identity on $C$, and sends the other half line to zero. Now $2\rho = 2\operatorname{sgn} n$, so

(d)                               $\mu_n + 2\rho_c - 2\rho = n - 2\operatorname{sgn} n.$

This belongs to $C$ (and so is equal to $\lambda_u(\mu_n)$) as long as $|n| \geq 2$. For smaller values of $n$ it does not belong to $C$, so the projection on $C$ is zero. To summarize:

(e)                 $\lambda_u(\mu_n) = \begin{cases} n - 2\operatorname{sgn} n, & \text{if } |n| \geq 2; \\ 0, & \text{if } |n| \leq 2. \end{cases}$

This verifies the assertions made in the introduction after (0.1). In particular, we see that the set of unitarily small $K$-types (or highest weights) is

(f)                               $B_u^z(\mathrm{SL}(2, \mathbb{R})) = \{\mu_n \mid |n| \leq 2\}.$

By contrast, the set of small $K$-types in the sense of [16] consists of the $\mu_n$ with $|n| \leq 1$.

*Example* 6.3.   Suppose $G = \mathrm{Sp}(4, \mathbb{R})$, the real symplectic group of rank 2. Then $K = \mathrm{U}(2)$ (the rank two unitary group), so $T = \mathrm{U}(1) \times \mathrm{U}(1)$, and

(a)                               $\widehat{T}^c \simeq \mathbb{Z}^2 \subset \mathbb{R}^2 \simeq i(\mathfrak{t}_0^c)^*.$

The set of roots of $T^c$ in $\mathfrak{g}$ is

(b)              $\Delta(\mathfrak{g}, \mathfrak{t}^c) = \{(\pm 2, 0), (0, \pm 2), \pm(1, 1), \pm(1, -1)\};$

all but the last pair are noncompact. As a positive root in $K$ we choose $(1, -1) = 2\rho_c$; the dominant weights $\mu_{(p,q)}$ are then parametrized by decreasing pairs of integers $(p, q)$. We write $\delta_{(p,q)}$ for the irreducible representation of $\mathrm{U}(2)$ of highest weight $\mu_{(p,q)}$. Calculation of $\lambda_u(\mu_{(p,q)})$ falls into four cases according to the positive root system defined by $\mu_{(p,q)} + 2\rho_c$. We look at just one of these cases carefully. Suppose that

(c)
$p + 1 \geq 1 - q \geq 0$,   $\Delta^+(\mathfrak{g}, \mathfrak{t}^c) = \{(2, 0), (0, -2), (1, 1), (1, -1)\}$,   $2\rho = (4, -2)$.

Then one can calculate

(d)    $\lambda_u(\mu_{(p,q)}) = \begin{cases} (p - 3, q + 1), & \text{if } p - 3 \geq -1 - q \geq 0; \\ (p - 3, 0), & \text{if } p - 3 \geq 0 \text{ and } -1 - q \leq 0; \\ (\frac{p-q-4}{2}, \frac{-(p-q-4)}{2}), & \text{if } p - q \geq 4 \text{ and } p - 3 \leq -1 - q; \\ (0, 0) & \text{if } p - 3 \leq 0 \text{ and } p - q \leq 4. \end{cases}$

Combining this information with similar calculations in the other three cases, we find

(e) $$B_u^z(\mathrm{Sp}(4,\mathbb{R})) = \{\delta_{(p,q)} \mid 3 \geq p \geq q \geq -3, \ p - q \leq 4\}.$$

There are 25 representations of $K$ in this set. The set of small representations of $K$ consists of just five of them:

(f) $$\{\delta_{(p,q)} \mid 1 \geq p \geq q \geq -1, \ p - q \leq 1\}.$$

*Example* 6.4. In order to understand the role of the compact part of the center of $G$, we include one example where that is nontrivial. So let $G = \mathrm{U}(1,1)$, the unitary group of the standard indefinite Hermitian form on $\mathbb{C}^2$. We can take $K = T^c = \mathrm{U}(1) \times \mathrm{U}(1)$, and

(a) $$\widehat{T}^c \simeq \mathbb{Z}^2 \subset \mathbb{R}^2 \simeq i(\mathfrak{t}_0^c)^*.$$

The set of roots of $T^c$ in $\mathfrak{g}$ is

(b) $$\Delta(\mathfrak{g}, \mathfrak{t}^c) = \{\pm(1,-1)\};$$

these roots are noncompact. Of course this example is very close to Example 6.2, and we use parallel notation. A calculation like the one given there shows that

(c) $$\lambda_u(\mu_{(p,q)}) = \begin{cases} (p,q) - (\mathrm{sgn}\,(p-q))(1,-1), & \text{if } |p-q| \geq 2; \\ ((p+q)/2, (p+q)/2), & \text{if } |p-q| \leq 2. \end{cases}$$

The set of unitarily small parameters is therefore

(d) $$\Lambda_u^z(\mathrm{U}(1,1)) = \{(m/2, m/2) \mid m \in \mathbb{Z}\}.$$

Notice that these weights are precisely the restrictions to the center of arbitrary weights in $\widehat{T}^c$. The unitarily small $K$-types are

(e) $$B_u^z(\mathrm{U}(1,1)) = \{\delta_{(p,q)} \mid |p-q| \leq 2\}.$$

With these examples in hand, we turn to a general description of unitarily small $K$-types. In addition to the notation from the introduction, we need to introduce

(6.1a) $$Z = \text{identity component of the center of } G.$$

This group is preserved by $\theta$, and so has a direct product decomposition

(6.1b) $$Z = Z_c Z_h \qquad (Z_c = Z \cap K, \ Z_h = \exp(\mathfrak{z}_0 \cap \mathfrak{p}_0)).$$

(The subscript $h$ stands for "hyperbolic.") Write

(6.1c) $$G_s = \text{derived group of } G_0, \quad T_s = G_s \cap T^c.$$

On the Lie algebra level we have direct sum decompositions

(6.1d) $$\mathfrak{g}_0 = \mathfrak{g}_{s,0} + \mathfrak{z}_0, \quad \mathfrak{t}_0^c = \mathfrak{t}_{s,0} + \mathfrak{z}_{c,0}.$$

(The subspaces $\mathfrak{t}_{s,0}$ and $\mathfrak{z}_{c,0}$ correspond to what was called $V_s$ and $V_z$ in (1.2).) On the group level the corresponding decompositions are not direct: $Z \cap G_s$ can be nontrivial. (Because we have taken $G_s$ to be connected, the multiplication map from $G_s \times Z$ to $G$ is never surjective unless $G$ is connected. In fact this problem persists even if we replace $G_s$ by the full derived group of $G$.)

We propose to use the decompositions in (6.1) to reduce many problems to the semisimple case. There is a technical complication, however: the subgroup $G_s$ need not be closed in $G$. It is preserved by the Cartan involution, and has a well-behaved Cartan decomposition

$$G_s = (G_s \cap K) \cdot \exp(\mathfrak{g}_{s,0} \cap \mathfrak{p}_0);$$

but the (connected reductive) group $G_s \cap K$ need not be compact. When this happens, $G_s$ is not in Harish-Chandra's class. Nevertheless all the results we have proved apply to it, after minor adjustments in the formulations. (For example, one considers only unitary representations of $G_s \cap K$, and only unitary characters of $T_s$.) Alternatively, one can strengthen the hypotheses on $G$ to exclude the problem. For example, one can assume that $G$ is a finite cover of a real reductive algebraic group. This assumption is inherited by $G_s$. We will therefore ignore the problem, and speak as if $G_s$ were again in Harish-Chandra's class.

LEMMA 6.5.    *With notation as above, suppose $\mu \in \widehat{T}^c$ is a $K$-dominant weight. Write $\mu_z$ and $\mu_s$ for its restrictions to $Z_c$ and $T_s$. Then*

$$\lambda_u(G)(\mu) = \lambda_u(G_s)(\mu_s) + \mu_z \in i\mathfrak{t}_{s,0}^* + i\mathfrak{z}_{c,0}^*.$$

*In particular, $\mu$ is unitarily small for $G$ if and only if $\mu_s$ is unitarily small for $G_s$.*

This is obvious from the definitions. In the same way we can reduce matters to the various simple factors of $G_s$.

LEMMA 6.6.    *With notation as above, list the simple factors of $G_s$ as $G_s^1, \ldots G_s^r$. Put $T_s^i = T_s \cap G_s^i$, so that there is a direct sum decomposition*

$$\mathfrak{t}_s^* = (\mathfrak{t}_s^1)^* + \cdots + (\mathfrak{t}_s^r)^*.$$

*Suppose $\mu_s \in \widehat{T}_s$ is a $(G_s \cap K)$-dominant weight. Write $\mu_s^i$ for its restriction to $T_s^i$. Then*

$$\lambda_u(G_s)(\mu_s) = \lambda_u(G_s^1)(\mu_s^1) + \cdots + \lambda_u(G_s^r)(\mu_s^r).$$

*In particular, $\mu_s$ is unitarily small if and only if all the $\mu_s^i$ are.*

Again we omit the simple proof.

Here is the main result of this section.

THEOREM 6.7.    *Use the notation of Definition 6.1 and (6.1). Suppose*
$\delta \in \widehat{K}$ *has highest weight* $\mu \in \widehat{T}^c$. *Write* $\mu_z \in i\mathfrak{z}_{c,0}^*$ *for the differential of the*
*restriction of* $\mu$ *to* $Z_c$. *Define* $\lambda_u(\mu)$ *and* $\lambda_a(\mu)$ *as in (0.5) and Proposition 2.3.*
*Then the following conditions are equivalent.*

a) *The* $K$-*type* $\delta$ *is unitarily small* (Definition 6.1).

b) *The parameter* $\lambda_u(\mu)$ *is equal to* $\mu_z$.

c) *The parameter* $\lambda_a(\mu)$ *belongs to* $\mu_z$ *plus the convex hull of the Weyl group*
*orbit of* $\rho$.

d) *Let* $\Delta^+(\mathfrak{g}, \mathfrak{t}^c)$ *be a positive root system making* $\mu + 2\rho_c$ *dominant. List the*
*fundamamental weights for* $\Delta^+$ *as* $\{\xi_1, \ldots, \xi_l\}$ (cf. (1.2)). *Write* $2\rho_n$ *for*
*the sum of the noncompact positive roots. Then*

$$\langle \xi_i, \mu \rangle \leq \langle \xi_i, 2\rho_n \rangle \qquad (1 \leq i \leq l).$$

e) *Let* $(\Delta^+)'(\mathfrak{g}, \mathfrak{t}^c)$ *be any positive root system containing* $\Delta^+(\mathfrak{k}, \mathfrak{t}^c)$. *List the*
*fundamamental weights for* $(\Delta^+)'$ *as* $\{\xi_1', \ldots, \xi_l'\}$ (cf. (1.2)). *Write* $2\rho_n'$ *for*
*the sum of the noncompact positive roots. Then*

$$\langle \xi_i', \mu \rangle \leq \langle \xi_i', 2\rho_n' \rangle \qquad (1 \leq i \leq l).$$

f) *There is a system of positive roots* $(\Delta^+)'(\mathfrak{g}, \mathfrak{t}^c)$ *containing* $\Delta^+(\mathfrak{k}, \mathfrak{t}^c)$, *with*
*the property that*

$$\mu = \mu_z + \sum_{\beta \in (\Delta^+)'(\mathfrak{p}, \mathfrak{t}^c)} c_\beta \beta \qquad (0 \leq c_\beta \leq 1).$$

g) *There is an expression*

$$\mu = \mu_z + \sum_{\beta \in \Delta(\mathfrak{p}, \mathfrak{t}^c)} b_\beta \beta \qquad (0 \leq b_\beta \leq 1).$$

*Proof.* The equivalence of (a) and (b) is clear from Lemma 6.5 and the
definitions. For the equivalence of (b) and (c), recall from (3.1c) that $\lambda_u(\mu) =$
$T_\rho((\lambda_a)(\mu))$. Now the equivalence is clear from Proposition 1.7.

For the rest, the conditions in (d) through (g) are clearly satisfied if and
only if the analogous conditions on $\mu_s$ are satisfied. Lemma 6.5 therefore
allows us to assume that $G = G_s$ is connected and semisimple. Under this
assumption, we will prove that (a) implies (d); that (d) implies (e); that (e)
implies (f); that (f) implies (g); and finally that (g) implies (a). The most
difficult of these is (e) implies (f).

So suppose that $G$ is semisimple and that $\mu$ is unitarily small; that is, that
$\lambda_u(\mu) = 0$. According to (0.5) and Proposition 1.7, this means that $\mu + 2\rho_c - 2\rho$

belongs to $-C^o$. Since the fundamental weights $\xi_i$ belong to the positive Weyl chamber $C$, it follows that

$$\langle \xi_i, \mu + 2\rho_c - 2\rho \rangle \leq 0.$$

Because $2\rho = 2\rho_c + 2\rho_n$, this is equivalent to

$$\langle \xi_i, \mu \rangle \leq \langle \xi_i, 2\rho_n \rangle.$$

This is (d).

Next, suppose that (d) holds. Fix $(\Delta^+)'$ as in (e); say

$$(6.2a) \qquad\qquad (\Delta^+)' = w\Delta^+,$$

for some $w$ in the Weyl group. We may assume that the fundamental roots are labelled in such a way that $\xi_i' = w\xi_i$. Then

$$(6.2b) \qquad \langle \xi_i', \mu \rangle = \langle \xi_i', \mu + 2\rho_c \rangle - \langle \xi_i', 2\rho_c \rangle = \langle w\xi_i', \mu + 2\rho_c \rangle - \langle \xi_i', 2\rho_c \rangle.$$

By Lemma 1.8,

$$(6.2c) \qquad\qquad \langle w\xi_i, \mu + 2\rho_c \rangle \leq \langle \xi_i, \mu + 2\rho_c \rangle.$$

By our hypothesis (d), it follows that

$$(6.2d) \qquad \langle w\xi_i, \mu + 2\rho_c \rangle \leq \langle \xi_i, 2\rho \rangle = \langle w\xi_i, 2w\rho \rangle = \langle \xi_i', 2\rho' \rangle.$$

Inserting the inequality (6.2d) in (6.2b) yields

$$(6.2e) \qquad\qquad \langle \xi_i', \mu \rangle \leq \langle \xi_i', 2\rho' \rangle - \langle \xi_i', 2\rho_c \rangle = \langle \xi_i', 2\rho_n' \rangle.$$

This is (e).

Because the proof that (e) implies (f) is the most difficult part of the argument, we postpone it to the end. That (f) implies (g) is trivial. Assume therefore that (g) holds (still assuming $G$ is semisimple); we will deduce (a). Choose a positive system $\Delta^+$ making $\mu + 2\rho_c$ dominant. The formula for $\mu$ may then be written as

$$\mu = \sum_{\beta \in \Delta^+(\mathfrak{p}, \mathfrak{t}^c)} (b_\beta - b_{-\beta})\beta.$$

Here the coefficients $c_\beta = b_\beta - b_{-\beta}$ lie between $-1$ and $1$. We may therefore write

$$\mu = 2\rho_n + \sum_{\beta \in \Delta^+(\mathfrak{p}, \mathfrak{t}^c)} (c_\beta - 1)\beta;$$

the coefficients $a_\beta = c_\beta - 1$ in this formula lie between $0$ and $-2$. It follows that

$$\mu + 2\rho_c - 2\rho = \sum_{\beta \in \Delta^+(\mathfrak{p}, \mathfrak{t})} a_\beta \beta \qquad (-2 \leq a_\beta \leq 0).$$

Obviously the right side belongs to $-C^o$, the cone spanned by negative roots. By Proposition 1.7, $P(\mu + 2\rho_c - 2\rho) = 0$; that is, $\lambda_u(\mu) = 0$. This is (a).

Finally, we turn to the proof that (e) implies (f). We proceed by induction on the dimension of $G$. Recall that we are assuming that $G = G_s$ is semisimple. Lemma 6.6 actually allows us to assume that $G$ is simple. (More precisely, if $\mathfrak{g}_0$ has a nontrivial decomposition into simple factors, then the summands have lower dimension, so we already know the implication for them. Lemma 6.6 allows us to deduce the implication for $G$.)

LEMMA 6.8.    *Suppose $G$ is compact and simple. Then the only weight $\mu$ satisfying either condition* (e) *or condition* (f) *of Theorem 6.7 is $\mu = 0$.*

*Proof.* The condition in (e) is

$$\langle \xi, \mu \rangle \leq 0$$

for any fundamental weight $\xi$ (with respect to the fixed positive system making $\mu$ dominant). Because fundamental weights have nonnegative inner product, this implies $\mu = 0$. The argument for (f) is even simpler.    □

We may therefore assume that $G$ is noncompact and simple.

LEMMA 6.9.    *Suppose $G$ is noncompact and simple. Then the noncompact roots span $(\mathfrak{t}^c)^*$. In particular, in the setting of Theorem 6.7(e),*

$$\langle \xi_i', 2\rho_n' \rangle > 0.$$

*Proof.* Define

(6.3a)         $\Delta^0 = \{\alpha \in \Delta(\mathfrak{g}, \mathfrak{t}^c) \mid \langle \alpha, \beta \rangle = 0 \ (\beta \in \Delta(\mathfrak{p}, \mathfrak{t}))\},$

(6.3b)         $\Delta^1 = $ roots in the integer span of $\Delta(\mathfrak{p}, \mathfrak{t}^c).$

Obviously $\Delta^0$ and $\Delta^1$ are mutually orthogonal root subsystems in $\Delta$. We claim that

(6.3c)         $\Delta = \Delta^0 \cup \Delta^1.$

Assume this for a moment. The orthogonality of the decomposition allows us to conclude that root vectors for $\Delta^0$ and $\Delta^1$ must commute with each other, and therefore that we get a corresponding decomposition

(6.3d)         $\mathfrak{g}_0 = \mathfrak{g}_0^0 + \mathfrak{g}_0^1.$

Because $\mathfrak{g}_0$ is simple, this decomposition must be trivial. Since $G$ is assumed to be noncompact, the second factor is nontrivial; so the first must be trivial. Therefore $\Delta = \Delta^1$, which is the first claim in the lemma.

To prove (6.3c), suppose $\gamma$ is any root belonging neither to $\Delta^0$ nor to $\Delta^1$. Since $\Delta^1$ includes the noncompact roots, $\gamma$ must be compact. Since it is not in $\Delta^0$, there is a noncompact root $\beta$ with which it has nonzero inner product. Possibly replacing $\beta$ by $-\beta$, we get

$$\langle \gamma, \beta \rangle < 0.$$

It follows (by the representation theory of the SU(2) generated by $\gamma$) that $\beta + \gamma$ is a noncompact root. Therefore $\gamma = (\beta + \gamma) - \beta$ belongs to $\Delta^1$, a contradiction.

To prove the last inequality, notice that $\xi_i'$ has nonnegative inner product with each positive root. Therefore

(6.3e)

$$\Delta^+(\mathfrak{p}, \mathfrak{t}^c) = \{\beta \in \Delta(\mathfrak{p}, \mathfrak{t}) \mid \langle \xi_i', \beta \rangle > 0\} \cup \{\text{certain roots orthogonal to } \xi_i'\}.$$

So

(6.3f)

$$\langle \xi_i', 2\rho_n \rangle = \sum_{\substack{\beta \in \Delta(\mathfrak{p}, \mathfrak{t}^c) \\ \langle \xi_i', \beta \rangle > 0}} \langle \xi_i', \beta \rangle.$$

This is a sum of strictly positive terms; and it is nonempty by the first part of the lemma.  $\square$

We continue now with the proof that (e) implies (f) in Theorem 6.7. Recall that we are proceeding by induction on the dimension of $G$, and that we have reduced to the case when $G$ is simple and noncompact. We propose to apply the inductive hypothesis to a maximal Levi subgroup of $G$. In order to construct it, we choose a positive root system $(\Delta^+)'(\mathfrak{g}, \mathfrak{t}^c) \supset \Delta^+(\mathfrak{k}, \mathfrak{t}^c)$ and a fundamental weight $\xi_i'$ so as to maximize the quotient

(6.4a)

$$c_i = \langle \xi_i', \mu \rangle / \langle \xi_i', 2\rho_n' \rangle.$$

Notice that the denominator is strictly positive by Lemma 6.9. We can make the numerator nonnegative by an appropriate choice of $(\Delta^+)'$ (for example by making $\mu$ dominant); so our maximum $c_i$ must be nonnegative. By hypothesis (e) of Theorem 6.7, $c_i$ is bounded above by 1. Consequently

(6.4b)

$$0 \le c_i \le 1.$$

The weight $\xi_i'$ belongs to $i(\mathfrak{t}_0^c)^*$, so by (0.6) it defines a proper Levi subgroup

(6.4c)

$$L_i = G(\xi_i') \supset T^c.$$

The root system of $L_i$ is

(6.4d)

$$\Delta_i = \{\gamma \in \Delta \mid \langle \xi_i', \gamma \rangle = 0\} = \text{span of } \{\alpha_j' \mid j \ne i\}.$$

Our choice of positive roots for $K$ gives a positive system

(6.4e)

$$\Delta_i^+(\mathfrak{l}_i \cap \mathfrak{k}, \mathfrak{t}^c) = \Delta_i \cap \Delta^+(\mathfrak{k}, \mathfrak{t}^c).$$

The construction of (0.6) also provides a unipotent subalgebra

$$(6.4f) \qquad \mathfrak{u}_i = \mathfrak{u}(\xi_i'), \qquad \Delta(\mathfrak{u}_i, \mathfrak{t}) = \{\gamma \in \Delta \mid \langle \xi_i', \gamma \rangle > 0\}.$$

Recall from (0.4d) that $2\rho(\mathfrak{u}_i \cap \mathfrak{p})$ denotes the sum of the roots of $\mathfrak{t}^c$ in $\mathfrak{u}_i \cap \mathfrak{p}$; that is, the sum of the noncompact roots having positive inner product with $\xi_i'$. Define

$$(6.4g) \qquad \mu_i = \mu - 2c_i\rho(\mathfrak{u}_i \cap \mathfrak{p}).$$

LEMMA 6.10. *Use the notation of* (6.4).

a) *The weight $\mu_i$ belongs to the $\mathbb{R}$-span of the roots of $L_i$.*
b) *The weight $\mu_i$ is dominant integral for $L_i \cap K$. More precisely, if $\alpha$ is any root in $L_i \cap K$, then*

$$\langle \alpha, \mu_i \rangle = \langle \alpha, \mu \rangle.$$

c) *The weight $\mu_i$ satisfies condition* (e) *of Theorem* 6.7 *with respect to the group $L_i$.*

Let us assume this lemma for a moment, and complete our proof that (e) implies (f). Because $L_i$ has strictly lower dimension than $G$, we may assume by inductive hypothesis that (e) implies (f) for $L_i$. Parts (a) and (b) of Lemma 6.10 guarantee that Theorem 6.7 makes sense for $L_i$ and $\mu_i$; and part (c) says that Theorem 6.7(e) is satisfied. By induction, Theorem 6.7(f) must be satisfied as well. That is, we can find a positive system $(\Delta_i^+)''$ for $L_i$, containing $\Delta_i^+(\mathfrak{l}_i \cap \mathfrak{k}, \mathfrak{t}^c)$, so that

$$(6.5a) \qquad \mu_i = \sum_{\beta \in (\Delta_i^+)''(\mathfrak{l}_i \cap \mathfrak{p}, \mathfrak{t}^c)} c_\beta \beta \qquad (0 \le c_\beta \le 1).$$

Now define a positive root system for $G$ by

$$(6.5b) \qquad (\Delta^+)'' = (\Delta_i^+)'' \cup \Delta(\mathfrak{u}_i, \mathfrak{t}^c).$$

Then (6.4g) and (6.5a) show that

$$(6.5c) \qquad \mu = 2c_i\rho(\mathfrak{u}_i \cap \mathfrak{p}) + \mu_i$$
$$= \sum_{\beta \in \Delta(\mathfrak{u}_i \cap \mathfrak{p})} c_i\beta + \sum_{\beta \in (\Delta_i^+)''(\mathfrak{l}_i \cap \mathfrak{p}, \mathfrak{t}^c)} c_\beta \beta \qquad (0 \le c_\beta \le 1).$$

Because of (6.4b), this expression has the form required in Theorem 6.7(f).

It remains to prove Lemma 6.10. For (a), (1.2) shows that the $\mathbb{R}$-span of the roots of $L_i$ is precisely the orthogonal complement of $\xi_i'$. We calculate

$$\langle \xi_i', \mu_i \rangle = \langle \xi_i', \mu \rangle - c_i \langle \xi_i', 2\rho(\mathfrak{u}_i \cap \mathfrak{p}) \rangle.$$

Using (6.3f), we see that the second inner product on the right is just $\langle \xi_i', 2\rho_n' \rangle$. Therefore

$$\langle \xi_i', \mu_i \rangle = \langle \xi_i', \mu \rangle - c_i \langle \xi_i', 2\rho_n' \rangle.$$

Inserting the definition (6.4a) of $c_i$, we see immediately that the right side is zero, as desired.

For (b), it is enough to prove the formula. In light of (6.4g), the formula amounts to

$$\langle \alpha, 2\rho(\mathfrak{u}_i \cap \mathfrak{p}) \rangle = 0 \qquad (\alpha \in \Delta(\mathfrak{l}_i \cap \mathfrak{k}, \mathfrak{t}^c)).$$

Because the adjoint action of the reductive group $L_i \cap K$ preserves $\mathfrak{u}_i \cap \mathfrak{p}$, the sum of the weights of the latter must be orthogonal to the roots of $L_i \cap K$.

We turn now to (c). Suppose therefore that $(\Delta_i^+)''$ is a system of positive roots for $L_i$ containing $\Delta_i^+(\mathfrak{l}_i \cap \mathfrak{k}, \mathfrak{t}^c)$. Just as in (6.5) we define a positive root system for $G$ by

$$(6.6a) \qquad\qquad (\Delta^+)'' = (\Delta_i^+)'' \cup \Delta(\mathfrak{u}_i, \mathfrak{t}^c).$$

Write

$$(6.6b) \qquad\qquad \Pi'' = \{\alpha_1'', \ldots, \alpha_l''\}$$

for the set of simple roots, and $\{\xi_j''\}$ for the dual basis of fundamental weights. By construction, the two positive root systems $(\Delta^+)'$ and $(\Delta^+)''$ differ by an element of the Weyl group of $\Delta_i$. It follows that

$$(6.6c) \qquad\qquad \xi_i' = \xi_i''.$$

Furthermore the constants $c_i$ defined using $(\Delta^+)'$ and $(\Delta^+)''$ coincide, as do the weights $2\rho(\mathfrak{u}_i \cap \mathfrak{p})$. Consequently there is no loss of generality in assuming henceforth that $(\Delta^+)' = (\Delta^+)''$.

We are trying to establish an inequality relating $\mu_i$ and the fundamental weights for the positive system $(\Delta_i^+)'$. Let us write $(\xi_j^i)'$ for these fundamental weights. According to (1.2), they are characterized by belonging to the span of $\Delta_i$, and by

$$(6.7a) \qquad\qquad \langle (\xi_j^i)', \alpha_k' \rangle = \delta_{jk} \qquad (k, j \neq i).$$

From these requirements we see that $(\xi_j^i)'$ must be the projection of $\xi_j'$ orthogonal to $\xi_i'$. That is,

$$(6.7b) \qquad\qquad (\xi_j^i)' = \xi_j' - (\langle \xi_i', \xi_j' \rangle / \langle \xi_i', \xi_i' \rangle) \xi_i'.$$

The sum of the noncompact positive roots for $L_i$ is

$$(6.7c) \qquad\qquad 2(\rho_n^i)' = 2\rho_n' - 2\rho(\mathfrak{u}_i \cap \mathfrak{p}).$$

The inequality of Theorem 6.7(e) that we are trying to prove is

$$(6.7d) \qquad\qquad \langle (\xi_j^i)', \mu_i \rangle \leq \langle (\xi_j^i)', 2(\rho_n^i)' \rangle.$$

We begin by examining the right side. Because $2(\rho_n^i)'$ is a sum of roots in $L_i$, it is orthogonal to $\xi_i'$. So (6.7b) allows us to replace $(\xi_j^i)'$ by $\xi_j'$ on the right side. Inserting (6.7c) gives

$$(6.7e) \qquad \langle (\xi_j^i)', 2(\rho_n^i)' \rangle = \langle \xi_j', 2\rho_n' - 2\rho(\mathfrak{u}_i \cap \mathfrak{p}) \rangle.$$

Next we examine the left side of (6.7d). By construction $\mu_i$ is orthogonal to $\xi_i'$, so (6.7b) allows us to replace $(\xi_j^i)'$ by $\xi_j'$ here as well. Just as in (6.4a), let us define

$$(6.7f) \qquad c_j = \langle \xi_j', \mu \rangle / \langle \xi_j', 2\rho_n' \rangle.$$

Using this definition and (6.4g), we get

$$(6.7g) \qquad \begin{aligned} \langle (\xi_j^i)', \mu_i \rangle &= \langle \xi_j', \mu \rangle - c_i \langle \xi_j', 2\rho(\mathfrak{u}_i \cap \mathfrak{p}) \rangle \\ &= c_j \langle \xi_j', 2\rho_n' \rangle - c_i \langle \xi_j', 2\rho(\mathfrak{u}_i \cap \mathfrak{p}) \rangle. \end{aligned}$$

Recall now that we chose $i$ to maximize $c_i$. In particular we must have $c_j \le c_i$. Since the coefficient of $c_j$ on the right in (6.7g) is positive (Lemma 6.9), we deduce

$$(6.7h) \quad \langle (\xi_j^i)', \mu_i \rangle \le c_i \left( \langle \xi_j', 2\rho_n' \rangle - \langle \xi_j', 2\rho(\mathfrak{u}_i \cap \mathfrak{p}) \rangle \right) = c_i \langle \xi_j', 2\rho_n' - 2\rho(\mathfrak{u}_i \cap \mathfrak{p}) \rangle.$$

Now (6.7e) and the fact that $c_i \le 1$ (cf. (6.4b)) show that (6.7h) implies the inequality (6.7d). This completes the proof of Lemma 6.10, and therefore of Theorem 6.7. $\qquad\qquad\square$

We want next to explore some formal consequences of the characterizations of unitarily small $K$-types in Theorem 6.7. We will apply the following lemma to the set of roots of $T^c$ in $\mathfrak{p}$.

LEMMA 6.11.    Use the notation of (0.3). Suppose $S \subset i(\mathfrak{t}_0^c)^*$ is an arbitrary finite subset invariant under the Weyl group of $K$:

$$w \cdot S \subset S \qquad (w \in W(K,T)).$$

Define

$$R(S) = \left\{ \sum_{s \in S} b_s s \qquad (0 \le b_s \le 1) \right\} \subset i(\mathfrak{t}_0^c)^*.$$

Then $R(S)$ is convex and invariant under $W(K)$. The following conditions on an irreducible representation $\delta \in \widehat{K}$ are equivalent.

a) Every weight of $\delta$ belongs to $R(S)$.
b) Every extremal weight of $\delta$ belongs to $R(S)$.
c) Every highest weight of $\delta$ belongs to $R(S)$.
d) Some highest weight of $\delta$ belongs to $R(S)$.
e) Some extremal weight of $\delta$ belongs to $R(S)$.

We will say that a representation of $K$ (even reducible) is *of type $S$* if these conditions are satisfied.

*Proof.* That $R(S)$ is convex and $W(K)$ invariant is obvious, as are the implications (a) $\Rightarrow$ (b) $\Rightarrow$ (c) $\Rightarrow$ (d) $\Rightarrow$ (e). We will complete the proof by showing that (e) $\Rightarrow$ (b) $\Rightarrow$ (a). So assume (e) holds. The set of extremal weights of $\delta$ is a single orbit of $W(K)$, so (b) follows by the $W(K)$ invariance of $R(S)$. Now assume (b). Every weight must lie in the convex hull of the extremal weights. (For connected $K$ this is a consequence of Proposition 1.7, and the general case follows at once.) So the convexity of $R(S)$ shows that (b) implies (a). $\qquad\square$

COROLLARY 6.12.   *Suppose $G$ is a real reductive group in Harish-Chandra's class, and $K$ is a maximal compact subgroup; use the notation of (6.1). The following conditions on an irreducible representation $\delta \in \widehat{K}$ are equivalent.*

a) *The $K$-type $\delta$ is unitarily small.*
b) *The restriction $\delta_s$ of $\delta$ to $G_s \cap K$ is of type $\Delta(\mathfrak{p}, \mathfrak{t}^c)$ (Lemma 6.11).*
c) *Write $\mu_z \in i\mathfrak{z}^*_{c,0}$ for the differential of the character of $Z_c$ on $\delta$. Then $\delta$ is of type $\{\mu_z\} \cup \Delta(\mathfrak{p}, \mathfrak{t}^c)$.*

*Proof.* In light of Lemmas 6.11 and 6.5, this is a reformulation of the equivalence of (a) and (g) in Theorem 6.7. $\qquad\square$

Before we embark on a more serious investigation of unitarily small $K$-types, we want to describe the unitarily small parameters of Definition 6.1. The result (Proposition 6.13) is fairly simple, but the proof is a little involved. We need to extend somewhat the notation introduced in (6.1). Recall that $\mathfrak{g}_{s,0} = [\mathfrak{g}_0, \mathfrak{g}_0]$ is the derived algebra of $\mathfrak{g}_0$. It is a direct sum of simple ideals (Lemma 6.6). Write

(6.8a)        $\mathfrak{g}_{c,0}$ = sum of the compact simple factors of $\mathfrak{g}_{s,0}$.

The corresponding subgroup $G_c$ of $G_s$ is compact, connected, and semisimple. Therefore it has finite center; so in particular the intersection

(6.8b)                    $F_c = Z(G_c) \cap Z \subset Z_c$

of its center with the identity component of the center of $G$ must be finite.

PROPOSITION 6.13.   *Assume the setting of Definition* 6.1; *use also the notation of* (6.1) *and* (6.8). *The set $\Lambda^z_u(G)$ of unitarily small parameters for $G$ consists precisely of the differentials of characters of $Z_c/F_c$.*

*Proof.* Suppose $\delta \in B^z_u(G)$ is a unitarily small representation of $K$. It is clear from conditions (f) or (g) of Theorem 6.7 that $\delta$ must be trivial on $G_c$.

(This was essentially proved in Lemma 6.8.) The restriction of $\delta$ to $Z_c$ therefore defines a character $\mu_z$ of $Z_c/F_c$, and Lemma 6.5 shows that the corresponding unitarily small parameter is the differential of $\mu_z$.

Conversely, suppose $\mu_z$ is a character of $Z_c/F_c$. We must show how to extend $\mu_z$ to a unitarily small representation of $K$. The problem is easily reduced to the identity component of $G$, so we assume henceforth that $G$ is connected. Define

(6.9a)          $G_n = $ identity component of the centralizer in $G$ of $G_c$.

Clearly this is a closed connected reductive subgroup of $G$. On the level of Lie algebras it is easy to prove a direct sum decomposition

(6.9b)          $$\mathfrak{g}_0 = \mathfrak{g}_{c,0} + \mathfrak{g}_{n,0}.$$

From this it follows that the multiplication map defines a surjective covering homomorphism

(6.9c)          $$G_c \times G_n \to G.$$

Define

(6.9d)          $$F = G_c \cap G_n = Z(G_c) \cap Z(G_n).$$

This is a finite central subgroup of both factors. Its intersection with the identity component of the center of $G_n$ is just the group $F_c$ defined in (6.8b). The kernel of the multiplication map is

(6.9e)          $$F_\Delta = \{(f, f^{-1}) \mid f \in F\} \subset G_c \times G_n.$$

Now (6.9) allows us to describe the representation theory of $G$ completely in terms of the representation theory of $G_c$ and $G_n$ separately. Here is an example.

LEMMA 6.14.          *Suppose $G$ is connected; use the notation of* (6.9). *The irreducible representations $\delta$ of $K$ are in one-to-one correspondence with pairs*

$$(\delta_c, \delta_n) \in \widehat{G}_c \times (G_n \cap K)\widehat{\phantom{o}},$$

*subject to the requirement that*

$$\delta_c|_F = \delta_n|_F.$$

*Proof.* Because of (6.9), the multiplication map defines an isomorphism

$$K \simeq (G_c \times (G_n \cap K))/F_\Delta.$$

The representations of the direct product $G_c \times (G_n \cap K)$ are tensor products, and so correspond to pairs $(\delta_c, \delta_n)$. The requirement in the lemma is exactly what is needed to make $\delta_c \otimes \delta_n$ trivial on $F_\Delta$.          $\square$

We now set out to construct a unitarily small representation of $K$ from our character $\mu_z$ of $Z_c/F_c$. As a first step, consider the subgroup $FZ_c$ of $G_n$. Because $F \cap Z_c = F_c$,

$$(6.10a) \qquad\qquad FZ_c/F \simeq Z_c/F_c.$$

So we begin by regarding $\mu_z$ as a character of $FZ_c/F$. Because this is a compact central subgroup of $G_n$, we can find an irreducible representation $\delta_n^1$ of $G_n \cap K$ such that

$$(6.10b) \qquad\qquad \delta_n^1|_{FZ_c} = \mu_z.$$

If $\delta_n^1$ were unitarily small, then $\mathbb{C} \otimes \delta_n^1$ would be a unitarily small representation of $K$ restricting to $\mu_z$, and we would be done. But $\delta_n^1$ may not be unitarily small, so we need to work a little harder.

Write $T_n$ for our fixed maximal torus in $G_n \cap K$ (the intersection of $T^c$ with $G_n$), and $\mu_n^1 \in \widehat{T}_n$ for any weight of $\delta_n^1$. By hypothesis

$$(6.10c) \qquad\qquad \mu_n^1|_{FZ_c} = \mu_z.$$

On the other hand, $G_n$ has no compact simple factors. By Lemma 6.9, this implies that the noncompact roots span the semisimple part of $i\mathfrak{t}_{n,0}^*$. We may therefore write

$$(6.10d) \qquad\qquad \mu_n^1 = \mu_z + \sum_{\beta \in \Delta(\mathfrak{p}, \mathfrak{t}_n)} b_\beta^1 \beta \qquad (b_\beta^1 \in \mathbb{R}).$$

For each noncompact root $\beta$, choose an integer $n_\beta$ so that

$$(6.10e) \qquad\qquad b_\beta = b_\beta^1 + n_\beta \in [0, 1].$$

Now define $\mu_n \in \widehat{T}^c$ to be $\mu_n^1$ twisted by the sum of the various $n_\beta \beta$:

$$(6.10f) \qquad\qquad \mu_n = \mu_n^1 + \sum_{\beta \in \Delta(\mathfrak{p}, \mathfrak{t}_n)} n_\beta \beta.$$

Because the roots are characters of $T_n$, this is well-defined. Because the center of $G_n$ acts trivially on the roots, (6.10c) implies

$$(6.10g) \qquad\qquad \mu_n|_{FZ_c} = \mu_z.$$

Finally, (6.10d–f) give

$$(6.10h) \qquad\qquad \mu_n = \mu_z + \sum_{\beta \in \Delta(\mathfrak{p}, \mathfrak{t}_n)} b_\beta \beta \qquad (0 \le b_\beta \le 1).$$

Now define $\delta_n$ to be the unique irreducible representation of $G_n \cap K$ of extremal weight $\mu_n$. By Lemma 6.11, the highest weight of $\delta_n$ has an expression of the same form as in (6.10h). By Theorem 6.7, $\delta_n$ is unitarily small. By Lemma 6.14, $\delta = \mathbb{C} \otimes \delta_n$ is a unitarily small representation of $K$, and the restriction of $\delta$ to $Z_c$ is $\mu_z$ by (6.10g).

## 7. Unitarily small $K$-types and the Dirac operator

PROPOSITION 7.1.    *In the setting of Definition* 6.1, *any irreducible representation of $K$ appearing in the exterior algebra $\bigwedge \mathfrak{p}$ (notation* (0.2)) *is unitarily small. More precisely, the representation $\bigwedge \mathfrak{p}$ is of type $\Delta(\mathfrak{p}, \mathfrak{t}^c)$ (cf. Lemma* 6.11).

*Proof.* The weights of the exterior algebra of $\mathfrak{p}$ are those of the form

$$\sum_{\beta \in \Delta(\mathfrak{p}, \mathfrak{t}^c)} b_\beta \beta, \qquad (b_\beta = 0 \text{ or } 1).$$

Therefore $\bigwedge \mathfrak{p}$ is of type $\Delta(\mathfrak{p}, \mathfrak{t})$. Now apply Corollary 6.12.    □

In the theory of unitary representations, the interaction between the $K$-types of $\bigwedge \mathfrak{p}$ and certain spinor representations plays an important role. In this section we will outline this theory, and show how to extend parts of it to all unitarily small representations.

Our development is taken from [4, Chapter II] and [2, II.6]. We will indicate particularly large departures from their conventions. We begin with a finite-dimensional real inner product space $(V_0, \langle \ , \ \rangle)$. Set

(7.1a)        $\mathrm{SO}(V_0) =$ special orthogonal group of $V_0$,        $V = V_0 \otimes_\mathbb{R} \mathbb{C}$.

The *Clifford algebra* $C(V_0)$ is the real associative algebra with unit generated by $V_0$ and subject to the relations

(7.1b)                $v^2 = -\langle v, v \rangle \cdot 1 \qquad (v \in V_0)$.

(For Chevalley this is the Clifford algebra of the form $-\langle \ , \ \rangle$.) Sometimes it is convenient to write the defining relations in the equivalent form

(7.1c)                $vw + wv = -2\langle v, w \rangle \qquad (v, w \in V_0)$.

Obviously $\mathrm{SO}(V_0)$ (and even the full orthogonal group $\mathrm{O}(V_0)$) acts on $C(V_0)$ by algebra automorphisms. The action of $-1 \in \mathrm{O}(V_0)$ defines a $\mathbb{Z}/2\mathbb{Z}$-grading

(7.1d)                $C(V_0) = C(V_0)_{\text{even}} + C(V_0)_{\text{odd}}$.

There is an algebra antiautomorphism *transpose* of $C(V_0)$, characterized by

(7.1e)        ${}^t v = -v, \qquad {}^t(ab) = {}^t b \, {}^t a \qquad (v \in V_0, \ a, b \in C(V_0))$.

(Chevalley emphasizes the antiautomorphism $\alpha$ that acts by $+1$ on $V_0$. It differs from transpose by the action of $-1 \in \mathrm{O}(V_0)$.) The Clifford algebra is filtered by assigning the generators degree 1. Write $C_p(V_0)$ for the $p^{\text{th}}$ level of the filtration. Then there is a natural isomorphism

$$(7.1f) \qquad C_p(V_0)/C_{p-1}(V_0) \simeq \bigwedge^p V_0$$

respecting the algebra structure and the $\mathbb{Z}/2\mathbb{Z}$-grading. It follows easily that there is a natural isomorphism

$$(7.1g) \qquad C_1(V_0)_{\text{odd}} \simeq V_0;$$

the map from right to left just sends the space $V_0$ of generators into $C(V_0)$.

Write $C(V_0)^\times$ for the group of invertible elements of $C(V_0)$. Because $C(V_0)$ is a finite-dimensional associative algebra with a unit, $C(V_0)^\times$ is an open subset. It is a Lie group; the Lie algebra is $C(V_0)$ with the commutator bracket. This group acts on $C(V_0)$ by conjugation:

(7.2a)
$$\text{Ad}: C(V_0)^\times \to \text{Aut}(C(V_0)), \quad \text{Ad}(u)(a) = uau^{-1} \quad (u \in C(V_0)^\times, a \in C(V_0)).$$

The differential of this representation is the action of $C(V_0)$ on itself by bracket. The subgroup $C(V_0)^\times_{\text{even}}$ acts in $\text{Ad}$ to preserve the $\mathbb{Z}/2\mathbb{Z}$ grading. We now define the *spin group of $V_0$* as

$$(7.2b) \qquad \text{Spin}(V_0) = \{g \in C(V_0)^\times_{\text{even}} \mid \text{Ad}(g)(V_0) \subset V_0, \ {}^t g = g^{-1}\}.$$

This is the group that Chevalley calls the *reduced Clifford group* $\Gamma_0^+$. The first condition here may also be expressed as requiring $\text{Ad}(g)$ to preserve the filtration of 7.1(f). By construction, $\text{Spin}(V_0)$ is an algebraic group, and the adjoint action provides an algebraic homomorphism

$$\tau: \text{Spin}(V_0) \to \text{GL}(V_0).$$

Because $\text{Ad}$ preserves the algebra structure, the relation (7.1b) shows that the image of $\tau$ is contained in $\text{O}(V_0)$:

$$(7.2c) \qquad \tau: \text{Spin}(V_0) \to \text{O}(V_0).$$

LEMMA 7.2. *The image of the map $\tau$ of (7.2c) is the special orthogonal group* $\text{SO}(V_0)$. *The kernel consists of the two scalar elements $\pm 1$ in the center of $C(V_0)$. If $\dim V_0 \geq 2$, this kernel is contained in the identity component of* $\text{Spin}(V_0)$; *so in that case $\text{Spin}(V_0)$ is a connected double cover of $\text{SO}(V_0)$.*

PROPOSITION 7.3. *The Clifford algebra $C(V_0)$ is a semisimple algebra of dimension $2^{\dim V_0}$. Write $Z = Z(C(V_0))$ for the center.*

*Suppose that $\dim V_0 = 2r$ is even. Then $Z$ consists only of the scalars $C_0$. Consequently $C(V_0)$ and its complexification are both central simple algebras.*

In particular, $C(V_0)$ has a unique complex simple module $(\gamma, S)$ (up to equivalence), of dimension $2^r$. Regard $\gamma$ as a representation of the group $C(V_0)^\times$. Write $\sigma$ for the restriction of $\gamma$ to the subgroup $\mathrm{Spin}(V_0)$. If $r \geq 1$, then $\sigma$ is the direct sum of two inequivalent irreducible representations $(\sigma_\pm, S_\pm)$, each of dimension $2^{r-1}$.

Suppose that $\dim V_0 = 2r + 1$ is odd. Then $Z$ is two-dimensional; it is spanned by the scalars and an odd element $z_-$, satisfying

(a) $$(z_-)^2 = \pm 1 \in C_0.$$

The subalgebra $C(V_0)_{\mathrm{even}}$ has dimension $2^{2r}$ and one-dimensional center; it is central simple, and generated as an algebra by $\mathrm{Spin}(V_0)$. The entire Clifford algebra is a tensor product

(b) $$C(V_0) \simeq Z \otimes_{\mathbb{R}} C(V_0)_{\mathrm{even}}.$$

There is a unique complex simple module $(\gamma_{\mathrm{even}}, S)$ for $C(V_0)_{\mathrm{even}}$ (up to equivalence), of dimension $2^r$. The restriction of $\gamma$ to $\mathrm{Spin}(V_0)$ is an irreducible representation $(\sigma, S)$. There are exactly two extensions $\gamma$ and $\gamma'$ of $\gamma_{\mathrm{even}}$ to a module for $C(V_0)$ (on the same space $S$); they are related by

(c) $$\gamma(v) = -\gamma'(v) \qquad (v \in V_0).$$

Any complex simple module for $C(V_0)$ is isomorphic to exactly one of $\gamma$ and $\gamma'$.

In either case, suppose $(\gamma, S)$ is a simple module for $C(V_0)$. Then there is a positive definite invariant Hermitian form $\langle \ , \ \rangle$ on $S$. This means that $\gamma(a)^* = \gamma({}^t a)$, or explicitly that

(d) $$\langle \gamma(a)s, s' \rangle = \langle s, \gamma({}^t a)s' \rangle.$$

This form is unique up to a positive scalar multiple. In particular, the representation $(\sigma, S)$ of $\mathrm{Spin}(V_0)$ is unitary.

We should make two remarks about the notation introduced in this proposition. There is no preferred choice between the two irreducible constituents of the spin representation in the even-dimensional case; to call one of them $\sigma_+$ is to make a choice. Similarly, there is no preferred extension of $\gamma_{\mathrm{even}}$ in the odd-dimensional case; calling one $\gamma$ is making a choice.

This proposition summarizes several results in Chapter 2 of [4]. (The formulation in II.6.2 of [2] is not quite correct in the odd-dimensional case, when it overlooks the nonuniqueness of the extension of $\gamma_{\mathrm{even}}$ to $C(V_0)$.)

*Definition* 7.4. Suppose we are in the setting 7.1–7.2. A *space of spinors for* $V_0$ is a complex simple module $(\gamma, S)$ for the Clifford algebra $C(V_0)$, endowed with a positive definite invariant Hermitian form $\langle \ , \ \rangle$, as described in

**Proposition 7.3.** We write $(\sigma, S)$ for the corresponding unitary representation of $\mathrm{Spin}(V_0)$; we will sometimes call $\sigma$ a *spin representation*. By 7.2(a) and the definition of $\sigma$,

$$\gamma(\mathrm{Ad}(g)(a)) = \sigma(g)\gamma(a)\sigma(g^{-1}) \qquad (g \in \mathrm{Spin}(V_0), a \in C(V_0)).$$

*Example* 7.5. Suppose $V_0$ is a one-dimensional space $\mathbb{R}v_0$, and that $v_0$ has length one. Then $\mathrm{SO}(V_0)$ is trivial. The Clifford algebra $C(V_0)$ is generated by $v_0$ with the relation $v_0^2 = -1$; that is, $C(V_0)$ is isomorphic to the algebra $\mathbb{C}$ of complex numbers. There are two possible isomorphisms, sending $v_0$ to $i$ or to $-i$. The even subalgebra is

$$C(V_0)_{\mathrm{even}} = \mathbb{R} \subset \mathbb{C},$$

and the odd subspace is $i\mathbb{R}$. The transpose automorphism is complex conjugation, so we see that

$$\mathrm{Spin}(V_0) = \pm 1 \subset \mathbb{R}^{\times}.$$

The unique complex simple module for $C(V_0)_{\mathrm{even}}$ is just $\mathbb{C}$; the two possible extensions to $C(V_0)$ have $v_0$ acting by $\pm i$. (The distinction between these is not just a matter of choosing $i \in \mathbb{C}$; we must also choose between $v_0$ and $-v_0$ as a basis element of $V_0$.)

*Example* 7.6. Suppose $V_0$ is two-dimensional, with an orthonormal basis $\{j, k\}$. Then $\mathrm{SO}(V_0)$ consists of the $2 \times 2$ matrices

(a) $$\mathrm{SO}(V_0) = \left\{ r(\theta) = \begin{pmatrix} \cos\theta & -\sin\theta \\ \sin\theta & \cos\theta \end{pmatrix} \right\}.$$

The Clifford algebra has basis $\{1, j, k, jk\}$, with relations

(b) $$j^2 = k^2 = -1, \quad jk = -kj.$$

This is evidently isomorphic to the algebra of quaternions

(c) $$\mathbb{H} = \{a + bi + cj + dk\}$$

by sending $jk$ to $i$. The transpose antiautomorphism is the usual conjugation of quaternions (which is $-1$ on $i$, $j$, and $k$). The even subalgebra is

(d) $$C(V_0)_{\mathrm{even}} = \mathbb{C} \subset \mathbb{H}.$$

It is not difficult to check from (d) and the definition (7.2b) that

(e) $$\mathrm{Spin}(V_0) = \{s(\phi) = \cos\phi + i\sin\phi \mid \phi \in \mathbb{R}\}.$$

The conjugation action on $V_0$ is

(f) $$\tau(s(\phi)) = r(2\phi).$$

To see that, one calculates for example

$$\begin{aligned}
\tau(s(\phi))j &= (\cos\phi + i\sin\phi)j(\cos\phi - i\sin\phi)\\
&= j\cos^2\phi - iji\sin^2\phi + (ij - ji)\sin\phi\cos\phi\\
&= j(\cos^2\phi - \sin^2\phi) + k(2\sin\phi\cos\phi)\\
&= j\cos 2\phi + k\sin 2\phi\\
&= r(2\phi)j.
\end{aligned}$$

Therefore $\text{Spin}(V_0)$ is the double cover of the circle $\text{SO}(V_0)$. The unique complex simple module for $C(V_0)$ may be realized on $S = \mathbb{C}^2$ by

$$\text{(g)}\qquad \gamma(j) = \begin{pmatrix} 0 & 1 \\ -1 & 0 \end{pmatrix}, \quad \gamma(k) = \begin{pmatrix} 0 & i \\ i & 0 \end{pmatrix}, \quad \gamma(i) = \begin{pmatrix} i & 0 \\ 0 & -i \end{pmatrix}.$$

The natural decomposition $\mathbb{C}^2 = \mathbb{C} + \mathbb{C}$ gives a $\text{Spin}(V_0)$-invariant decomposition $S = S_+ + S_-$. Combining (e) with the third formula in (g), we see that

$$\text{(h)}\qquad\qquad \sigma_{\pm}(r(\phi)) = \exp(\pm i\phi).$$

That is, the two half-spin representations $\sigma_{\pm}$ are the first two characters of the circle $\text{Spin}(V_0)$.

One reason for considering these two examples in such detail is that many aspects of general Clifford algebras and spin representations can be reduced to dimensions one and two. In order to understand how, we need a little more notation.

Suppose first that $W_0 \subset V_0$ is any subspace. We immediately get a natural map of Clifford algebras

$$\text{(7.3a)}\qquad\qquad C(W_0) \to C(V_0)$$

respecting $\mathbb{Z}/2\mathbb{Z}$ gradings, filtrations, and the transpose. By (7.1f) this map is one-to-one; we use it to identify $C(W_0)$ as a subalgebra of $C(V_0)$. It is clear from the definition in 7.2 that (7.3a) induces an inclusion of groups

$$\text{(7.3b)}\qquad\qquad \text{Spin}(W_0) \hookrightarrow \text{Spin}(V_0).$$

Suppose now that $U_0$ is a second subspace of $V_0$, and that

$$\text{(7.3c)}\qquad\qquad U_0 \text{ is orthogonal to } W_0.$$

Then (7.1c) says that

$$\text{(7.3d)}\qquad\qquad wu + uw = 0 \qquad (u \in U_0,\ w \in W_0).$$

From this we can get complete commutation or anticommutation rules for even and odd elements of $C(U_0)$ and $C(V_0)$; but we will be content with

$$\text{(7.3e)}\qquad\qquad ab - ba = 0 \qquad (a \in C(W_0)_{\text{even}}, b \in C(U_0)).$$

In particular,

(7.3f)                        Spin($W_0$) and Spin($U_0$) commute.

The isomorphism (7.1f) also shows that $C(W_0)$ and $C(U_0)$ meet only in their common scalars. Lemma 7.2 therefore shows that

(7.3g)                        $\mathrm{Spin}(W_0) \cap \mathrm{Spin}(U_0) = \{\pm 1\},$

$$\mathrm{Spin}(W_0) \times_{\{\pm 1\}} \mathrm{Spin}(U_0) \hookrightarrow \mathrm{Spin}(V_0).$$

Here the notation means that we divide the product group by the diagonal copy of $\{\pm 1\}$.

The formula (7.3f) suggests that an orthogonal decomposition of $V_0$ ought to give rise to a tensor product decomposition of the Clifford algebra. This is complicated by the strange shape of the commutation relation (7.3d). The following bit of extra structure will allow us to circumvent this.

*Definition* 7.7.   Suppose we are in the setting 7.1–7.2, and that dim $V_0$ is even. A *special element* $\varepsilon \in \mathrm{Spin}(V_0)$ is one satisfying

$$\tau(\varepsilon) = -1 \in \mathrm{SO}(V_0).$$

Such an element exists and is unique up to sign by Lemma 7.2. Because $\tau(\varepsilon^2) = 1$, we must have $\varepsilon^2 = \pm 1$. In Example 7.6, $\varepsilon = s(\pi/2) = i$ is a special element of square $-1$. If $V_0 = 0$, then $\varepsilon = \pm 1$ are special elements of square 1. It turns out that $\varepsilon^2 = (-1)^r$ when dim $V_0 = 2r$, but we will not need this.

PROPOSITION 7.8 ([4, II.2.5]).   *In the setting 7.1–7.2, suppose that there is an orthogonal decomposition*

(a)                        $V_0 = W_0 + U_0,$

*and that* dim $W_0$ *is even. Choose a special element* $\varepsilon \in \mathrm{Spin}(W_0) \subset C(W_0)_{\mathrm{even}}$.

*Suppose first that* $\varepsilon^2 = 1$. *Then there is an isomorphism of algebras*

(b)                        $j_{\mathbb{R}} \colon C(W_0) \otimes_{\mathbb{R}} C(U_0) \to C(V_0),$

*defined by*

(c)          $j_{\mathbb{R}}(a \otimes b) = \begin{cases} ab, & \text{if } a \in C(W_0) \text{ and } b \in C(U_0)_{\mathrm{even}}; \\ a\varepsilon b, & \text{if } a \in C(W_0) \text{ and } b \in C(U_0)_{\mathrm{odd}}. \end{cases}$

*This isomorphism complexifies to*

(d)                        $j_{\mathbb{C}} \colon C(W_0)_{\mathbb{C}} \otimes_{\mathbb{C}} C(U_0)_{\mathbb{C}} \to C(V_0)_{\mathbb{C}}.$

*Suppose next that* $\varepsilon^2 = -1$. *Then there is an isomorphism of complexified algebras as in* (d), *given by*

(e)          $j_{\mathbb{C}}(a \otimes b) = \begin{cases} ab, & \text{if } a \in C(W_0)_{\mathbb{C}} \text{ and } b \in C(U_0)_{\mathbb{C},\mathrm{even}}; \\ ia\varepsilon b, & \text{if } a \in C(W_0)_{\mathbb{C}} \text{ and } b \in C(U_0)_{\mathbb{C},\mathrm{odd}}. \end{cases}$

*Suppose* $(\gamma(W_0), S(W_0))$ *and* $(\gamma(U_0), S(U_0))$ *are complex simple modules for* $C(W_0)$ *and* $C(U_0)$ *(Proposition 7.3). Then*

$$(f) \qquad S = S(W_0) \otimes S(U_0)$$

*is a simple module for* $C(V_0)$ *under the isomorphism* $j_{\mathbb{C}}$ *of* (d). *In this way the restriction of the spin representation* $\sigma(V_0)$ *of* $\mathrm{Spin}(V_0)$ *to the subgroup* $\mathrm{Spin}(W_0)\mathrm{Spin}(U_0)$ *(Def. 7.4 and (7.3g)) is identified with the tensor product of the spin representations of the factors:*

$$(g) \qquad \sigma(V_0)|_{\mathrm{Spin}(W_0)\mathrm{Spin}(U_0)} \simeq \sigma(W_0) \otimes \sigma(U_0).$$

COROLLARY 7.9. *In the setting 7.1–7.2, suppose* $\dim V_0 = 2r + \delta$, *with* $\delta$ *equal to 0 or 1. Choose a collection of* $r$ *orthogonal two-dimensional subspaces* $W_1, \ldots, W_r$ *of* $V_0$. *Recall from Example 7.6 that each group* $\mathrm{Spin}(W_i)$ *is a double cover of the circle group* $\mathrm{SO}(W_i)$, *with two-element kernel* $\pm 1_i$:

$$1 \to \{\pm 1_i\} \to \mathrm{Spin}(W_i) \to \mathrm{SO}(W_i) \to 1.$$

*Write*

$$T(W_1, \ldots, W_r) = \mathrm{SO}(W_1) \times \cdots \times \mathrm{SO}(W_r),$$

$$\widetilde{T}(W_1, \ldots, W_r) = \mathrm{Spin}(W_1) \times_{\{\pm 1\}} \cdots \times_{\{\pm 1\}} \mathrm{Spin}(W_r).$$

*Here the second definition means the quotient of the direct product by the subgroup generated by all elements* $(-1_i)(-1_j)$ $(1 \le i, j \le r)$. *This is a double cover of the corresponding product of orthogonal groups:*

$$1 \to \{\pm 1\} \to \widetilde{T}(W_1, \ldots, W_r) \to T(W_1, \ldots, W_r) \to 1.$$

a) *The group* $T(W_1, \ldots, W_r)$ *is a maximal torus in* $\mathrm{SO}(V_0)$.

b) *The group* $\widetilde{T}(W_1, \ldots, W_r)$ *may be naturally identified as a subgroup of* $\mathrm{Spin}(V_0)$; *it is a maximal torus.*

c) *Write* $\pm\beta_i$ *for the two weights of* $\mathrm{SO}(W_i)$ *on* $(W_i)_{\mathbb{C}}$; *extend them to all of* $T(W_1, \ldots, W_r)$ *by making them trivial on the other factors* $\mathrm{SO}(W_j)$. *Then the weights of* $T(W_1, \ldots, W_r)$ *on* $(V_0)_{\mathbb{C}}$ *are* $\{\pm\beta_1, \ldots, \pm\beta_r\}$, *together with* 0 *if* $\delta = 1$.

d) *The weights of the spin representation* $\sigma(V_0)$ *of* $\mathrm{Spin}(V_0)$ *are the* $2^r$ *weights* $\{1/2(\pm\beta_1 \pm \cdots \pm \beta_r)\}$.

COROLLARY 7.10 ([2, Lemma II.6.5]). *In the setting 7.1–7.2, choose a spin representation* $(\sigma, S)$ *(Def. 7.4). Then the tensor product representation* $\sigma \otimes \sigma$ *of* $\mathrm{Spin}(V_0)$ *descends to* $\mathrm{SO}(V_0)$ *by the map of (7.2c). As representations of* $\mathrm{SO}(V_0)$,

$$\bigwedge(V_0)_{\mathbb{C}} \simeq \begin{cases} \sigma \otimes \sigma & \text{if } \dim V_0 \text{ is even;} \\ \sigma \otimes \sigma + \sigma \otimes \sigma & \text{if } \dim V_0 \text{ is odd.} \end{cases}$$

Finally we introduce the Dirac operator. In the setting (0.2), recall that our fixed bilinear form $\langle\ ,\ \rangle$ is positive definite on $\mathfrak{p}_0$. We will make the constructions based on 7.1–7.2 with $\mathfrak{p}_0$ playing the role of $V_0$. The adjoint action of $K$ on $\mathfrak{p}_0$ preserves the form, and so defines a homomorphism

(7.4a)                          $\mathrm{Ad}\colon K \to \mathrm{O}(\mathfrak{p}_0).$

Of course the identity component $K_0$ is mapped to $\mathrm{SO}(\mathfrak{p}_0)$. Recall from (7.2c) and Lemma 7.2 that we also have a spin double cover of $\mathrm{SO}(\mathfrak{p}_0)$. This covering pulls back by Ad to a double cover of $K_0$. More precisely, we define

(7.4b)          $\widetilde{K}_0 = \{(k,g) \in K_0 \times \mathrm{Spin}(\mathfrak{p}_0) \mid \mathrm{Ad}(k) = \tau(g)\}.$

Then projection on the first factor defines a short exact sequence

(7.4c)                    $1 \to \{\pm 1\} \to \widetilde{K}_0 \to K_0 \to 1.$

Similarly, projection on the second factor defines a homomorphism

(7.4d)                          $\widetilde{\mathrm{Ad}}\colon \widetilde{K}_0 \to \mathrm{Spin}(\mathfrak{p}_0).$

Fix now a space of spinors $(\gamma, S)$ for $\mathfrak{p}_0$ (Definition 7.4). By composition with $\widetilde{\mathrm{Ad}}$, the spin representation defines a representation

(7.4e)                          $(\sigma \circ \widetilde{\mathrm{Ad}}, S)$

of $\widetilde{K}_0$ on $S$. We call this *the spin representation of* $\widetilde{K}_0$, and often denote it simply $\sigma$.

Suppose now that $X$ is a $(\mathfrak{g}, K)$-module. Often we can understand very well the action of $K$ on $X$, in terms of the Cartan-Weyl theory of finite-dimensional representations. What remains, therefore, is to understand the action of $\mathfrak{p}_0$ on $X$. The Dirac operator provides a way to compare this action with the action of $\mathfrak{p}_0$ on $S$ by Clifford multiplication. Because the latter action is also well understood, the problem of understanding $X$ is to some extent reduced to understanding the Dirac operator.

*Definition* 7.11. Suppose $G$ is a reductive group as in (0.2). Use the notation of (7.4), so that in particular we have chosen a space $(\gamma, S)$ of spinors for $\mathfrak{p}_0$. Suppose that $(\pi, X)$ is a representation of $\mathfrak{g}$. The *Dirac operator for* $X$ is an endomorphism $D$ of $X \otimes S$, defined as follows. Let $\{Y_i\}$ be a basis of $\mathfrak{p}_0$, and $\{Y^j\}$ the dual basis:

(a)                          $\langle Y_i, Y^j \rangle = \delta_{ij}.$

Then the Dirac operator on $X \otimes S$ is

(b)                          $D = \sum_i \pi(Y_i) \otimes \gamma(Y^i).$

The construction of a Dirac operator therefore requires a choice of a space of spinors. We saw in Proposition 7.3 that if $\mathfrak{p}_0$ is odd-dimensional, there are two inequivalent choices for such a space. It is natural to ask how changing the choice affects the Dirac operator. Proposition 7.3 shows that the two choices $\gamma$ and $\gamma'$ may be realized on a common space $S$, and that they are related by $\gamma(Y) = -\gamma'(Y)$ (for $Y \in \mathfrak{p}_0$). The two Dirac operators $D$ and $D'$ on $X \otimes S$ are then related by $D = -D'$. We will be interested mostly in the eigenvalues of $D^2$; they are of course unaffected by such a change.

PROPOSITION 7.12 (Parthasarathy [14, Lemma 2.5]; see also [13, Prop. 3.1], and [2, Lemma II.6.11]).

a) *In the setting of Definition 7.11, the operator $D$ is independent of the choice of basis $\{Y_i\}$.*

b) *Write $\Omega_G$ for the Casimir operator of $G$, and $\Omega_K$ for the Casimir operator of $K$ ([8, p. 209). Then*

$$D^2 = -\pi(\Omega_G) \otimes 1 + (\pi \otimes \sigma)(\Omega_K) - (\langle \rho, \rho \rangle - \langle \rho_c, \rho_c \rangle).$$

*Here $\rho_c$ is as defined in (0.4e), and $\rho$ is half the sum of any set of positive roots for $\mathfrak{g}$.*

c) *Suppose that $X$ is a Hermitian representation of $\mathfrak{g}_0$. Then $D$ is formally self-adjoint with respect to the tensor product Hermitian structure on $X \otimes S$:*

$$\langle Dx, y \rangle = \langle x, Dy \rangle \qquad (x, y \in X \otimes S).$$

d) *Suppose $X$ is a $(\mathfrak{g}, K)$-module. Then $X \otimes S$ is a locally finite representation of $\widetilde{K}_0$ for the tensor product action $\pi \otimes \sigma$. The Dirac operator $D$ intertwines this action.*

e) *Suppose $X$ is a unitary $(\mathfrak{g}, K)$-module, and that $\pi(\Omega_G)$ is a scalar operator. Then $D^2$ has nonnegative eigenvalues on $X \otimes S$. In particular, suppose $X$ has an infinitesimal character corresponding by the Harish-Chandra parametrization to a weight $\phi(X) \in \mathfrak{h}^*$. Let $\widetilde{\mu} \in (\mathfrak{t}^c)^*$ be the highest weight of any representation of $\widetilde{K}_0$ on $X \otimes S$. Then*

$$\langle \widetilde{\mu} + \rho_c, \widetilde{\mu} + \rho_c \rangle \geq \langle \phi(X), \phi(X) \rangle.$$

*Proof.* Part (a) is entirely elementary; the argument is omitted in any of the references mentioned. A version of (b) (essentially for $X$ equal to $C^\infty(G)$) is proved in [13]; the observation that Parthasarathy's argument applies in the present setting appears in all the other references. For part (c), the Hermitian property of $\pi$ means that the operators $\pi(Y_i)$ are skew-adjoint. By (7.1e) and (7.3d), the Clifford multiplication operators $\gamma(Y^i)$ are skew-adjoint. By (b) of Definition 7.11, the Dirac operator $D$ is self-adjoint.

For (d), the first assertion is trivial. For the second, suppose $\widetilde{k} \in \widetilde{K}_0$; write $k \in K_0$ for its image by (7.4c). Then

$$(7.5a) \quad \left(\pi(k) \otimes \sigma(\widetilde{k})\right) D \left(\pi(k^{-1}) \otimes \sigma(\widetilde{k}^{-1})\right) = \sum_i \pi(\mathrm{Ad}(k)Y_i) \otimes \gamma(\mathrm{Ad}(k)Y^i).$$

Here we have used for $\pi$ the compatibility of the group and Lie algebra representations, and for $\gamma$ the condition in Definition 7.4. Because $\{\mathrm{Ad}(k)Y_i\}$ is a basis of $\mathfrak{p}_0$ and $\{\mathrm{Ad}(k)Y^i\}$ is the dual basis (by the $K$-invariance of our bilinear form), it follows from (a) that the right side of (7.5a) is $D$, proving (d).

For (e), part (b) shows that $D^2$ is diagonalizable. If $x$ is a nonzero eigenvector of eigenvalue $c$, then it follows from (c) that

$$(7.5b) \qquad\qquad c = \langle Dx, Dx \rangle / \langle x, x \rangle.$$

Because the form on $X$ is positive, the form on $X \otimes S$ is positive as well; so the denominator in (7.5b) is positive and the numerator nonnegative. Therefore $c \geq 0$, as we wished to show. For the last inequality, the eigenvalues of $D^2$ are computed by (b). Under the assumptions in (e), we find first of all that

$$(7.5c) \qquad\qquad \pi(\Omega_G) = \langle \phi(X), \phi(X) \rangle - \langle \rho, \rho \rangle.$$

(This standard formula is proved for example in [7, Ex. 4, p. 134].) Let $\widetilde{\delta}$ be a representation of $\widetilde{K}_0$ of highest weight $\widetilde{\mu}$; then for the same reason we have

$$(7.5d) \qquad\qquad \widetilde{\delta}(\Omega_K) = \langle \widetilde{\mu} + \rho_c, \widetilde{\mu} + \rho_c \rangle - \langle \rho_c, \rho_c \rangle.$$

Suppose now that $x \in X \otimes S$ transforms according to the representation $\widetilde{\delta}$. Then the formulas in (b) and (7.5c) and (7.5d) show that $x$ is an eigenvector for $D^2$. The eigenvalue $c$ is equal to

$$(7.5e)$$
$$-(\langle \phi(X), \phi(X) \rangle - \langle \rho, \rho \rangle) + (\langle \widetilde{\mu} + \rho_c, \widetilde{\mu} + \rho_c \rangle - \langle \rho_c, \rho_c \rangle) - (\langle \rho, \rho \rangle - \langle \rho_c, \rho_c \rangle)$$
$$= \langle \widetilde{\mu} + \rho_c, \widetilde{\mu} + \rho_c \rangle - \langle \phi(X), \phi(X) \rangle.$$

The nonnegativity of $c$ is therefore the last inequality in (e). $\qquad\square$

The inequality in Proposition 7.12(e) is sufficiently important that we state it separately, in a slightly strengthened form.

PROPOSITION 7.13 (Parthasarathy's Dirac operator inequality; see [14, Lemma 2.5]). *Suppose $G$ is a reductive group as in (0.2), and $X$ is an irreducible Hermitian $(\mathfrak{g}, K)$-module. Suppose that $X$ has infinitesimal character corresponding by the Harish-Chandra parametrization to a weight $\phi(X) \in \mathfrak{h}^*$. Recall from Definition 5.6 the canonical real part $\mathrm{RE}\,\phi(X) \in i\mathfrak{t}_0^* + \mathfrak{a}_0^*$. Fix a space $S$ of spinors for $\mathfrak{p}_0$, and let $\widetilde{\mu}$ be the highest weight of a representation $\widetilde{\delta}$ of $\widetilde{K}_0$ occurring in $X \otimes S$.*

*List the representations of $K_0$ occurring in $\widetilde{\delta} \otimes \sigma$ as $\{\delta_1, \ldots, \delta_r\}$. Assume that*

*the Hermitian form on $X$ is positive on the $K_0$-types $\delta_i$.*

(*This is automatic if $X$ is unitary.*) *Then*

$$\langle \widetilde{\mu} + \rho_c, \widetilde{\mu} + \rho_c \rangle \geq \langle \operatorname{RE} \phi(X), \operatorname{RE} \phi(X) \rangle.$$

The improvement over the inequality in Proposition 7.12(e) takes two forms. First, the positivity of the Hermitian form on $X$ is assumed only on a specific finite set of $K_0$-types, and not on the entire representation. Second, we have replaced $\phi(X)$ by its canonical real part. To see that this is an improvement, write $\phi(X) = \operatorname{RE} \phi(X) + i \operatorname{IM} \phi(X)$. Then the squared length of $\phi(X)$ is

(7.6a)   $\langle \operatorname{RE} \phi(X), \operatorname{RE} \phi(X) \rangle - \langle \operatorname{IM} \phi(X), \operatorname{IM} \phi(X) \rangle + 2i \langle \operatorname{RE} \phi(X), \operatorname{IM} \phi(X) \rangle.$

The inner product is real and positive definite on the image of RE and IM; so the first term here is positive, the second negative, and the third purely imaginary. On a Hermitian representation the Casimir operator must have a real eigenvalue, so the third term is zero; and we get

(7.6b)   $$\langle \phi(X), \phi(X) \rangle \leq \langle \operatorname{RE} \phi(X), \operatorname{RE} \phi(X) \rangle.$$

This says that the inequality in Proposition 7.13 implies the one in Proposition 7.12(e).

*Proof.* Write $X_0$ for the $K_0$-invariant subspace of $X$ generated by the isotypic subspaces of type $\delta_i$. Suppose $x \in X \otimes S$ is a nonzero vector transforming according to the representation $\widetilde{\delta}$. Then it follows from the definition of $\{\delta_i\}$ that $x \in X_0 \otimes S$. By Proposition 7.12(d), $Dx$ also transforms according to $\widetilde{\delta}$, so $Dx \in X_0 \otimes S$. By hypothesis, our Hermitian form is positive definite on $X_0 \otimes S$; so

(7.7a)   $$\langle Dx, Dx \rangle / \langle x, x \rangle \geq 0.$$

Now the same calculation as in (7.5b) and (7.5e) shows that

(7.7b)   $$\langle \widetilde{\mu} + \rho_c, \widetilde{\mu} + \rho_c \rangle \geq \langle \phi(X), \phi(X) \rangle.$$

To get the canonical real part into the picture, we use Theorem 16.10 in [8]. The result is of sufficient interest to state separately.

PROPOSITION 7.14 ([8, Th. 1 6.10]).   *Suppose $G$ is a real reductive group in Harish-Chandra's class, and $(X, \langle \ , \ \rangle)$ is an irreducible Hermitian $(\mathfrak{g}, K)$-module. Suppose $H$ is a $\theta$-stable Cartan subgroup of $G$, and that the infinitesimal character of $X$ is given in the Harish-Chandra parametrization by a weight* $\phi(X) \in \mathfrak{h}^*$. Then there exist a parabolic subgroup $P_1 = M_1 A_1 N_1$ of $G$, an irreducible Hermitian Harish-Chandra module $X_1$ for $M_1$, and a unitary character $\nu_1 \in i\mathfrak{a}_{1,0}^*$, with the following three properties.

a) The infinitesimal character of $X_1$ has a Harish-Chandra parameter that is conjugate by $\mathrm{Ad}(\mathfrak{g})$ to $\mathrm{RE}\,\phi(X)$.

b) The weight $\nu_1$ is conjugate by $\mathrm{Ad}(\mathfrak{g})$ to $i\mathrm{IM}\,\phi(X)$.

c) There is an isomorphism of Hermitian representations

$$X \simeq \mathrm{Ind}_{M_1 A_1 N_1}^{G}(X_1 \otimes e^{\nu_1} \otimes 1).$$

Consider the Hermitian representation

$$X' = \mathrm{Ind}_{M_1 A_1 N_1}^{G}(X_1 \otimes 1 \otimes 1);$$

this is a $(\mathfrak{g}, K)$-module of finite length.

d) The representations $X$ and $X'$ are isomorphic as Hermitian representations of $K$. In particular, their signature characters coincide (Definition 5.3).

e) The infinitesimal character of $X'$ is just that of $X_1 \otimes 1$ composed with a Harish-Chandra homomorphism, and is therefore given in the Harish-Chandra parametrization by the weight

$$\phi(X') = \mathrm{RE}\,\phi(X).$$

*Proof.* Parts (a)–(c) are essentially proved in [8, Th. 16.10]. Part (d) follows immediately from the construction of the Hermitian form on a parabolically induced representation in the "compact picture." Part (e) is a standard calculation of infinitesimal characters of induced representations ([8, Prop. 8.22]). □

Returning to the proof of Proposition 7.12, we now apply Proposition 7.12 to the representation $X$. By (d), the representation $X'$ inherits the positivity assumption made on the $K$-types of $X$. We may also identify the element $x$ with an element of $X' \otimes S$, still transforming by $\tilde{\delta}$. Applying the inequality (7.7b) to the representation $X'$ therefore gives the conclusion we want.

Now we will relate these ideas to the notion of unitarily small $K$-types. Here is a classical fact that hints at what we want.

PROPOSITION 7.15. *Suppose $G$ is a connected reductive group as in (0.2), and $(\sigma, S)$ is a spin representation for $\mathfrak{p}_0$ (a representation of $\widetilde{K}_0$; see (7.4)). Suppose $\delta$ is a finite-dimensional irreducible representation of $K$. Then the following conditions on $\delta$ are equivalent.*

a) *The representation $\delta$ occurs in $\bigwedge \mathfrak{p}$.*

b) *There is an irreducible constituent $\sigma_1$ of $\sigma$ so that $\delta \otimes \sigma_1$ contains a $\widetilde{K}_0$-type in $S$.*

c) *There are irreducible constituents $\sigma_1$ and $\sigma_2$ of $\sigma$ so that $\delta$ occurs in $\sigma_1 \otimes \sigma_2$.*

This is an easy consequence of Corollary 7.10; we omit the argument. By Proposition 7.1, a representation $\delta$ satisfying these conditions is unitarily small. What we want to do is weaken the conditions until they are precisely equivalent to unitarily small. We begin with some properties of the spin representation.

LEMMA 7.16 ([2, Lemma II.6.9]). *Suppose $G$ is a connected reductive group as in (0.2), and $(\sigma, S)$ is a spin representation for $\mathfrak{p}_0$ (a representation of $\widetilde{K}_0$; see (7.4)).*

a) *The representation $\sigma$ is of type $\frac{1}{2}\Delta(\mathfrak{p}, \mathfrak{t}^c)$ (Lemma 6.11). That is, every weight of $\sigma$ is a sum of noncompact roots with coefficients between $-1/2$ and $1/2$.*

b) *Suppose $\Delta^+(\mathfrak{g}, \mathfrak{t}^c)$ is a system of positive roots containing $\Delta^+(\mathfrak{k}, \mathfrak{t}^c)$. Then*

$$\rho_n = 1/2 \sum_{\beta \in \Delta^+(\mathfrak{p}, \mathfrak{t}^c)} \beta$$

*is a highest weight of $\sigma$. Similarly $-\rho_n$ is a lowest weight of $\sigma$.*

c) *Every highest weight of $\sigma$ is of the form in (b).*

PROPOSITION 7.17. *Suppose $G$ is a semisimple group as in (0.2), and $(\sigma, S)$ is a spin representation for $\mathfrak{p}_0$. Suppose $\delta$ is an irreducible representation of $K$. In the notation of Definition 6.1, the following conditions on $\delta$ are equivalent.*

a) *The representation $\delta$ is unitarily small.*

b) *There is an irreducible constituent $\sigma_1$ of $\sigma$ and an irreducible constituent $\sigma_2$ of $\delta \otimes \sigma_1$ so that $\sigma_2$ is of type $\frac{1}{2}\Delta(\mathfrak{p}, \mathfrak{t}_c)$ (Lemma 6.11). That is, every weight of $\sigma_2$ is a sum of noncompact roots with coefficients between $-1/2$ and $1/2$.*

c) *There is an irreducible constituent $\sigma_1$ of $\sigma$ and an irreducible representation $\sigma_2$ of type $\frac{1}{2}\Delta(\mathfrak{p}, \mathfrak{t}_c)$, with the property that $\delta|_{K_0}$ contains an irreducible constituent of $\sigma_1 \otimes \sigma_2$.*

d) *There are irreducible representations $\sigma_1$ and $\sigma_2$ of $\widetilde{K}_0$, of type $\frac{1}{2}\Delta(\mathfrak{p}, \mathfrak{t}_c)$, with the property that $\delta|_{K_0}$ contains an irreducible constituent of $\sigma_1 \otimes \sigma_2$.*

*Proof.* We will prove that (a)$\Rightarrow$(b)$\Rightarrow$(c)$\Rightarrow$(d)$\Rightarrow$(a). So suppose (a) holds. Let $\mu \in i(\mathfrak{t}_0^c)^*$ be a highest weight of $\delta$. By Theorem 6.7(f), we can find a positive root system $\Delta^+(\mathfrak{g}, \mathfrak{t}^c)$ so that

$$\mu = \sum_{\beta \in \Delta^+(\mathfrak{p}, \mathfrak{t}^c)} c_\beta \beta \qquad (0 \le c_\beta \le 1).$$

Write $\delta_0$ for an irreducible constituent of $\delta$ (as a representation of $K_0$) of highest weight $\mu$, and $\sigma_1$ for an irreducible constituent of $\sigma$ of lowest weight $-\rho_n$ (Lemma 7.16(b)). By a theorem of Parthasarathy, Rao, and Varadarajan ([12, Cor. 1 to Th. 2.1]), there is an irreducible constituent $\sigma_2$ of $\delta_0 \otimes \sigma_1$ of extremal weight

$$\mu - \rho_n = \sum_{\beta \in \Delta^+(\mathfrak{p}, \mathfrak{t}^c)} (c_\beta - 1/2)\beta \qquad (-1/2 \leq c_\beta - 1/2 \leq 1/2).$$

By Lemma 6.11, $\sigma_2$ is of type $\frac{1}{2}\Delta(\mathfrak{p}, \mathfrak{t}_c)$. This is (b). The implication (b)$\Rightarrow$(c) is formal from the self-duality of $S$. That (c)$\Rightarrow$(d) is a formal consequence of Lemma 7.16(a). Finally, assume (d). The weights of $\sigma_1 \otimes \sigma_2$ are all of the form $\gamma_1 + \gamma_2$, with $\gamma_i$ a weight of $\sigma_i$. It follows immediately that $\sigma_1 \otimes \sigma_2$ is of type $\Delta(\mathfrak{p}, \mathfrak{t}_c)$. Now (a) follows from Theorem 6.7(g).          $\square$

We can now prove a result in the direction of Conjecture 5.7; more precisely, of Conjecture 5.7′.

PROPOSITION 7.18.     *Suppose $G$ is a semisimple group, and suppose $X \in \Pi_h^0(G)$ is an irreducible Hermitian $(\mathfrak{g}, K)$-module (Definition 5.1). Let $\mathfrak{h}$ be a $\theta$-stable Cartan subalgebra of $\mathfrak{g}$, and $\phi \in \mathfrak{h}^*$ a weight parametrizing the infinitesimal character of $X$. Assume that*

$$\langle \mathrm{RE}\,\phi, \mathrm{RE}\,\phi \rangle > \langle \rho, \rho \rangle.$$

*Then the Hermitian form on $X$ must be indefinite on $K$-types in $B_u^0(G)$.*

*Proof.* Suppose not. Then we may assume that the form is positive definite on the $K$-types in $B_u^0(G)$. By hypothesis $X$ contains a unitarily small $K$-type $\delta \in B_u^0(G)$. By Proposition 7.17, we can find a representation $\sigma_2$ of type $\frac{1}{2}\Delta(\mathfrak{p}, \mathfrak{t}^c)$ occurring in $\delta \otimes \sigma$, and therefore in $X \otimes S$. We are going to apply Proposition 7.12. We therefore list the representations of $K$ containing constituents of $\sigma_2 \otimes \sigma$ as $\{\delta_1, \ldots, \delta_r\}$. By Proposition 7.17(d), these are all unitarily small; so by hypothesis the Hermitian form on $X$ is positive on them. Write $\mu_2$ for the highest weight of $\sigma_2$. According to Proposition 7.12,

$$(7.8a) \qquad\qquad \langle \mu_2 + \rho_c, \mu_2 + \rho_c \rangle \geq \langle \mathrm{RE}\,\phi, \mathrm{RE}\,\phi \rangle.$$

Fix a positive root system $\Delta^+(\mathfrak{g}, \mathfrak{t}^c)$ containing $\Delta^+(\mathfrak{k}, \mathfrak{t}^c)$, and write $\rho_n$ for half the sum of the corresponding noncompact positive roots. Because $\sigma_2$ is of type

$\frac{1}{2}\Delta(\mathfrak{p}, \mathfrak{t}^c)$, we can write

$$(7.8b) \qquad \mu_2 = \sum_{\beta \in \Delta^+(\mathfrak{p}, \mathfrak{t}^c)} (c_\beta - c_{-\beta})\beta \qquad (0 \le c_{\pm\beta} \le 1/2)$$

$$= \sum_{\beta \in \Delta^+(\mathfrak{p}, \mathfrak{t}^c)} a_\beta\beta \qquad (-1/2 \le a_\beta \le 1/2)$$

$$= \rho_n - \sum_{\beta \in \Delta^+(\mathfrak{p}, \mathfrak{t}^c)} b_\beta\beta \qquad (0 \le b_\beta \le 1).$$

Adding $\rho_c$ gives

$$(7.8c) \qquad \mu_2 + \rho_c = \rho - \sum_{\beta \in \Delta^+(\mathfrak{p}, \mathfrak{t}^c)} b_\beta\beta \qquad (0 \le b_\beta \le 1).$$

By Proposition 1.10, it follows that $\mu_2 + \rho_c$ belongs to the convex hull of the Weyl group orbit of $\rho$. Consequently

$$(7.8d) \qquad \langle \mu_2 + \rho_c, \mu_2 + \rho_c \rangle \le \langle \rho, \rho \rangle.$$

Combining this with (7.8a) gives $\langle \rho, \rho \rangle \ge \langle \operatorname{RE}\phi, \operatorname{RE}\phi \rangle$, contradicting the hypothesis in the proposition. This contradiction completes the proof. $\qquad\square$

This argument seems to us to be rather natural and powerful; so it is disappointing that we cannot deduce Conjecture 5.7′ from it. Perhaps the problem lies in Proposition 7.13. The following sharpening of that proposition would (by the argument above) immediately imply Conjecture 5.7′.

*Conjecture* 7.13. Suppose we are in the setting of Proposition 7.13; make the same hypotheses as for that proposition. Then $\operatorname{RE}\phi(X)$ is conjugate by $\operatorname{Ad}(\mathfrak{g})$ to a weight in the convex hull of the Weyl group orbit of $\widetilde{\mu} + \rho_c$.

We are being a little bit vague about the Weyl group in this conjecture; probably it is enough to use the Weyl group $W(\mathfrak{g}, \mathfrak{t}^c)$ of the restricted root system, but the conjecture would still imply Conjecture 5.7′ even if we take the full Weyl group $W(\mathfrak{g}, \mathfrak{h}^c)$.

## 8. Topology of the unitary dual

In this section we consider the relationship between the parametrization of the unitary dual given (conjecturally) by Theorem 5.8 and the Fell topology. A careful treatment of this topic would require a careful foundational discussion of the Fell topology, which we prefer to avoid. Our arguments will therefore be sketchy at several crucial points.

We begin with a discussion of infinitesimal characters. Fix a $\theta$-stable positive root system $\Delta^+(\mathfrak{g}, \mathfrak{h}^c)$ for the Cartan subalgebra $\mathfrak{h}^c$ of (0.4a); this choice is equivalent to the choice of a system $\Delta^+(\mathfrak{g}, \mathfrak{t}^c)$ of positive restricted roots. We may assume that $\Delta^+(\mathfrak{g}, \mathfrak{t}^c)$ contains our fixed system $\Delta^+(\mathfrak{k}, \mathfrak{t}^c)$ of positive roots for $\mathfrak{k}$.

*Definition* 8.1.   With notation as above, suppose $\phi \in (\mathfrak{h}^c)^*$ is a weight. We say that $\phi$ is *complex dominant* if it satisfies either of the following equivalent conditions.

a) The canonical real part $\mathrm{RE}\,\phi$ (Definition 5.6) is dominant for $\Delta^+(\mathfrak{g}, \mathfrak{h}^c)$; and the canonical imaginary part $\mathrm{IM}\,\phi$ is dominant for those positive roots orthogonal to $\mathrm{RE}\,\phi$.

b) For every positive root $\alpha \in \Delta^+(\mathfrak{g}, \mathfrak{h}^c)$, we have either $\mathrm{Re}\,\langle \phi, \alpha \rangle > 0$; or $\mathrm{Re}\,\langle \phi, \alpha \rangle = 0$, and $\mathrm{Im}\,\langle \phi, \alpha \rangle > 0$; or $\langle \phi, \alpha \rangle = 0$.

The set of complex dominant weights is denoted

$$C_{\mathbb{C}}(\Delta^+(\mathfrak{g}, \mathfrak{h}^c)) \subset (\mathfrak{h}^c)^*.$$

It is a fundamental domain for the action of $W(\mathfrak{g}, \mathfrak{h}^c)$ on $(\mathfrak{h}^c)^*$:

$$C_{\mathbb{C}}(\Delta^+(\mathfrak{g}, \mathfrak{h}^c)) \simeq (\mathfrak{h}^c)^*/W(\mathfrak{g}, \mathfrak{h}^c).$$

We endow $C_{\mathbb{C}}(\Delta^+(\mathfrak{g}, \mathfrak{h}^c))$ with the quotient topology from this bijection (and *not* with the subspace topology from $\mathfrak{h}^c$, which has more neighborhoods of purely imaginary elements).

The ordinary positive Weyl chamber $C_{\mathbb{R}}(\Delta^+(\mathfrak{g}, \mathfrak{h}^c))$ consisting of dominant elements of $(i\mathfrak{t}_0^c + \mathfrak{a}_0^c)^*$ is a subset of $C_{\mathbb{C}}(\Delta^+(\mathfrak{g}, \mathfrak{h}^c))$. We will use the fact that the topology on $C_{\mathbb{R}}$ coming from this inclusion (that, is the quotient topology from

$$C_{\mathbb{R}}(\Delta^+(\mathfrak{g}, \mathfrak{h}^c)) \simeq (i\mathfrak{t}_0^c + \mathfrak{a}_0^c)^*/W(\mathfrak{g}, \mathfrak{h}^c)$$

coincides with the subspace topology from the inclusion

$$C_{\mathbb{R}}(\Delta^+(\mathfrak{g}, \mathfrak{h}^c)) \subset (i\mathfrak{t}_0^c + \mathfrak{a}_0^c)^*.$$

We now write

(8.1a)                    $$\mathfrak{Z}(\mathfrak{g}) = \text{center of } U(\mathfrak{g}).$$

The Harish-Chandra isomorphism is

(8.1b)                    $$\mathfrak{Z}(\mathfrak{g}) \simeq S(\mathfrak{h}^c)^{W(\mathfrak{g}, \mathfrak{h}^c)}.$$

An *infinitesimal character* is by definition a homomorphism from $\mathfrak{Z}(\mathfrak{g})$ to the complex numbers, or equivalently a maximal ideal in $\mathfrak{Z}(\mathfrak{g})$. By analogy with the notation we are using for representations, the set of infinitesimal characters is written

(8.1c)                    $$\Pi(\mathfrak{Z}(\mathfrak{g})) = \mathrm{Max}\,\mathfrak{Z}(\mathfrak{g}).$$

The isomorphism (8.1b) identifies these infinitesimal characters with $W$ orbits on $\mathfrak{h}^*$:

$$(8.1d) \qquad \Pi(\mathfrak{Z}(\mathfrak{g})) \simeq (\mathfrak{h}^c)^*/W(\mathfrak{g}, \mathfrak{h}^c) \simeq C_{\mathbb{C}}(\Delta^+(\mathfrak{g}, \mathfrak{h}^c));$$

the last equivalence comes from Definition 8.1. We will describe this map a little more carefully in (8.2) below.

The set of infinitesimal characters carries a "Fell topology" of its own: the weakest topology making each of the evaluation maps

$$(8.2a) \qquad e_z \colon \Pi(\mathfrak{Z}(\mathfrak{g})) \to \mathbb{C}, \qquad e_z(\chi) = \chi(z) \qquad (z \in \mathfrak{Z}(\mathfrak{g}))$$

continuous. (We use the usual analytic topology on $\mathbb{C}$; if instead we used the Zariski topology, we would get the restriction to Max $\mathfrak{Z}(\mathfrak{g})$ of the Zariski topology on Spec $\mathfrak{Z}(\mathfrak{g})$.) We can restate the definition of the topology as follows. A set $U \subset \Pi(\mathfrak{Z}(\mathfrak{g}))$ is a neighborhood of the element $\chi_0 \in U$ if and only if there are elements $z_1, \ldots, z_n \in \mathfrak{Z}(\mathfrak{g})$ and positive numbers $\varepsilon_1, \ldots, \varepsilon_n$ with the property that $U$ contains

$$(8.2b) \quad N(z, \varepsilon; \chi_0) = \{\chi \in \Pi(\mathfrak{Z}(\mathfrak{g})) \mid |\chi(z_i) - \chi_0(z_i)| < \varepsilon_i \quad (i = 1, \ldots, n)\}.$$

Now Chevalley's theorem says that $\mathfrak{Z}(\mathfrak{g})$ is isomorphic to a polynomial ring $\mathbb{C}[x_1, \ldots, x_m]$. (The elements $x_1, \ldots, x_m$ may be thought of as generators of $\mathfrak{Z}(\mathfrak{g})$; by the Harish-Chandra isomorphism (8.1b), they correspond to generators for the algebra of $W$-invariants in $S(\mathfrak{h}^c)$.) Each element $\xi = (\xi_1, \ldots, \xi_m) \in \mathbb{C}^m$ defines an element $\chi_\xi \in \Pi(\mathfrak{Z}(\mathfrak{g}))$; the characteristic property is $\chi_\xi(x_i) = \xi_i$. By the Nullstellensatz, this correspondence gives an isomorphism

$$(8.2c) \qquad \Pi(\mathfrak{Z}(\mathfrak{g})) \simeq \mathbb{C}^m.$$

It is very easy to check (by the description of neighborhoods in (8.2b)) that this isomorphism sends the topology on $\Pi(\mathfrak{Z}(\mathfrak{g}))$ to the standard topology on $\mathbb{C}^m$. We want to know something slightly more subtle, however. We may identify $S(\mathfrak{h}^c)$ with the algebra of polynomial functions on $(\mathfrak{h}^c)^*$. Each $\lambda \in (\mathfrak{h}^c)^*$ therefore defines an algebra homomorphism

$$(8.2d) \qquad \xi_\lambda \colon S(\mathfrak{h}^c) \to \mathbb{C}, \qquad p \mapsto p(\lambda).$$

We use the same symbol to denote the composition with the Harish-Chandra homomorphism:

$$(8.2e) \qquad \xi_\lambda \colon \mathfrak{Z}(\mathfrak{g}) \to \mathbb{C}.$$

In this way we get a map

$$(8.2f) \qquad (\mathfrak{h}^c)^* \to \Pi(\mathfrak{Z}(\mathfrak{g})), \qquad \lambda \mapsto \xi_\lambda.$$

This is the map of (8.1d): it is surjective, and the fibers are precisely the orbits of $W(\mathfrak{g}, \mathfrak{h}^c)$.

It is more or less obvious from the definitions that the map of (8.2f) is continuous (from the standard topology on $(\mathfrak{h}^c)^*$ to the "Fell topology" on $\Pi(\mathfrak{Z}(\mathfrak{g}))$). What we need is the stronger assertion that the topology on $\Pi(\mathfrak{Z}(\mathfrak{g}))$ is actually the quotient topology; equivalently, that the map (8.2f) is open. Here is an explicit formulation.

LEMMA 8.2.    *In the setting of* (8.2), *suppose* $(p_1, \ldots, p_m)$ *is a collection of generators for* $S(\mathfrak{h}^c)^{W(\mathfrak{g}, \mathfrak{h}^c)}$. *Fix* $\lambda_0 \in (\mathfrak{h}^c)^*$ *and* $\varepsilon > 0$. *Then there is a positive number* $\delta$ *so that if* $\lambda \in (\mathfrak{h}^c)^*$ *and*

$$|p_i(\lambda) - p_i(\lambda_0)| < \delta \qquad (i = 1, \ldots, m),$$

*then there is a* $w \in W(\mathfrak{g}, \mathfrak{h}^c)$ *such that*

$$|\lambda - w \cdot \lambda_0| < \varepsilon.$$

In the case of $\mathrm{GL}(n)$, this lemma in turn can be reformulated as follows. Suppose $f_0$ is a monic polynomial of degree $n$ with roots $(\lambda_0^1, \ldots, \lambda_0^n)$, and $\varepsilon > 0$. Then there is a $\delta > 0$ so that if $f$ is a monic polynomial whose coefficients are within $\delta$ of those of $f_0$, then the roots $(\lambda^1, \ldots, \lambda^n)$ of $f$ are each (after permutation) within $\varepsilon$ of the corresponding roots of $f_0$. That is, it says that the collection of roots of a polynomial depends continuously on the coefficients. This is of course well known. In keeping with the incomplete nature of this section, we omit the argument for the general case.

COROLLARY 8.3.    *The bijection of* (8.1d) *is a homeomorphism from the "Fell topology" to the quotient topology from* $(\mathfrak{h}^c)^*$.

*Definition* 8.3.    The *infinitesimal character map*

$$\Pi_u(G) \to \Pi(\mathfrak{Z}(\mathfrak{g})), \qquad \pi \mapsto \xi(\pi)$$

sends an irreducible unitary representation $\pi$ to the homomorphism $\xi(\pi)$ by which $\mathfrak{Z}(\mathfrak{g})$ acts on the space of smooth vectors of $\pi$.

THEOREM 8.5 (Bernat-Dixmier [1]).    *The infinitesimal character map is continuous from the Fell topology on the unitary dual to the "Fell topology" on* $\Pi(\mathfrak{Z}(\mathfrak{g}))$ (cf. (8.2)).

COROLLARY 8.6.    *The map*

$$\Pi_u(G) \to C_{\mathbb{C}}(\Delta^+), \qquad \pi \mapsto \phi(\pi; \Delta^+)$$

(*obtained by composing the infinitesimal character map of Theorem* 8.5 *with the isomorphisms of* (8.1d)) *is continuous from the Fell topology on the unitary dual to the quotient topology* (Definition 8.1) *on the cone of complex dominant weights.*

This is immediate from Theorem 8.5 and Corollary 8.3.

We are now in a position to draw some conclusions about the topology of the unitary dual. A little more notation is helpful. In the setting of Definition 8.1, suppose $\lambda \in i(\mathfrak{t}_0^c)^*$. Define $D_{\Delta^+}(\lambda)$ to be the unique $W(\mathfrak{g}, \mathfrak{t}^c)$ conjugate of $\lambda$ belonging to the positive Weyl chamber $C$ of (0.5c). (Notice that $C$ is contained in $C_{\mathbb{R}}(\Delta^+)$.) The map

$$(8.3a) \qquad i(\mathfrak{t}_0^c)^* \to C, \qquad \lambda \mapsto D_{\Delta^+}(\lambda)$$

is continuous and finite-to-one. On the other hand, the canonical real part (Definition 5.6) defines a continuous map

$$(8.3b) \qquad (\mathfrak{h}^c)^* \to i(\mathfrak{t}_0^c)^* + (\mathfrak{a}_0^c)^*, \qquad \phi \mapsto \mathrm{RE}\,\phi.$$

This map is $W$-equivariant, and therefore descends to a continuous map for the quotient topologies

$$(8.3c) \qquad C_{\mathbb{C}}(\Delta^+) \to C_{\mathbb{R}}(\Delta^+), \qquad \phi \mapsto \mathrm{RE}\,\phi.$$

Similarly, the map $T_\rho$ of Proposition 1.4 is continuous:

$$(8.3d) \qquad i(\mathfrak{t}_0^c)^* + (\mathfrak{a}_0^c)^* \to i(\mathfrak{t}_0^c)^* + (\mathfrak{a}_0^c)^*, \qquad \phi \mapsto T_\rho(\phi).$$

To see this, recall that $T_\rho$ is defined separately on each closed Weyl chamber $C'$: if $\rho'$ is half the sum of the positive roots for $C'$, and $P'$ is the projection on $C'$ (Definition 1.2), then

$$(8.3e) \qquad T_\rho(\phi) = P'(\phi - \rho').$$

This mapping is obviously continuous on $C'$, because $P'$ is a contraction. Since $C'$ is closed, and the definitions agree on the intersections of different chambers (Proposition 1.4), $T_\rho$ is continous.

THEOREM 8.6. *In the setting of Definition 8.1, suppose $\pi$ is an irreducible unitary representation, and*

$$\phi(\pi; \Delta^+) \in C_{\mathbb{C}}(\Delta^+) \subset (\mathfrak{h}^c)^*$$

*comes from the infinitesimal character as in Corollary 8.6. Define*

$$\lambda(\pi; \Delta^+) = T_\rho(\mathrm{RE}\,\phi(\pi; \Delta^+))$$

(*notation as in* Proposition 1.4).

a) *The weight $\lambda(\pi; \Delta^+)$ belongs to the positive Weyl chamber $C_{\mathbb{R}}(\Delta^+) \subset i(\mathfrak{t}_0^c)^* + (\mathfrak{a}_0^c)^*$.*

b) *The map*

$$\Pi_u(G) \to C_{\mathbb{R}}(\Delta^+), \qquad \pi \mapsto \lambda(\pi; \Delta^+)$$

*is continuous from the Fell topology to the standard topology on the cone of dominant weights.*

c) *Assume that Conjecture 5.7 holds, and that* $\pi \in \Pi_u^{\lambda_u}(G)$ *(Definition 0.2). Then*

$$D_{\Delta^+}(\lambda_u) = \lambda(\pi; \Delta^+)$$

*(notation as in* (8.3a)).

*Proof.* Because $\mathrm{RE}\,\phi(\pi; \Delta^+)$ is dominant for $\Delta^+(\mathfrak{g}, \mathfrak{h}^c)$, $T_\rho$ may be computed by subtracting $\rho$ and projecting on the positive Weyl chamber. This proves (a). Part (b) follows from Corollary 8.6, the continuity of the maps in (8.3c and d), and the identification noted in Definition 8.1 of the quotient and subspace topologies on $C_{\mathbb{R}}(\Delta^+)$. Part (c) is a reformulation of Theorem 5.8(c). □

COROLLARY 8.8.    *Fix* $\lambda_u \in \Lambda_u$ (cf. (0.5e)). *Choose representatives* $\lambda_u^1, \ldots, \lambda_u^r$ *for the orbits of* $R(G)$ *on*

$$W(\mathfrak{g}, \mathfrak{t}^c) \cdot \lambda_u \cap \Lambda_u.$$

*Assume that Conjecture 5.7 holds.*

a) *The set* $\bigcup_{i=1}^r \Pi_u^{\lambda_u^i}(G)$ *of unitary representations is open and closed in the Fell topology.*

b) *Suppose* $i \neq j$, *and* $\pi \in \Pi_u^{\lambda_u^i}(G)$. *Then* $\pi$ *contains no $K$-types in* $B_u^{\lambda_u^j}(G)$.

c) *Each set* $\Pi_u^{\lambda_u^i}(G)$ *is open and closed in the Fell topology.*

d) *The mapping*

$$\Pi_u(G) \to \Lambda_u/R(G), \qquad \pi \mapsto \lambda_u(\pi)$$

*(notation as in* Definition 0.2) *is continuous from the Fell topology to the natural (discrete) topology on the image.*

The statement of this corollary is a little cumbersome; part (c), for example, is obviously stronger than (a). The formulation is chosen to make the proof clearer.

One might expect that (d) should be more or less a formal consequence of the definitions. To see why it is not, consider the analogous map

(8.4)                    $$\Pi_u(G) \to \Lambda_a/R(G), \qquad \pi \mapsto \lambda_a(\pi)$$

(notation  as in (2.1) and Definition 2.5). The definition of $\lambda_a$ is formally very similar to that of $\lambda_u$; but the mapping in (8.4) is *not* continuous in the Fell topology. To see why, suppose $G = \mathrm{SL}(2, \mathbb{R})$, and let $\pi$ be the first holomorphic discrete series representation. In the notation of Example 6.2, the lowest $K$-type of $\pi$ is $\mu_2$, and $\lambda_a(\pi) = 1$. But there is a sequence $\{\pi_j\}$ of spherical complementary series representations having $\pi$ as a limit point. (The reason is that $\pi$ appears as a composition factor in the reducible spherical principal series representation at the end of the complementary series.) Each $\pi_j$ has lowest $K$-type $\mu_0$, so $\lambda_a(\pi_j) = 0$. So $\lambda_a(\pi_j)$ fails to converge to $\lambda_a(\pi)$. (On the other hand, Example 6.2 shows that $\lambda_u(\pi) = 0 = \lambda_u(\pi_j)$, so there is no contradiction to (d) in the corollary.)

*Proof.* By Theorem 8.6(c), we have a continuous map

$$\Pi_u(G) \to D_{\Delta^+}(\Lambda_u), \qquad \pi \mapsto \lambda(\pi, \Delta^+) = D_{\Delta^+}(\lambda_u(\pi)).$$

The range is discrete because $\Lambda_u$ is (cf. (0.5e)), and so carries the discrete topology. The set of unitary representations in (a) is just the preimage of $D_{\Delta^+}(\lambda_u)$, and so must be open and closed. Part (b) is immediate from Lemma 5.5.

For (c), it suffices to prove that each $\Pi_u^{\lambda_i}(G)$ is closed. So suppose $\{\pi_n\}$ is a sequence of representations in $\Pi_u^{\lambda_i}(G)$, and that $\pi$ is a limit point. By (a), $\pi$ must belong to some $\Pi_u^{\lambda_j}(G)$; we want to show that $i = j$. Suppose not. By Definition 0.2, $\pi$ must contain a $K$-type $\tau$ in $B_u^{\lambda_j}(G)$. We can therefore form a corresponding nonzero matrix coefficient $f_\pi$; it transforms on the right and left by $\tau$ under $K$. By the definition of the Fell topology, $f_\pi$ may be uniformly approximated on compact sets by matrix coefficients of the representations $\pi_n$. By (b), these representations do not contain $\tau$; and it follows easily that their matrix coefficients cannot approximate $f_\pi$. This contradiction proves (c). Part (d) is immediate from (c). □

COROLLARY 8.9. *In the setting of Theorem 5.8, assume that Conjecture 5.7 holds. Then the bijection*

$$\mathcal{L}_{S_u}(\lambda_u) \colon \Pi_u^{\lambda_u}(G(\lambda_u)) \to \Pi_u^{\lambda_u}(G)$$

*is a homeomorphism in the Fell topology; the image is open and closed in $\Pi_u(G)$.*

*Sketch of proof.* According to fundamental (and partly unpublished) work of Miličić, the Fell topology may be described in terms of the reducibility of certain parabolically induced representations. To be a little more precise, suppose $P = MAN$ is a parabolic subgroup of $G$, $\delta \in \Pi_a(M)$ is an admissible representation, and

(8.5a) $$\nu \colon [0, 1] \to \mathfrak{a}^*$$

is a continuous map. Define a family of admissible representations by

(8.5b) $$I(t) = \mathrm{Ind}_P^G(\delta \otimes \nu(t) \otimes 1), \qquad (t \in [0, 1]).$$

Now we make the following assumption:

(8.5c)  For all $t \in [0, 1)$, the representation $I(t)$ is irreducible and unitary.

Here by "unitary" we mean "infinitesimally equivalent to a unitary representation." List the irreducible composition factors of the representation $I(1)$ as $I_1, I_2, \ldots, I_r$. As a consequence of assumption 8.5(c) and a theorem of Miličić, (see [11]), all of the representations $I_j$ are unitary; and they are precisely the limit points of $I(t)$ as $t$ approaches 1:

$$(8.5d) \qquad \lim_{t \to 1} I(t) = \{I_1, \ldots, I_r\}.$$

Finally, every convergent sequence in the Fell topology has a subsequence constructed in this way.

By this description of the Fell topology, and the good behavior of $\mathcal{L}_{S_u}(\lambda_u)$ (in particular Theorem 5.4(e)), it is not difficult to prove the corollary. We omit the details.

We make one final remark. All of the interesting results in this section depend on Conjecture 5.7; it is natural to ask what can be said absolutely. In the setting of Theorems 5.4 and 5.8, we can consider a certain subset

$$(8.6a) \qquad \Pi_{u,\text{good}}^{\lambda_u}(G(\lambda_u)) \subset \Pi_u^{\lambda_u}(G(\lambda_u))$$

consisting of unitary representations for which the infinitesimal character $\phi(Z)$ satisfies

$$(8.6b) \qquad \text{RE } \phi(Z) \in \lambda_u + \{\text{convex hull of } W \cdot \rho\};$$

here both the Weyl group and $\rho$ come from $G(\lambda_u)$. The content of Conjecture 5.7 is that this set is all of $\Pi_u^{\lambda_u}(G(\lambda_u))$; but in any case it is a well-defined closed subset. The proof of Theorem 5.8 shows that $\mathcal{L}_{S_u}(\lambda_u)$ defines a bijection from this set to a closed subset

$$(8.6c) \qquad \Pi_{u,\text{good}}^{\lambda_u}(G) \subset \Pi_u^{\lambda_u}(G).$$

We define

$$(8.6d) \qquad \Pi_{u,\text{good}} = \cup_{\lambda_u \in \Lambda_u} \Pi_{u,\text{good}}^{\lambda_u}(G),$$

a closed subset of $\Pi_u(G)$. Everything we have said about $\Pi_u(G)$ using Conjecture 5.7 applies unconditionally to $\Pi_{u,\text{good}}(G)$.

Mathematical Sciences, New Mexico State University, Las Cruces, NM
*E-mail address*: ssalaman@nmsu.edu
Massachusetts Institute of Technology, Cambridge, MA
*E-mail address*: dav@math.mit.edu